\documentclass[11p,reqno]{amsart}
\usepackage[toc,title,page]{appendix}
\usepackage{etoolbox}

\patchcmd{\abstract}{\null\vfil}{}{}{}

\usepackage[dvipsnames]{xcolor}
\usepackage{tikz}
\usepackage{amsmath,bm,bbm,amsthm, amssymb}
\usepackage[foot]{amsaddr}
\usepackage{mathtools}
\usepackage{fullpage}
\usepackage{color}
\usepackage{dsfont}
\usepackage{etoolbox}
\usepackage{comment}
\usepackage{pgf}
\usetikzlibrary{calc}
\usetikzlibrary{patterns}
\usetikzlibrary{arrows}
\usetikzlibrary{decorations.pathreplacing}
\usepackage[utf8]{inputenc}
\usepackage{pgfplots}

\usepackage{tcolorbox}
\usepackage{adjustbox}
\usepackage[final]{pdfpages}
\usepackage[ruled,vlined,linesnumbered,algosection,resetcount]{algorithm2e}
\usepackage{blkarray}
\usepackage{arydshln}
\usepackage{wrapfig}
\usepackage{enumitem}

\theoremstyle{plain}
\usepackage{amsthm}
\newtheorem{theorem}{Theorem}
\newtheorem{proposition}[theorem]{Proposition}
\newtheorem{corollary}[theorem]{Corollary}
\newtheorem{lemma}[theorem]{Lemma}

\theoremstyle{definition}
\newtheorem{definition}[theorem]{Definition}

\newtheorem{example}[theorem]{Example}

\theoremstyle{remark}
\newtheorem{remark}[theorem]{Remark}

\numberwithin{theorem}{section}

\theoremstyle{plain}

\newtheorem*{assumptionS*}{Assumption S}

\newtheorem*{assumptionT*}{Assumption T}

\newtheorem*{assumptionU*}{Assumption U}

\newtheorem*{assumptionR*}{Assumption R}

\newtheorem*{assumptionP*}{Assumption P}

\newtheorem*{assumptionC1*}{Assumption C1}

\newtheorem*{assumptionC2*}{Assumption C2}

\newtheorem*{modelA*}{Model A}

\newtheorem*{modelB*}{Model B}

\def\R{\mathbb R}

\def\Z{\mathbb Z}

\def\P{\mathbb P}

\def\E{\mathbb E}
\def\eps{\varepsilon}

\def\00{\mathbf 0}

\def\bee{\begin{example}}
\def\ene{\end{example}}
\def\bet{\begin{theorem}}
\def\ent{\end{theorem}}
\def\bec{\begin{corollary}}
\def\enc{\end{corollary}}
\def\bep{\begin{proof}}
\def\enp{\end{proof}}

\def\om{\omega}

\def\la{\lambda}

\def\N{\mathbb N}
\def\mc{ \mathcal}

\def\one{\mathds1}

\def\d{\mathrm d}

\renewcommand\leq{\leqslant}

\renewcommand\geq{\geqslant}

\def\La{\Lambda}

\DeclarePairedDelimiter{\abs}{\lvert}{\rvert}

\def\bel{\begin{lemma}}
\def\enl{\end{lemma}}
\def\im{\item}
\def\been{\begin{enumerate}}
\def\enen{\end{enumerate}}

\def\k_d{\kappa_d}

\def\k{\kappa}

\def\bepr{\begin{proposition}}
\def\enpr{\end{proposition}}
\def\bede{\begin{definition}}
\def\ende{\end{definition}}
\def\bere{\begin{remark}}
\def\enre{\end{remark}}

\def\bfxi{\boldsymbol{\xi}}
\def\bfXi{\boldsymbol{\Xi}}
\newcommand{\oneJ}[1]{\mathbf{#1}_1^J}
\newcommand{\onenJ}[2]{\mathbf{#1}_{1,#2}^J}

\begin{document}

\author[C. Hirsch]{Christian Hirsch$^{1,2}$}
\email{hirsch@math.au.dk}

\author[M. Petr\'akov\'a]{Martina Petr\'akov\'a$^{3}$}
\email{petrakova@karlin.mff.cuni.cz}

\address{$^1$DIGIT Center, Aarhus University, Finlandsgade 22, 8200 Aarhus N, Denmark}
\address{$^2$Department of Mathematics, Aarhus University, Ny Munkegade 118, 8000 Aarhus C, Denmark}
\address{$^3$Department of Probability and Mathematical Statistics, Faculty of Mathematics and Physics, Charles University, Sokolovsk\'a 83, 186 74 Prague 8, Czech Republic}
\title{Large-deviation analysis for canonical Gibbs measures}
\thanks{{\it Authors' note:} This version of the article has been accepted for publication, after peer review but is not the Version of Record and does not reflect post-acceptance improvements, or any corrections. The Version of Record is available online at: https://doi.org/10.1007/s10955-025-03451-8}
\maketitle
\thispagestyle{empty}

\begin{abstract}
In this paper, we present a large-deviation theory developed for functionals of canonical Gibbs processes, i.e., Gibbs processes with respect to the binomial point process. We study the regime of a fixed intensity in a sequence of increasing windows. Our method relies on the traditional large-deviation result for local bounded functionals of Poisson point processes noting that the binomial point process is obtained from the Poisson point process by conditioning on the point number. Our main methodological contribution is the development of coupling constructions allowing us to handle  delicate and unlikely pathological events. The presented results cover three types of Gibbs models - a~model given by a~bounded local interaction, a~model given by a non-negative possibly unbounded increasing local interaction and the hard-core interaction model. The derived large deviation principle is formulated for the distributions of individual empirical fields driven by canonical Gibbs processes, with its special case being a large deviation principle for local bounded observables of the canonical Gibbs processes. We also consider unbounded non-negative increasing local observables, but the price for treating this more general case is that we only get large-deviation bounds for the tails of such observables. Our primary setting is the one with periodic boundary condition, however, we also discuss generalizations for different choices of the boundary condition.
\end{abstract}

\section{Introduction}
\label{sec:intro}
Gibbs processes are among the most fundamental and vigorously investigated models in statistical mechanics. In the lattice setting, they include the celebrated Ising and Potts models. In the continuum, the Widom-Rowlinson model or multibody-systems interacting with Lennart-Jones potential are key examples. We refer the reader to \cite{friedli,georgiiBook,ruelle} for an excellent introductions into the topic as well as to 
\cite{coll,vogel,chan,dickson,jansen,zass} for further recent results in this research area.

%
%
Moreover, through their definition, the understanding of the behavior of Gibbs processes has close connections to the mathematical theory of large deviations. This has also ramifications for central quantities in statistical physics such as the free energy or the pressure. In particular, in the setting of lattice-based Gibbs processes, {large deviation principles (LDPs) are shown under very general conditions~\cite{GeorgiiLDPZ,rassoul}}.

%
%
For Gibbs processes in the continuum, the large-deviation behavior has also received some research interest \cite{houdebert,georgii3,GeorgiiLDP}. In comparison to the classical theory, new phenomena arise that are completely different. {For lattice-Gibbs processes with nearest-neighbor interactions, the number of neighbors of a site is fixed. For continuum processes, the number of relevant neighbors is random and can be unbounded. This leads to novel localization effects such as the ones described in \cite{harel}.} Hence, {any investigation is} substantially more delicate. This is because, for instance, in the continuum pathological configurations could arise where many points cluster in a small volume of the space. Moreover, the original large-deviation results provided in \cite{georgii3,GeorgiiLDP} apply to a rather restricted class of functionals, namely, loosely speaking, averages of scores with a deterministic dependence range. This has motivated more recent large-deviation works where this condition is replaced by less restrictive assumptions \cite{HJT20LowerLDP,yukLDP,yukLDP2}.

%
%
In addition to the restrictions on the admissible observables, all of the above results concern the grand-canonical setting. That is, the Gibbs point processes are considered with respect to a Poisson point process. However, the case of canonical Gibbs point processes is equally important. Here, the Gibbs point processes are constructed with respect to the binomial point processes, i.e., with respect to a~fixed number of points. In particular, the delicate arguments for the central examples of Coulomb and Riesz gases often sensitively depend on the number of points being fixed \cite{dd1,dd2,serf1,serfBook}. With the exception of the very specific example of Coulomb and Riesz gases, we are not aware of a systematic study of {LDPs} for canonical Gibbs processes \cite{serfInv}.

%
%
Our main results are the first steps {towards} closing the gap between the {large-deviation} analysis of grand-canonical and canonical Gibbs processes. More precisely, we prove the following three main results.
\been
\im We establish the LDP for bounded local observables of two types of canonical Gibbs processes - processes driven by a local bounded interaction (see Theorem \ref{thm:ldp}) and processes driven by a possibly unbounded non-negative increasing local interaction (see Theorem \ref{thm:ldpunbounded}). Particularly, Theorem \ref{thm:ldp} can be seen as the analog of the grand-canonical LDPs from  
\cite{GeorgiiLDP} to the setting of canonical Gibbs processes.

\im Despite the importance of these theorems, in applications in statistical physics and stochastic geometry, many functionals are not bounded. Hence, we also discuss some extensions for non-negative observables such as edge counts in the random geometric graph. The price for treating this more general case is that we are only able to describe the large-deviation upper bound for the lower tails and large-deviation lower bound for both tails of the distributions of these unbounded observables. Corollaries \ref{lemma:UBunbounded1} and \ref{cor:LBunboundedscore} consider the case with bounded interaction, while Corollaries \ref{cor:UnbUBUnbscore} and \ref{cor:UnbLBUnbscore} deal with the model from Theorem \ref{thm:ldpunbounded}.

\im A drawback of Theorem \ref{thm:ldpunbounded} is that it assumes finiteness of the interaction function. Hence, the important case of hard-core interaction is excluded. However, our methodology is strong enough to deal also with such questions. To illustrate this potential, in Theorem \ref{cor:HCLDP} we show the LDP for the particularly important case of the hard-core Strauss process under some necessary constraints on the intensity.
\enen

%
%
The main idea of the proof is to rely on known results for the grand-canonical setting by using that the binomial point process can be considered as a Poisson point process conditioned on the number of points. While in general, the event of a Poisson point process having a specific number of points is highly unlikely, this probability is exactly of the right scale for a large-deviation analysis. Since the event of having a fixed number of points is closed in the considered $\tau_{\mc{L}^o}$-topology (which differs from the topology of weak convergence), this interpretation is extremely helpful for deriving the large-deviation upper bounds. However, the analysis of the lower bounds is more delicate. There we require a coupling construction to show that the considered conditioning can be approximated arbitrarily closely with open sets.  

%
%
This construction is particularly delicate for Theorem \ref{thm:ldpunbounded} and Corollaries \ref{cor:LBunboundedscore} and \ref{cor:UnbLBUnbscore}, where we consider not necessarily bounded observables of Gibbs processes with not necessarily bounded interaction {functions}. The challenge is that the corresponding result in the grand-canonical setting \cite{HJT20LowerLDP} heavily relies on a coupling construction involving thinning and sprinkling operations of the Poisson point process. However, such operations are not available in the canonical setting where the number of points is held fixed. Therefore, we need a new idea in the canonical setting. Here, our approach is to construct a coupling by replacing the addition/deletion of points with a ``move'' operation where we transfer points from densely occupied parts of the sampling windows into more sparse areas. This operation is highly delicate since we not only need to guarantee that the move operation establishes the desired configurations in the regions where points are moved to, but, simultaneously, we need to guarantee that the move operation does not cause issues in the regions where the points were taken from. Such a move operation was also successfully used in the analysis of Gibbsian systems with Riesz-type interactions \cite{dd2}.

For hard-core interactions, the construction needs to be refined because configurations with close points are not just improbable, they are forbidden. For the prototypical example of the Strauss process with hard-core interaction, we explain how to circumvent this problem. First, we construct a suitable coupling of the binomial and the Poisson point process. Then, once we can work with the Poisson process, we approximate the possibly infinite interaction by its bounded counterpart and implement a thinning to ensure that inadmissible configurations are avoided in the coupling.

The rest of the manuscript is organized as follows. Our main results are stated in Section \ref{sec:model}. Section~\ref{sec:proof} contains the proof of the LDP for the model with bounded interaction, while Sections \ref{sec:proofunbounded} and \ref{sec:prooflowunbounded} deal with the proof in the unbounded case. Section \ref{sec:HC} is devoted to the proof of the LDP for the hard-core process and Section \ref{sec:DifHam} concludes our paper with the proof of the LDP for processes with different boundary conditions. 
\section{Model and main results}\label{sec:model}

%
%
\subsection{Notation and model description} 

Let $(\Theta,\mathcal{A},\P)$ be a probability space. Let $\la > 0$ be {the} intensity and $d \geq 2$ {the} dimension. For $n \in \N$ denote by $W_n := \big[-\frac{1}{2}\left(\frac{n}{\la}\right)^{\frac{1}{d}},\frac{1}{2}\left(\frac{n}{\la}\right)^{\frac{1}{d}}\big]^d$ the centered cube with volume ${\abs{W_n} = } \frac{n}{\la}$ in $\R^d$. Clearly $W_n \uparrow \R^d$. Denote by $b(x,r)$ the closed ball with radius $r$  and center $x$ in $\R^d$ and $b_{r} := b(0,r)$, by $v_d$ the volume of the unit ball in $\R^d$ and by $Q_r(x) := x + \left[-\frac{r}{2},\frac{r}{2}\right)^d$, $Q_r:= Q_r(0)$. {Let $\mc{B}^d$ denote the Borel $\sigma$-algebra on $\R^d$.}  For $A \in \mc{B}^d$ and $a \in \R$ denote $aA := \{a x : x \in A \}$.  We will use $\abs{A}$ to denote both the number of points (if $A$ is countable) and the $d$-dimensional Lebesgue measure. The interpretation will always be clear from the context.

Denote by $(\Omega,\mc{F})$ the space of all {simple locally finite counting measures} on $(\R^d,\mc{B}^d)$ with the  smallest $\sigma$-algebra $\mc{F}$ such that the counting variable  $\om \rightarrow \om(B)$ is measurable for any $B \in \mc{B}^d$. {Here, $\om(B)$ denotes the number of points from $\om$ in the set $B$.} {As usual}, we treat $\om \in \Omega$ both as a measure on $\R^d$ and a subset of $\R^d$. Denote by $\Omega_f$ the set of finite measures {in} $\Omega$ and by $\Omega_\La$ the set of measures {in} $\Omega$ with support in $\La$, $\La \in \mc{B}^d$. For $\om \in \Omega$ denote by $\om^{(n)}:= \sum_{z \in w_n\mathbb{Z}^d} ((\om \cap W_n) + z)$, where $w_n: = \left(\frac{n}{\la}\right)^{\frac{1}{d}}$ is the side-length of $W_n$, the periodized version of $\om \cap W_n$. {For $x \in W_n$} denote by $b^{(n)}(x,r)$ the periodized version of $b(x,r)$ in~$W_n$ and similarly $Q^{(n)}_{r}(x)$.

\bede
Let $h: \Omega \rightarrow \R$ be a measurable function. Then, we say that it is
\begin{itemize}
    \item[$\bullet$] \textbf{$r$-local}  (or \textbf{local}) for $ r > 0$, if it holds that $h(\omega) = h(\omega \cap b_{r})$ for all $\om \in \Omega$, 
    \item[$\bullet$] \textbf{bounded}, if there exists $ c > 0$ such that $\abs{h(\om)} \leq c$ for all $\om \in \Omega$,
     \item[$\bullet$] \textbf{increasing} if $h(\om)\leq h(\om\cup \{x\})$ for any $\om \in \Omega$, $x \in \R^d$,
    \item[$\bullet$] \textbf{{quasi local bounded}}, if there exists $\{h_r\}_{r=1}^\infty$ such that $h_r \uparrow h$ point-wise and each $h_r$ is bounded and $r$-local,
    \item[$\bullet$] \textbf{cardinality-bounded by a sequence\footnote{Let us note that we will always assume that $\{M_b\}_{b \in \N}$ is non-decreasing and $\lim_{b \rightarrow \infty} M_b = \infty$. } $\{M_b\}_{b \in \N}$}, if for any $ b \in \N$ there exists $ M_b > 0$ such that $h(\om) \leq M_b$ for any $\om \in \Omega$ with $\abs{\om} \leq b$.
\end{itemize}
\ende

Denote by $Bi(n,p)$ the binomial distribution with parameters $n$ and $p$, by $Pois(n)$ the Poisson distribution with mean $n$ and by $Unif\,(0,1)$ the uniform distribution on the interval $(0,1)$. Furthermore, denote by $B_n$ the binomial point process in $W_n$ with $n$ points and by $\mathbb{B}_n$ its distribution, by $\pi_n$ the Poisson point process in $W_n$ with intensity $\la$ and by $\Pi_n$ its distribution. The Poisson point process with intensity $\la$ in $\R^d$ will be denote by $\pi$ and its distribution by $\Pi$. 

In this paper, we work with the following binomial Gibbs setting in order to be compatible with the formalism from \cite{GeorgiiLDP}. Binomial Gibbs models with different Hamiltonians (particularly with different boundary conditions) will be discussed in {Sections~\ref{secresults:DifHam} and \ref{sec:DifHam}}.

\bede \label{def:BGPP}
Let $n \in \N$. Let $V: \Omega \rightarrow \R$ be a measurable function called the \textbf{interaction function}. Define the \textbf{Hamiltonian} in $W_n$ as $H_n(\om) = \sum_{x \in \om\cap W_n} V(\om^{(n)}-x)$, $\om \in \Omega_{W_n}$. Then, the \textbf{binomial Gibbs point process} in {the} window $W_n$ with Hamiltonian $H_n$ is defined as a point process with distribution
$$
\d P_n (\om) = \frac{1}{\tilde{Z}_n} \exp(-H_n(\om)) \d \mathbb{B}_n(\om),
$$
where the normalizing constant $\tilde{Z}_n = \int \exp(-H_n(\om)) \d \mathbb{B}_n(\om) $ is called {the} \textbf{partition function}. The binomial Gibbs point process will be denoted as $\rho_n$.
\ende

Let us {note that, since} we have a fixed number of points, adding a constant to the interaction function $V$ does not change the Gibbs measure $P_n$.

{The crucial ingredient of the proofs is the well-known fact that {the} binomial point process is just the Poisson point process conditioned on the number of points, i.e.\,$\d \mathbb{B}_n(\om) = \frac{1}{\Pi_n(\abs{\om} = n)} \one \left[ \abs{\om} = n\right] \d \Pi_n(\om)$. Therefore, \begin{equation*}
    \d P_n (\om) = \frac{1}{Z_n} \exp(-H_n(\om)) \one \left[ \abs{\om} = n\right] \d \Pi_n(\om),
\end{equation*} where $Z_n = \int \exp(-H_n(\om)) \one \left[ \abs{\om} = n\right] \d \Pi_n(\om)= \Pi_n(\abs{\om} = n) \tilde{Z}_n$ will also be called the partition function.}

Without any assumptions on $V$, the measures $P_n$, {$n \in \N$}, may not be well defined. Let us now discuss the three settings for $V$ considered in this paper.
\begin{enumerate}
    \item Let $V$ be $r$-local and bounded. Then, $0< \tilde{Z}_n < \infty$ {holds for all $n \in \N$} and therefore $P_n$ is well-defined {for all $n \in \N$}. 
    \item Let $V \geq 0$ be {quasi local bounded}. Then, $\tilde{Z}_n < \infty$ {holds for all $n \in \N$}, but we are not able to show {for any $n \in \N$} (in general) that $\tilde{Z}_n > 0$. However, this setting is only considered in Lemma~\ref{lemma:UBunbounded2} whose claim still holds even if $\tilde{Z}_n = 0$ {for all $n \in \N$} (see more in Section \ref{sec:proofupunbounded}).
    \item Let $V \geq 0$ be $r$-local, increasing and cardinality-bounded. Then, $0< \tilde{Z}_n < \infty$ {holds} for any $n \in \N$ and therefore $P_n$ is well-defined {for all $n \in \N$}. 
\end{enumerate}

A key step in our paper will be to derive LDP for functions of the process $\rho_n$ {of} the following form.

\bede 
Let $\xi: \Omega \rightarrow \R$ be a measurable function called the \textbf{score function}. Then, we denote 
$$
\Xi_n(\om) := \frac{1}{n}\sum_{x \in \om\cap W_n} \xi(\om^{(n)}-x). 
$$
\ende

In this paper, we will usually work with a set of score functions $\{\xi_1,\dots,\xi_J\}$ for some $J \in \N$. In that case we denote $\Xi_{j,n}(\om) := \frac{1}{n}\sum_{x \in \om\cap W_n} \xi_j(\om^{(n)}-x)$, $j \in \{1, \dots, J\}$. 
To ease the notation, we will work with vectors of functions $\bfxi_1^J := \left(\xi_1,\dots,\xi_J\right)^T$ and $\bfXi_{1,n}^J:=\left(\Xi_{1,n},\dots,\Xi_{J,n}\right)^T$ and denote by $\mathbf{B}_1^J:= B_1 \times \dots \times B_J$ the Cartesian product of sets $B_j \subset \R$. Then, we can write $\bfxi_1^J \in \mathbf{B}_1^J$, instead of $\xi_j \in B_j$ for all $j \in \{1,\dots, J\}$, and, if $P$ is a measure on $(\Omega, \mathcal{F})$, then $P(\bfxi_1^J):= \left(P(\xi_1),\dots,P(\xi_J)\right)^T$.  { Furthermore, by calling $\mathbf{B}_1^J$ an open (closed) set, we mean that $\mathbf{B}_1^J= B_1 \times \dots \times B_J$, where $B_i$ are open (closed) subsets of $\R$.}

Let us finish this section by summarizing the {more general} notation from paper \cite{GeorgiiLDP}{,} which will be needed for our investigations. Firstly, they deal with more general case of marked point processes with marks from some space $E$. As we consider only the unmarked case here, we can choose $E$ to be any singleton. For simplicity, we omit the notation of the marks altogether.

First, denote by $\mc{M}$ the space of finite measures on $(\Omega, \mc{F})$, by $\mc{P}$ the set of probability measures on $(\Omega, \mc{F})$ and by $\mc{P}_s$ the set of stationary probability measures on $(\Omega, \mc{F})$ {with finite intensity}. For a stationary measure $P \in \mc{P}_s$ denote

\begin{itemize}
    \item[$\bullet$] by $P^o$ the \textbf{Palm measure} of $P$ (see Rem. 2.1 in \cite{GeorgiiLDP}). Particularly, $P^o$ is the unique finite measure on $(\Omega,\mathcal{F})$ such that for any {measurable function $f: \R^d\times\Omega\rightarrow \left[0,\infty\right)$} the following equality holds
\begin{equation*}
    \int \int f(x,\om-x)\, \d \om (x)\, \d P(\om) = \int \int f(x,\om)\, \d x\, \d P^o(\om). 
\end{equation*}
    \item[$\bullet$] by $I(P)$ the \textbf{specific entropy} of $P$, (see (2.12) in \cite{GeorgiiLDP}), defined as 
    \begin{equation*}
   I(P) = \lim_{n \rightarrow \infty} \frac{1}{\abs{W_n}} I(P_{W_n};\Pi_n), 
\end{equation*}
where $I(P;Q)$ is the relative entropy of the measure $P$ w.r.t.\,the measure $Q$ (for definition see (2.10) in \cite{GeorgiiLDP}) and $P_{W_n}$ is the restriction of $P$ to the window $W_n$. 
\end{itemize}

Throughout this paper, $\mc{M}$ will be equipped with {the} {\bf $\tau_{\mc{L}^o}$-topology}, which is defined as the smallest topology such that the mapping $e_g : \mc{M} \rightarrow \R$, $e_g(m) = m(g)$ is continuous for any function $g \in \mc{L}^o := \{g: \Omega \rightarrow \R \text{ measurable, local and bounded} \}$  (for more details see \cite{GeorgiiLDP}, after Remark 2.1). {The $\tau_{\mc{L}^o}$-topology has a basis consisting of the sets {of} the form 
\begin{equation}\label{rem:basis}
    U_{{\bf g}_1^J;m} :=\{\tilde{m} \in \mc{M}: \max_{1\leq j \leq J} \abs{m(g_j)-\tilde{m}(g_j)}< \delta\}
\end{equation}
where $J \in \N$, $\delta > 0$ and {$\oneJ{{\bf g}}$} are bounded and local functions \cite{EiSch98}.}

\bede 
Let  $n \in \N$ and $\om \in \Omega$. We define {the} \textbf{individual empirical field $R^o_{n,\om}$} in $W_n$ as 
$$R^o_{n,\om} :=\frac{1}{\abs{W_n}} \sum_{x\in \om \cap W_n} \delta_{\om^{(n)}-x} = \frac{\la}{n} \sum_{x\in \om \cap W_n} \delta_{\om^{(n)}-x}.$$
\ende

It holds that $R^o_{n,\om} \in \mc{M}$ and its random counterpart $R^o_{n,\rho_n}$, where $\rho_n$ is the binomial Gibbs point process, is a well-defined random element in {$(\mc{M},\mathcal{B}(\tau_{\mc{L}^o}))$, where $\mathcal{B}(\tau_{\mc{L}^o})$ is the Borel $\sigma$-algebra generated by the  $\tau_{\mc{L}^o}$-topology}.

We conclude this section by recalling the definitions of semi-continuity. A function $G : \mc{M} \rightarrow \R \cup \{-\infty,+\infty\}$ is called \textbf{lower semi-continuous (lsc)}, if the level sets $L_b := \{m \in \mc{M}: G(m) \leq b  \}$ are closed in $(\mc{M},\tau_{\mc{L}^o})$ for all $b \in \mathbb{R}$, and it is called \textbf{upper semi-continuous (usc)}, if $-G$ is {lsc}. For any function $G : \mc{M} \rightarrow \R \cup \{-\infty,+\infty\}$ we denote by $G_{lsc}$ its lower semi-continuous version, i.e.\,the largest {lsc} function $F$ such that $G \geq F$, and  by $G^{usc}$ its upper semi-continuous version, i.e., the smallest {usc} function $F$ such that $G \leq F$.

Throughout this paper we will use the following properties of lsc/usc functions:

\bel \label{rem:lsc/uscproperties} {It holds that
\begin{enumerate}
    \item {any} finite sum of lsc functions is lsc {and any} finite sum of usc functions is usc,
    \item a continuous function is both {lsc} and {usc},
    \item if $G$ is lsc, then $G_{lsc} = G$, if $G$ is usc, then $G^{usc} = G$,
   \item if $U_b = \{m \in \mc{M}: G(m) \geq b  \}$ is closed in $(\mc{M},\tau_{\mc{L}^o})$ for any $b \in \R$, then $G$ is usc,
   \item if $C \subset \mc{M}$ is closed in the $\tau_{\mc{L}^o}$-topology, then $G(m) = +\infty \cdot \one \left[m \in C^c\right]$ is lsc.  
\end{enumerate}}
\enl

\bep
All of these follow easily from the definitions. 
\enp

We will also need the following (easy to prove) observation.

\bel \label{rem:lsc/uscproperties2} {
Let $h: \Omega \rightarrow \R$ be a measurable local bounded function. Define the function $G: \mc{M} \rightarrow \R \cup \{-\infty,+\infty\}$ as $G(m) := +\infty \cdot \one \left[m(h) \in A\right]$ for some $A \subset \R$ (with $+\infty \cdot 0 = 0$). Then, 
   \begin{enumerate}
       \item if $A$ is closed, then $G$ is {usc},
   \item if $A$ is open, then $G$ is {lsc}.    \end{enumerate}}
\enl

\bep
{These points follow easily from Points 4 and 5 in Lemma \ref{rem:lsc/uscproperties}.}
\enp

To finish this section, we note that in the proofs of our results, we will use frequently the following (traditional) observations. Let $a_n,b_n \geq 0$, $n \in \N$, then 
\begin{alignat*}{3}
     \limsup \frac{\la}{n} \log a_n &\geq \limsup \frac{\la}{n} \log b_n &\implies \limsup \frac{\la}{n} \log (a_n + b_n) &= \limsup \frac{\la}{n} \log a_n, \\
    \liminf \frac{\la}{n} \log a_n &> \limsup \frac{\la}{n} \log b_n &\implies \liminf \frac{\la}{n} \log (a_n - b_n) &= \liminf \frac{\la}{n} \log a_n.
\end{alignat*}

\subsection{Main results}\label{sec:mainresults}

Our main result is the {LDP} for the sequence of random individual empirical fields $R^o_{n,\rho_n}$ driven by binomial Gibbs processes $\rho_n$. We consider two cases - the case with bounded local interaction function $V$ and the case where $V$ is bounded from below, local, increasing and cardinality bounded. In both of these cases, however, $V$ can only attain finite values. {A large-deviation result for a process with hard-core interaction (i.e., a process where $V$ can possibly attaining the value $+\infty$) is presented in Section \ref{secresults:HC}. As our primary setting is a one with the periodic boundary condition, we also generalize our results for the cases with other boundary conditions (see Section~\ref{secresults:DifHam}), albeit only in the bounded case.}

{Our first result is the LDP for the Gibbs model with local bounded interaction $V$.}

\bet[{The} LDP in the bounded case]
\label{thm:ldp}
Assume that the interaction function $V$ is $r$-local and bounded. Then, the sequence of random individual empirical fields $R^o_{n,\rho_n}$ driven by the binomial Gibbs processes $\rho_n$ satisfies an LDP with speed {$\abs{W_n}$} and rate function $\mc{J}:\mc{M} \rightarrow \R \cup \{+\infty\} $, i.e., 
\begin{align}
\begin{split}\label{def:LDP}
    \limsup{\frac{1}{\abs{W_n}}\log \P(R^o_{n,\rho_n} \in C)} &\leq - \inf_{m \in C} \mc{J}(m) \quad \text{ for any } C \subset \mc{M} \text{ closed,} \\
    \liminf{\frac{1}{\abs{W_n}}\log \P(R^o_{n,\rho_n} \in O)} &\geq - \inf_{m \in O} \mc{J}(m) \quad \text{ for any } O \subset \mc{M} \text{ open,} 
\end{split}
\end{align}
where 
\begin{equation} \label{def:J}
\mc{J}(m) := \tilde{\mc{J}}(m)-\inf_{m \in \mc{M}} \tilde{\mc{J}}(m)
\end{equation}
with
\begin{equation} \label{def:tildeJ}
   \tilde{\mc{J}}(m):=\begin{cases}
			I(P)+m(V), & \text{if } m = P^o \text{ for some } P \in \mc{P}_s \text{ and } m(\one) =  \la, \\
            + \infty, & \text{otherwise.}
		 \end{cases}
\end{equation}

\ent

\bere
Even though it is omitted in the notation, we should keep in mind that the rate function depends on the interaction function $V$.
\enre

{The proof of Theorem \ref{thm:ldp} is postponed to Section \ref{sec:proof}. It relies on proving the LDP for random vectors $\onenJ{\bfXi}{n}(\rho_n)$, with local and bounded score functions $\oneJ{\bfxi}$.} However, we can also derive large-deviation lower and upper bounds for tails of more general score functions~$\oneJ{\bfxi}$. While the large-deviation upper bound can be derived only for the lower tail, we can get the {large-deviation} lower bound for both upper and lower tails of the averages of possibly unbounded score functions. {Note that for unbounded score functions, localization effects can appear \cite{harel}, which imply that the large deviations are no longer of volume order. Therefore, it is not surprising that in such cases our methods do not give  upper and lower bounds for intervals.
}

\bec[]\label{lemma:UBunbounded1}
Let $J \in \N$ and $\mathbf{a}\in \R^J$. Assume that the interaction function $V$ is $r$-local and bounded and that the score functions $\bfxi_1^J$ are non-negative and {quasi local bounded}. Then,
\begin{align*}
    \limsup&{{ \frac{1}{\abs{W_n}} }\log P_n(\bfXi_{1,n}^J \leq \mathbf{a})} \\& \qquad \leq - \inf \{I(P) + P^o(V): P \in \mc{P}_s, P^o({\bf1}) = \la, P^o(\bfxi_1^J) \leq \la \mathbf{a} \} + \inf_{m \in \mc{M}} \tilde{\mc{J}}(m),  \\
    \liminf&{{ \frac{1}{\abs{W_n}} }\log P_n(\bfXi_{1,n}^J > \mathbf{a})} \\& \qquad \geq - \inf \{I(P) + P^o(V): P \in \mc{P}_s, P^o({\bf1}) = \la, P^o(\bfxi_1^J) > \la \mathbf{a} \} + \inf_{m \in \mc{M}} \tilde{\mc{J}}(m).
\end{align*}
\enc

\bep
We can write $\one \left[\bfXi_{1,n}^J(\rho_n) \leq \mathbf{a}\right] = \one \left[ R^o_{n,\rho_n} \in  \{m \in \mc{M}: m(\bfxi_1^J) \leq \la \mathbf{a} \} \right]$ and $$\one \left[\bfXi_{1,n}^J(\rho_n) > \mathbf{a}\right] = \one \left[ R^o_{n,\rho_n} \in  \{m \in \mc{M}: m(\bfxi_1^J) > \la \mathbf{a} \} \right],$$ since $\bfXi_{1,n}^J(\om) = \frac{1}{\la} R^o_{n,\om}(\bfxi_1^J)$.
Hence, it is enough to prove that the set $\{m \in \mc{M}: m(\bfxi_1^J) \leq \la \mathbf{a} \}$ is a closed subset of $\mc{M}$ and use Theorem \ref{thm:ldp}.

Denote $\bfxi_{1,r}^J:= \left( \xi_{1,r}, \dots, \xi_{J,r} \right)^T$ {the approximating sequence. Then, clearly
\begin{equation*}
   \{m \in \mc{M}: m(\bfxi_1^J) \leq \la \mathbf{a} \} \subset \bigcap_{r = 1}^\infty \{m \in \mc{M}: m(\bfxi_{1,r}^J) \leq \la \mathbf{a} \}.
\end{equation*}
Using the monotone convergence theorem, we can prove that if $m \in \mc{M}$ satisfies $\sup_{r \in \N}m(\bfxi_{1,r}^J) \leq \la \mathbf{a}$, then $m(\bfxi_1^J) \leq \la \mathbf{a}$. Therefore, 
\begin{equation*}
   \{m \in \mc{M}: m(\bfxi_1^J) \leq \la \mathbf{a} \} \supset \bigcap_{r = 1}^\infty \{m \in \mc{M}: m(\bfxi_{1,r}^J) \leq \la \mathbf{a} \}.
\end{equation*}
 Since $\xi_{j,r}$ are bounded and local, we get that the set $M_r:=\{m \in \mc{M}: m(\bfxi_{1,r}^J) \leq \la \mathbf{a} \}$ is closed for all $r \in \N$. Therefore, also the set $  \bigcap_{r = 1}^\infty M_r =\{m \in \mc{M}: m(\bfxi_1^J) \leq \la \mathbf{a} \} $ is closed, which finishes the proof.}
\enp

\bec[] \label{cor:LBunboundedscore} 
Let $J \in \N$ and $\mathbf{a} \in \R^J$. Assume that the interaction function $V$ is $r$-local and bounded by a constant $c$. Assume that the score functions $\bfxi_1^J$ are bounded from below by a constant $\tilde{c}$, $r$-local, increasing and cardinality-bounded by a sequence $\{M_b\}_{b \in \N}$. Then, the following lower bound holds
\begin{align*}
    \liminf&{{ \frac{1}{\abs{W_n}} }\log P_n(\bfXi_{1,n}^J < \mathbf{a})} \\& \qquad \geq - \inf \{I(P) + P^o(V): P \in \mc{P}_s, P^o({\bf1}) = \la, P^o(\bfxi_1^J) < \la \mathbf{a}\} + \inf_{m \in \mc{M}} \tilde{\mc{J}}(m).
\end{align*}
\enc

{While the proof of Corollary \ref{lemma:UBunbounded1} is a simple application of Theorem \ref{thm:ldp}, the large-deviation lower bound presented in Corollary \ref{cor:LBunboundedscore}} cannot be proved {analogously} by finding an open set $O \subset \mc{M}$ such that $\{\bfXi_{1,n}^J(\rho_n) < \mathbf{a}\} = \{R^o_{n,\rho_n} \in O\}$, {since countable intersection of open sets is not necessarily an open set. The proof of Corollary \ref{cor:LBunboundedscore} is postponed to Section \ref{sec:proofunbounded}.}

Let us move to the second case, where we prove the LDP for binomial Gibbs point process with unbounded interaction function $V$. Let us note that{,} while we weaken our demands in the "bounded" part, requesting only the cardinality-bounded assumption, we add the assumption that $V$ is increasing. Particularly, Theorem \ref{thm:ldp} is not a special case of Theorem \ref{thm:ldpunbounded}.

\bet[{The} LDP in the unbounded case]\label{thm:ldpunbounded}
Assume that the interaction function $V$ is bounded from below, increasing, $r$-local and cardinality-bounded by a sequence $\{M_b\}_{b\in \N}$. Then{,} $\lim { \frac{1}{\abs{W_n}} } \log Z_n$ exists and is equal to {$-\inf_{m \in \mc{M}} \tilde{\mc{J}}(m)$}, where $\tilde{\mc{J}}$ is defined in (\ref{def:tildeJ}) .

Furthermore, assume that $\lim { \frac{1}{\abs{W_n}} } \log Z_n > - \infty$. Then, the sequence of random individual empirical fields $R^o_{n,\rho_n}$ driven by binomial Gibbs point processes $\rho_n$ satisfies an {LDP with speed $\abs{W_n}$} and rate function $\mc{J}:\mc{M} \rightarrow \R \cup \{+\infty\}$ from (\ref{def:J}) ({i.e., it satisfies (\ref{def:LDP}))}. 
\ent

\bere
Let us recall that under the assumptions on $V$ from Theorem \ref{thm:ldpunbounded}, the Gibbs measures $P_n$ are always well-defined. Furthermore, we always have that  $\lim {\frac{1}{\abs{W_n}} \log Z_n   < +\infty}$. However, contrary to the bounded case, it can happen that $\lim { \frac{1}{\abs{W_n}} } \log Z_n = -\infty$. 
\enre

As in the bounded case, we can prove large-deviation upper and lower bound for tails of averages of more general score functions $\oneJ{\bfxi}$.

\bec[]\label{cor:UnbUBUnbscore}
Let $J \in \N$ and $\mathbf{a} \in \R^J$. Assume that the interaction function $V$ is bounded from below, increasing, $r$-local and cardinality-bounded by a sequence $\{M_b\}_{b\in \N}$ and {that $\lim \frac{1}{\abs{W_n}} \log Z_n$ is finite}. Assume that the score functions $\bfxi_1^J$ are non-negative and {quasi local bounded}. Then,
\begin{align*}
    \limsup&{{ \frac{1}{\abs{W_n}} }\log P_n(\bfXi_{1,n}^J \leq \mathbf{a})} \\& \quad \leq - \inf \{I(P) + P^o(V): P \in \mc{P}_s, P^o({\bf1}) = \la, P^o(\bfxi_1^J) \leq \la \mathbf{a} \} + \inf_{m \in \mc{M}} \tilde{\mc{J}}(m),  \\
    \liminf&{{ \frac{1}{\abs{W_n}} }\log P_n(\bfXi_{1,n}^J > \mathbf{a})} \\& \quad  \geq - \inf \{I(P) + P^o(V): P \in \mc{P}_s, P^o({\bf1}) = \la, P^o(\bfxi_1^J) > \la \mathbf{a} \} + \inf_{m \in \mc{M}} \tilde{\mc{J}}(m).
\end{align*}
\enc

\bep
This follows from Theorem \ref{thm:ldpunbounded} {similarly} as Corollary \ref{lemma:UBunbounded1} followed from Theorem \ref{thm:ldp}.
\enp

\bec[]\label{cor:UnbLBUnbscore}
Let $J \in \N$ and $\mathbf{a} \in \R^J$. Assume that the score functions $\bfxi_1^J$ and the interaction function $V$ are bounded from below, increasing, $r$-local and cardinality-bounded by a sequence $\{M_b\}_{b\in \N}$. Assume {that $\lim \frac{1}{\abs{W_n}} \log Z_n$ is finite}. Then,
\begin{align*}
    \liminf&{{ \frac{1}{\abs{W_n}} }\log P_n(\bfXi_{1,n}^J < \mathbf{a})} \\& \quad  \geq - \inf \{I(P) + P^o(V): P \in \mc{P}_s, P^o({\bf1}) = \la, P^o(\bfxi_1^J) < \la \mathbf{a} \} + \inf_{m \in \mc{M}} \tilde{\mc{J}}(m).
\end{align*}
\enc

The proofs of Theorem \ref{thm:ldpunbounded} and Corollary \ref{cor:UnbLBUnbscore} are postponed to Section \ref{sec:proofunbounded}.

\subsubsection{Hard-core process}\label{secresults:HC} 
In this section, we will deal with the process driven by a hard-core interaction. We will use the notion "$R$-neighbor of a point $x$", meaning a point that lies inside $b^{(n)}(x,R)$, where $n$ will be clear from the context. 

Let $R > 0$ be the hard-core interaction radius and denote by $N_n(\om)$ the number of points from $\om \cap W_n$ with at least one $R$-neighbor 
$$
N_n(\om) := \abs{\{x \in \om \cap W_n: \om^{(n)}(b(x,R)) \geq 2\}}, \, \om \in \Omega,\, n \in \N.
$$
We will work with the following interaction function
\begin{equation}\label{def:HCintfction}
    V(\om) =  \infty \cdot \one \left[\om(b(0,R)) \geq 2 \right].
\end{equation}
The corresponding Hamiltonian $H_n$ is equal to $$ H_n(\om) = \sum_{x \in \om\cap W_n} V(\om^{(n)}-x) = \infty \cdot \one \left[N_n(\om) \geq 1 \right].$$ Particularly, configurations with at least one point with an $R$-neighbor are forbidden under the corresponding Gibbs measure $P_n$.

\bere
The interaction function $V$ is clearly local, but not bounded and therefore we have to deal with it separately. Although we believe that a similar approach as in this section could be also applied to slightly more general setting with interaction function $V(\om)=  \infty \cdot \one \left[\om(b(0,R)) \geq 2 \right] + V_1(\om)$, where $V_1$ is a local bounded function, we do not pursue such generality in our paper. 
\enre

{To ensure that the hard-core binomial Gibbs process with $n$ points in $W_n$ (i.e., the Gibbs process with interaction function $V$ from \eqref{def:HCintfction}) is well defined, we must pose the following assumption on the intensity $\la$ (w.r.t.\,the hard-core interaction radius $R$)
\begin{equation} \label{ass:HCmodelLambda}
    \la R^d v_d < 1.
\end{equation}
This assumption guarantees that the partition function $Z_n$ corresponding to the interaction function $V$ satisfies $Z_n > 0$. It will also be used later in one of the auxiliary bounds (see Lemma \ref{lemma:HCAuxLBII}). 
 }

Our main result is the LDP for local bounded functionals of hard-core binomial Gibbs processes. Let us note that although we do not go to this generality in the hard-core case, we could again generalize {our result} as in the case of Theorem \ref{thm:ldp} to obtain LDP for the sequence of random individual empirical fields driven by hard-core binomial Gibbs processes.

\bet[{The} LDP for the hard-core process] \label{cor:HCLDP}
There exists the limit for the partition functions of the hard-core process
\begin{equation*} 
   \lim{{ \frac{1}{\abs{W_n}} }\log Z_n} = - \inf \{I(P) + P^o(V): P \in \mc{P}_s, P^o({\bf1}) = \la\} =: - A_{HC}.
\end{equation*} 
Let $J \in \N$, $\oneJ{C}\subset \R^J$ be a closed set and $\oneJ{U} \subset \R^J$ be an open set. Assume that the score functions $\oneJ{\bfxi}$ are $r$-local and bounded.  If $A_{HC} < \infty$, then 
\begin{align*}
    \limsup { \frac{1}{\abs{W_n}} } &\log P_n(\onenJ{\bfXi}{n} \in \oneJ{C}) \\ 
    &\leq - \inf \{I(P) + P^o(V): P \in \mc{P}_s, P^o({\bf1}) = \la, P^o(\oneJ{\bfxi}) \in \la \oneJ{C}\} + A_{HC},\\
    \liminf { \frac{1}{\abs{W_n}} } &\log P_n(\onenJ{\bfXi}{n} \in \oneJ{U}) \\ 
    &\geq - \inf \{I(P) + P^o(V): P \in \mc{P}_s, P^o({\bf1}) =\la, P^o(\oneJ{\bfxi}) \in \la \oneJ{U}\} + A_{HC}.
\end{align*}
\ent

\bere
Observe that $A_{HC} \geq 0$. Furthermore, observe that we can write
$$\inf \{I(P) + P^o(V): P \in \mc{P}_s, P^o({\bf1}) = \la\} = \inf \{I(P) : P \in \mc{P}_s, P^o({\bf1})= \la, P^o(\Omega_R)=0 \},$$ where $\Omega_R:= \{\om \in \Omega \,:\, \om(b(0,R)) \geq 2\}.$ 
\enre

\subsubsection{Different Hamiltonians}\label{secresults:DifHam}
In this section, we generalize the previous results for Gibbs models with different Hamiltonians. Particularly, we will consider models with different boundary conditions. For $n \in \N$ let $\hat{H}_n$ be our new Hamiltonian (a measurable mapping from $(\Omega,\mathcal{F})$ to $(\R,\mathbb{B})$) and consider the binomial Gibbs model
\begin{equation} \label{def:hatBGPP}
    \d \hat{P}_n(\om) = \frac{1}{\hat{Z}_n} \exp{(-\hat{H}_n(\om))}\, \d \mathbb{B}_n(\om).
\end{equation}
We assume, that we only work with such Hamiltonian $\hat{H}_n$, so that $\hat{P}_n$ is well defined. In the main result of this section, we consider the following two Hamiltonians $\hat{H}^1_n$ and $\hat{H}^2_n$. For any $\gamma \in \Omega$ and any interaction function $V$ denote
\begin{align} 
\begin{split}\label{def:difham}
    \hat{H}^1_n(\om)&:= \sum_{x \in \om\cap W_n} V(((\om\cap W_n)\cup(\gamma \cap W_n^c)) - x),\\
    \hat{H}^2_n(\om)&:= \sum_{x \in \om\cap W_n} V(((\om\cap W_n)\cup(\gamma \cap W_n^c)) - x) + \sum_{x \in \om\cap W_n} V((\gamma \cap W_n^c) - x).
\end{split}
\end{align}

Let us remark on the difference between the two Hamiltonians $ \hat{H}^1_n$ and $ \hat{H}^2_n$. {Imagine that $V$ is a $k$-wise interaction function (see Example \ref{ex:k-int}).} Then, the Hamiltonian $\hat{H}^2_n$ corresponds to the standard Gibbs process with conditional energy, i.e. it is the part of the energy of the whole system contributed by the points inside the window $W_n$. In the case of $ \hat{H}^1_n$, the interaction between the points inside the window is counted for each point (i.e.\,twice), whereas the interaction between a point inside and a point outside of the window $W_n$ is counted only once.

Now, we can present the main result, which is the large-deviation result for measures $\hat{P}_n$ for the Hamiltonians from (\ref{def:difham}).

\bet[] \label{lemma:difbc} 
Let $\gamma \in \Omega$ be a given boundary condition. Let $V$ be $r$-local and bounded by a constant $c_1$ and assume\footnote{This can always be assumed, since we can shift $V$ by a constant without making a change to the Gibbs measure it generates.} that $V(\{0\}) = 0$. Then, the Gibbs measures $\hat{P}_n$ given by {the} Hamiltonians $\hat{H}_n = \hat{H}^i_n$, $i = 1,2$, from (\ref{def:difham}), satisfy $$\lim_{n \rightarrow \infty} { \frac{1}{\abs{W_n}} } \log \hat{Z}_n = - \inf \{I(P) + P^o(V): P \in \mc{P}_s, P^o({\bf1}) = \la\} =: - A.$$  Furthermore,
\begin{align} \label{hatLDP}
    \begin{split}
         \limsup&{{ \frac{1}{\abs{W_n}} }\log \hat{P}_n(\onenJ{\bfXi}{n}  \in \oneJ{C})} \\
         &\leq - \inf \{I(P) + P^o(V): P \in \mc{P}_s, P^o({\bf1}) = \la, P^o(\oneJ{\bfxi})  \in \la \oneJ{C} \} + A,\\
         \liminf&{{ \frac{1}{\abs{W_n}} }\log \hat{P}_n(\onenJ{\bfXi}{n}  \in \oneJ{U})} \\
         &\geq - \inf \{I(P) + P^o(V): P \in \mc{P}_s, P^o({\bf1}) = \la, P^o(\oneJ{\bfxi}) \in\la \oneJ{U} \} + A,
    \end{split}
\end{align}
 hold for the following cases, {where $J \in \N$,}
 \begin{itemize}
 \item[(a)]  $\oneJ{C}\subset \R^J$ a closed set, $\oneJ{U} \subset \R^J$ an open set and $\oneJ{\bfxi}$ local bounded score functions,
 \item[(b)]  $\oneJ{U} = (-\infty, \mathbf{a})$ for any $\mathbf{a} \in \R^J$ and $\oneJ{\bfxi}$ bounded from below, increasing, $r$-local and cardinality-bounded by a sequence $\{M_b\}_{b\in \N}$,
 \item[(c)]  $\oneJ{C} = \left(-\infty, \mathbf{a}\right]$ and $\oneJ{U} = (\mathbf{a},\infty)$ for any $\mathbf{a} \in \R^J$ and $\oneJ{\bfxi}$ non-negative quasi local bounded. 
 \end{itemize} 
\ent

\bere
Under the assumptions of Theorem \ref{lemma:difbc} (particularly that $V$ is bounded), the binomial Gibbs model (\ref{def:hatBGPP}) with $\hat{H}_n = \hat{H}^i_n$, $i = 1,2$, is clearly well defined.
 \enre

The proof of Theorem \ref{lemma:difbc} is postponed to Section \ref{sec:DifHam} and it is based on a general result from Lemma \ref{lemma:hatLDP}. Similarly, we could use this lemma to derive large-deviation results for measures $\hat{P}_n$ for other suitable Hamiltonians. We also note that analogously as in the previous sections, one could generalize Theorem \ref{lemma:difbc} to full LDP for the sequence of individual empirical fields driven by point process $\hat{\rho}_n$ with distribution $\hat{P}_n$.

\subsection{Examples} \label{sec:Examples}

Before we proceed with the proofs, we present some examples of models and score functions to which our results can be applied.

Firstly, any non-negative score function $\xi$ is in fact {quasi local bounded}, since we can define $\xi_r(\om) = \min\{r,\xi(\om \cap b_r)\}$. Therefore, Corollaries \ref{lemma:UBunbounded1} and \ref{cor:UnbUBUnbscore} are always applicable, given that we have a suitable interaction function $V$. As examples of $V$ satisfying assumptions of Theorem \ref{thm:ldpunbounded}, we present the Strauss model and the bounded $k$-wise interaction model.

\bee[Binomial Strauss process] \label{ex:Strauss}
Let  $\gamma \in \left(0,1\right]$ and $r > 0$. The {\bf binomial Strauss process} (also called conditional Strauss process, see (6.15) in \cite{MollerWaagBook}) in {the} window $W_n$ with periodic boundary condition, $n$ points and parameters $\gamma$ and $r$, is the binomial Gibbs process with distribution 
\begin{equation*}\mathrm{d} P_n(\om) = \frac{1}{\tilde{Z}_n} \gamma^{S_r^{(n)}(\om)} \mathrm{d} \mathbb{B}_n(\om),\end{equation*} where $S_r^{(n)}(\om) := \sum_{\{x,y\} \subset \om, x\neq y} \one \left[x \in b^{(n)}(y,r)\right]$ is the number of pairs of points that are closer than $r$.

Let us make two remarks. Firstly, for $\gamma = 1$ we get the binomial point process. Secondly, in the standard Strauss model (i.e.\,with density w.r.t.\,the Poisson point process) we have an additional term $\beta^{\abs{\om}}$ in the density. However, in our case this term becomes a constant, so it can be omitted. It is an example of a~process with repulsive interactions, where the strength of the repulsion depends on the parameter $\gamma$. 

In our setting, we can say that the binomial Strauss process is the binomial Gibbs process $P_n$ given by {the} interaction functional 
\begin{equation*}
V(\om) = \frac{1}{2}\log(\frac{1}{\gamma}) \sum_{x \in \om} \one \left[0 < \abs{x}\leq r\right] = \frac{1}{2}\log(\frac{1}{\gamma}) (\om(b_r)-1\cdot \one \left[0 \in \om \right]).
\end{equation*} 
Since $\gamma \leq 1$, we get that $V$ is non-negative. Clearly it is also an $r$-local and increasing functional. Finally, we can bound $V(\om) \leq \frac{1}{2}\log(\frac{1}{\gamma}) \cdot \abs{\om} $ for $\om \in \Omega_f$, and therefore $V$ is cardinality bounded by the sequence $\{M_b\}_{b \in \N}$, where $M_b = \frac{1}{2}\log(\frac{1}{\gamma}) b$. Therefore, we know from Theorem \ref{thm:ldpunbounded} that the limit $\lim { \frac{1}{\abs{W_n}} } \log Z_n$ exists and is equal to 
\begin{equation*}
\lim { \frac{1}{\abs{W_n}} } \log Z_n = - \inf \{I(P) + P^o(V): P \in \mc{P}_s, P^o({\bf1}) = \la\}.
\end{equation*}
In order to prove the LDP for individual empirical fields of binomial Strauss process (and to be able to apply Corollaries \ref{cor:UnbUBUnbscore} and \ref{cor:UnbLBUnbscore}), we must also prove {that} $\lim { \frac{1}{\abs{W_n}} } \log Z_n > -\infty$. But this can be show as follows. Let $\Pi$ denote the distribution of the Poisson point process $\pi$ in $\R^d$ with intensity $\la$. Then, 
\begin{equation*}
     \inf \{I(P) + P^o(V): P \in \mc{P}_s, P^o({\bf1}) = \la\} \leq I(\Pi) + \Pi^o(V)
\end{equation*}
and $I(\Pi) = 0$. Furthermore, to finish the argument, we can bound 
\begin{align*}
  \Pi^o(V) &=  \la \E V(\pi \cup \{0\}) =  \la \E V((\pi\cap b_r) \cup \{0\})  \\ &=     \la  \sum_{l = 0}^\infty \E V((\pi\cap b_r) \cup \{0\}) \cdot \one[\abs{\pi(b_r)} = l]  \\
     &\leq    \la  \sum_{l = 0}^\infty  M_{l+1} \P[\pi(b_r) = l] = \la  \sum_{l = 0}^\infty  \frac{1}{2}\log(\frac{1}{\gamma})(l+1) \P[\pi(b_r) = l] < +\infty.
\end{align*}
\ene

 To conclude this example, we note that we deal with the hard-core case ($\gamma = 0$) separately in {Sections \ref{secresults:HC} and \ref{sec:HC}.} It is also worth noting that the binomial Strauss process is well defined also for $\gamma > 1$ (contrary to the standard Strauss process). However, our results do not carry over to this case, since the corresponding interaction function would not be bounded from below.

We now move to the model with $V$ driven by the interaction of $k$-tuples of points. Let $k \in \N$, $k \geq 2$.

\bee[Bounded $k$-wise interaction] \label{ex:k-int}
Let  $r > 0$ and let $\phi : (\R^d)^k \rightarrow \left[0,+\infty\right)$ be a bounded and symmetric measurable function. Define
\begin{equation}\label{def:k-wiseInt}
    V(\om) = \sum_{{\{x_1,\dots,x_{k-1}\}\subset \om \cap b_r:}\atop{ x_i \neq x_j,\, i \neq j} } \phi(x_1,\dots,x_{k-1},0),\, \om \in \Omega_f,
\end{equation}
with the convention that $V(\om) = 0$, if $\om(b_r) < k-1$. Then, $V$ is non-negative, $r$-local and increasing. Denote by $c$ the constant bounding the function $\phi$, then we can bound 
\begin{align*}
    V(\om) = \sum_{{\{x_1,\dots,x_{k-1}\}\subset \om \cap b_r:}\atop{ x_i \neq x_j,\, i \neq j} } \phi(x_1,\dots,x_{k-1},0) \leq c \abs{\om\cap b_r}^k \leq c \abs{\om}^k. 
\end{align*}
Therefore, $V$ is cardinality bounded by the sequence $\{cb^k\}_{b \in \N}$ and as it also holds that 
\begin{equation*}
    \sum_{l = 0}^\infty  c(l+1)^k \P[\pi(b_r) = l] < +\infty,
\end{equation*}
we get (similarly as in Example \ref{ex:Strauss}) that $\lim { \frac{1}{\abs{W_n}} } \log Z_n > -\infty$ for the process given by the bounded $k$-wise interaction (\ref{def:k-wiseInt}). 
\ene

{Finally}, we provide examples of unbounded score functions that fit in our setting.

\bee[Unbounded score functions]\label{ex:scores} Regarding the application of Corollaries \ref{cor:LBunboundedscore} and \ref{cor:UnbLBUnbscore}, we can note that any $\xi$ given by the {Formula \ref{def:k-wiseInt}} satisfies the necessary assumptions. This includes the clique counts and power-weighted edge lengths in random geometric
graphs, see Section 2.1 in \cite{HJT20LowerLDP} for more details.
\ene

\section{Proof of Theorem \ref{thm:ldp}}\label{sec:proof}

{The goal of this section is to prove Theorem \ref{thm:ldp}. An important step is to prove an LDP for random vectors $\onenJ{\bfXi}{n}(\rho_n)$, where $\oneJ{\bfxi}$ are bounded local functions (see Lemma \ref{cor:boundedLDP}).} 
{To do that, we first prove an unnormalized large-deviation upper bound.} As the event that there is exactly $n$ points corresponds to a closed set, this is just a direct application of Corollary 3.2 from \cite{GeorgiiLDP}.

\bel \label{lemma:aux1}
Let the interaction function $V$ be $r$-local and bounded.  Let $C \subset \mc{M}$ be closed. 
Then, the following inequality holds
\begin{align*}
\begin{split}
      \limsup { \frac{1}{\abs{W_n}} }\log \int &\exp(-H_n(\om)) \one \left[ \abs{\om} = n\right] \one \left[R^o_{n,\om} \in C\right]\d \Pi_n(\om) \\
    & \leq - \inf \{I(P) + P^o(V): P \in \mc{P}_s, P^o({\bf1}) = \la, P^o \in C\}.
    \end{split}
\end{align*}
\enl 

\bep
Using Corollary 3.2 from \cite{GeorgiiLDP}, we can write
\begin{align*}
\begin{split}
      \limsup &{ \frac{1}{\abs{W_n}} }\log \int \exp(-H_n(\om)) \one \left[ \abs{\om} = n\right] \one \left[R^o_{n,\om} \in C\right]\d \Pi_n(\om) \\
     = &\limsup{ { \frac{1}{\abs{W_n}} }\log\int \exp(-G_n(\om))\d \Pi_n(\om)} \leq  - \inf \{I(P) + G_{lsc}(P^o): P \in \mc{P}_s\},
    \end{split}
\end{align*}
where $G_n(\om) = \frac{n}{\la} G(R^o_{n,\om})$ with $G(m) = G_1(m) + G_2(m) + G_3(m)$ for $G_1(m):=m(V)$, $G_2(m):= \infty \cdot \one\left[m({\bf1}) \neq \la \right]$ and $ G_3(m):= \infty \cdot \one\left[m \in C^c \right].$ Since $V$ is local and bounded, $G_1$ is continuous.

Using Lemmas \ref{rem:lsc/uscproperties} and \ref{rem:lsc/uscproperties2}, it follows that $G$ is lsc and hence $G=G_{lsc}$. Therefore, 
$$\inf \{I(P) + G_{lsc}(P^o): P \in \mc{P}_s\} = \inf \{I(P) + P^o(V): P \in \mc{P}_s, P^o({\bf1}) = \la, P^o \in C\},$$
which finishes the proof.
\enp

{As a special case, we get the following unnormalized upper bound.}

\bec[Unnormalized upper bound] \label{proof:upbound}
Let $J \in \N$. Let $\oneJ{C} \subset \R^J$ be a closed set.  Assume that the interaction function $V$ and the score functions $\bfxi_1^J$ are $r$-local and bounded. Then, the following upper bound holds
\begin{align} \label{upperbound}
\begin{split}
     \limsup { \frac{1}{\abs{W_n}} }\log &\int \exp(-H_n(\om)) \one \left[ \abs{\om} = n\right] \one \left[\bfXi_{1,n}^J(\om) \in \mathbf{C}_1^J\right]\d \Pi_n(\om) \\ &\leq - \inf \{I(P) + P^o(V): P \in \mc{P}_s, P^o({\bf1}) = \la, P^o(\bfxi_1^J) \in \la \mathbf{C}_1^J\}.
\end{split}
\end{align}
\enc

\bep
We can write $\one \left[\bfXi_{1,n}^J(\om) \in \mathbf{C}_1^J\right] = \one \left[ R^o_{n,\om} \in  \{m \in \mc{M}: m(\bfxi_1^J) \in \la \mathbf{C}_1^J \} \right]$, using $\bfXi_{1,n}^J(\om) = \frac{1}{\la} R^o_{n,\om}(\bfxi_1^J) $. 
Since $\bfxi_1^J$ are bounded and local functions, we get that the set $\{m \in \mc{M}: m(\bfxi_1^J) \in \la \mathbf{C}_1^J \}$ is a closed subset of $\mc{M}$ and therefore (\ref{upperbound}) is just a special case of Lemma \ref{lemma:aux1}.
\enp

In order to apply Corollary 3.2 from \cite{GeorgiiLDP} in the proof of the lower bound, we must proceed more cautiously. Firstly, we prove that the binomial case can be approximated by the Poisson {case, where the number of points is} $\eps$-close to $n$ (bound (\ref{PoistoBinom}) in Lemma \ref{lemma:auxlb}). This allows us to derive an auxiliary unnormalized lower bound in Lemma \ref{auxlowerbound}, with the infimum ranging over a smaller open set.

\bel[] \label{lemma:auxlb}
Let $J \in \N$ and $r,c > 0$. Let $\oneJ{O},\oneJ{U} \subset \R^J$ be open sets such that there exists $\delta$ small enough so that $O_j \oplus b(0,5c\delta) \subset U_j$, $j \in \{1,\dots,J\}$. For this $\delta$, let $\eps \in \left(0,\min \Bigl\{\delta, \frac{\delta}{e^2 \la  r^d  v_d}\Bigr\}\right)$. Assume that the interaction function $V$ and the score functions $\oneJ{\bfxi}$ are $r$-local and bounded by the constant $c$. Then, for any $n$ large enough
\begin{align}
&\E e^{-H_n(\pi_n)} \one\left[\abs{\pi_n} \in (n \pm n\eps) \right]\one \left[\onenJ{\bfXi}{n}(\pi_n) \in \oneJ{O} \right] \label{PoistoBinom}  \\
    &\leq 2e^{3cn\eps}n\eps \left[ e^{2cn\delta} \E e^{-H_n(\pi_n)}  \one \left[\onenJ{\bfXi}{n}(\pi_n) \in \oneJ{U} \right] \one \left[ \abs{\pi_n} = n\right] + e^{3cn-\frac{n\delta}{4}\log \frac{\delta}{\eps \la  r^d  v_d}} \right]. \nonumber
\end{align}
\enl

\bel[Auxiliary unnormalized lower bound] \label{auxlowerbound}
Let $J \in \N$ and $r,c > 0$. Let $\oneJ{O},\oneJ{U} \subset \R^J$ be open sets such that for some $\Delta > 0$ we have $O_j \oplus b(0,\Delta) \subset U_j$, $j \in \{1,\dots,J\}$. Assume that the interaction function $V$ and the score functions $\oneJ{\bfxi}$ are $r$-local and bounded by the constant $c$. Then, the following lower bound holds
\begin{align} \label{weaklb}
\begin{split}
    \liminf { \frac{1}{\abs{W_n}} }\log\int &\exp(-H_n(\om)) \one \left[ \abs{\om} = n\right] \one \left[\onenJ{\bfXi}{n}(\om)  \in \oneJ{U} \right]\d \Pi_n(\om) \\
    & \geq - \inf \{I(P) + P^o(V): P \in \mc{P}_s, P^o({\bf1}) = \la, P^o(\oneJ{\bfxi}) \in \la \oneJ{O}\}.
    \end{split}
\end{align}
\enl

{The proofs of Lemmas \ref{lemma:auxlb} and \ref{auxlowerbound} are postponed to Section \ref{sec:AuxProof3}.} The last step is to prove the true unnormalized lower bound, i.e.\,(\ref{weaklb}) with $\oneJ{O}=\oneJ{U}$. Then, we can prove that ${ \frac{1}{\abs{W_n}} }\log Z_n$ has a limit as $n \rightarrow \infty$ and subsequently we prove the LDP for the random vectors $\onenJ{\bfXi}{n}(\rho_n)$.

\bel[]\label{cor:boundedLDP}
 Assume that the interaction function $V$ is $r$-local and bounded by a constant $c$. Then, the following limit exists
\begin{equation} \label{limitZn}
     \lim{{ \frac{1}{\abs{W_n}} }\log Z_n} = - \inf \{I(P) + P^o(V): P \in \mc{P}_s, P^o({\bf1}) = \la\} =: - A.
\end{equation} 
Let $J \in \N$. Let $\oneJ{C}\subset \R^J$ be a closed set and $\oneJ{U} \subset \R^J$ be an open set. Assume that the score functions $\oneJ{\bfxi}$ are $r$-local and bounded by the constant $c$. Then, 
\begin{align}
\begin{split}\label{upbound}
    \limsup&{{ \frac{1}{\abs{W_n}} }\log P_n(\onenJ{\bfXi}{n} \in \oneJ{C})} \\&\quad \leq - \inf \{I(P) + P^o(V): P \in \mc{P}_s, P^o({\bf1}) = \la, P^o(\oneJ{\bfxi}) \in \la \oneJ{C}\} + A, 
    \end{split}
\end{align} 
\begin{align}
\begin{split}\label{lowerbound}
   \liminf&{{ \frac{1}{\abs{W_n}} }\log P_n(\onenJ{\bfXi}{n} \in \oneJ{U})} \\&\quad \geq - \inf \{I(P) + P^o(V): P \in \mc{P}_s, P^o({\bf1}) =\la, P^o(\oneJ{\bfxi}) \in \la \oneJ{U}\} + A.
     \end{split}
\end{align}
\enl

{
\bere  
Since $I$ and $V$ are bounded from below, we get that $A > -\infty$. Furthermore, let $\Pi$ be the distribution of the Poisson point process in $ \R^d$ with intensity $\la$. Then, $A \leq I(\Pi) + \Pi^o(V) \leq 0 +\la c < \infty.$ Also, recalling the notation of Theorem \ref{thm:ldp}, we have that $ \inf_{m \in \mc{M}} \tilde{\mc{J}}(m) = A \in \R$. Particularly, the lower and upper bounds are always well defined (i.e.\,we do not get the expression "$\infty - \infty$").
\enre
}

\bep
{The first step} is to prove the unnormalized lower bound 
\begin{align}\label{lb}
\begin{split}
    \liminf { \frac{1}{\abs{W_n}} }\log&\int \exp(-H_n(\om)) \one \left[ \abs{\om} = n\right] \one \left[\onenJ{\bfXi}{n}(\om)  \in \oneJ{U} \right]\d \Pi_n(\om) \\
    & \geq - \inf \{I(P) + P^o(V): P \in \mc{P}_s, P^o({\bf1}) = \la, P^o(\oneJ{\bfxi}) \in \la \oneJ{U}\}.
    \end{split}
\end{align}

To show that (\ref{lb}) holds, one needs to show that 
\begin{equation*}
    \liminf { \frac{1}{\abs{W_n}} }\log\int \exp(-H_n(\om)) \one \left[ \abs{\om} = n\right]  \one \left[\onenJ{\bfXi}{n}(\om)  \in \oneJ{U} \right]\d \Pi_n(\om) \geq - (I(P) + P^o(V))
\end{equation*}
holds for any $P \in \mc{P}_s$ such that $P^o({\bf1}) = \la$ and $P^o(\oneJ{\bfxi}) \in \la \oneJ{U}$. But this is true {since, for} any such $P \in \mc{P}_s$ there exists $\alpha > 0$ such that $b(P^o(\xi_j),\alpha) \subset \la U_j$ for any $j \in \{1,\dots,J\}$. Now, it is enough to use Lemma \ref{auxlowerbound} with  $U_j$, $O_j=b(\frac{1}{\la}P^o(\xi_j),\frac{ \alpha}{2 \la })$, $j \in \{1,\dots,J\}$, and $\Delta = \frac{\alpha}{2 \la}$. 

Now, we can prove (\ref{limitZn}). Recall that $Z_n = \int \exp(-H_n(\om)) \one \left[ \abs{\om} = n\right] \d \Pi_n(\om)$. Using Corollary \ref{proof:upbound} and (\ref{lb}), with $J = 1$, $\xi_1 = \mathbf{1}$, $C_1 = \{1\}$ and $U_1 = (\frac{1}{2},\frac{3}{2})$, we can show that
\begin{align*}
  - \inf &\{I(P) + P^o(V): P \in \mc{P}_s, P^o({\bf1}) = \la\}\leq \liminf{{ \frac{1}{\abs{W_n}} }\log Z_n} \\
  &\qquad \leq \limsup{{ \frac{1}{\abs{W_n}} }\log Z_n} \leq - \inf \{I(P) + P^o(V): P \in \mc{P}_s, P^o({\bf1}) =\la\}.
\end{align*}
Therefore, (\ref{limitZn}) holds. Now, it is enough to realize that for any $D_j \in \mathcal{B}^d$, $j \in \{1,\dots,J\}$,
\begin{align}\label{finaleq}
\begin{split}
   &{ \frac{1}{\abs{W_n}} }\log P_n(\onenJ{\bfXi}{n} \in \oneJ{D}) \\
  &= { \frac{1}{\abs{W_n}} }\log \frac{1}{Z_n} + { \frac{1}{\abs{W_n}} }\log\int \exp(-H_n(\om)) \one \left[ \abs{\om} = n\right] \one \left[\onenJ{\bfXi}{n}(\om) \in\oneJ{D}\right]\d \Pi_n(\om). 
  \end{split}
\end{align}
Taking $D_j = C_j$, $j \in \{1,\dots,J\}$, $\limsup$ on both sides of (\ref{finaleq}) and using Corollary \ref{proof:upbound} proves (\ref{upbound}) and taking $D_j = U_j$, $j \in \{1,\dots,J\}$, and $\liminf$ on both sides of (\ref{finaleq}) and using (\ref{lb}) proves (\ref{lowerbound}). 
\enp

Now, we are ready to prove { Theorem \ref{thm:ldp}}. 

\bep[Proof of Theorem \ref{thm:ldp}]
First, we prove the upper bound. Let $C \subset \mc{M}$ be closed (in the $\tau_{\mc{L}^o}$ topology). Lemma \ref{cor:boundedLDP} tells us that there exists the limit $\lim{{ \frac{1}{\abs{W_n}} }\log Z_n}$ and it is clearly equal to $ - \inf_{m \in \mc{M}} \tilde{\mc{J}}(m)$. Using Lemma \ref{lemma:aux1}, we get that 
\begin{align*}
\begin{split}
    & \limsup{{ \frac{1}{\abs{W_n}} }\log \P(R^o_{n,\rho_n} \in C)} \\
    &= \lim{{ \frac{1}{\abs{W_n}} }\log \frac{1}{Z_n}} + \limsup { \frac{1}{\abs{W_n}} }\log \int \exp(-H_n(\om)) \one \left[ \abs{\om} = n\right] \one \left[R^o_{n,\om} \in C\right]\d \Pi_n(\om) \\ 
    &\leq \inf_{m \in \mc{M}} \tilde{\mc{J}}(m) - \inf \{I(P) + P^o(V): P \in \mc{P}_s, P^o({\bf1}) = \la, P^o \in C\}  = - \inf_{m \in C} \mc{J}(m).
    \end{split}
\end{align*}

Now, we prove the lower bound. We have proved in Lemma \ref{cor:boundedLDP} that for any $J \in \N$, any $\xi_1,\dots,\xi_J$ local and bounded score functions and any $U_1,\dots,U_J$ open subsets of $\R$ 
\begin{align*}
\begin{split}
    \liminf&{{ \frac{1}{\abs{W_n}} }\log P_n(\bfXi_{1,n}^J \in \mathbf{U}_1^J)} \\
    &\geq - \inf \{I(P) + P^o(V): P \in \mc{P}_s, P^o({\bf1}) =\la, P^o(\bfxi_1^J) \in \la \mathbf{U}_1^J\} + \inf_{m \in \mc{M}} \tilde{\mc{J}}(m).
    \end{split}
\end{align*}
In other words, for any open set $O =  \{m \in \mc{M}: m(\bfxi_1^J) \in \la \mathbf{U}_1^J\}$ we have that
\begin{align}
    \begin{split}\label{ineq}
        \liminf&{{ \frac{1}{\abs{W_n}} }\log \P\left(R^o_{n,\rho_n} \in O \right)} \geq - \inf_{m \in O} \mc{J}(m).
    \end{split}
\end{align}
Now, let $U \subset \mc{M}$ be a general open set (in the $\tau_{\mc{L}^o}$ topology) and take any $m \in U$. If $\mc{J}(m) < +\infty$, then {(see (\ref{rem:basis}))} there exist $J \in \N$, $\delta > 0$ and $\bfxi_1^J$ bounded local functions such that $U_{\bfxi_1^J;m} \subset U$. Since we can write $U_{\bfxi_1^J;m} =  \{\tilde{m} \in \mc{M}: \tilde{m}(\bfxi_1^J) \in \la \mathbf{U}_1^J\} $ for $\la U_j = (m(\xi_j)-\delta,m(\xi_j)+\delta)$, we get that 
\begin{align*} 
    \liminf{{ \frac{1}{\abs{W_n}} }\log \P\left(R^o_{n,\rho_n} \in U \right)} &\geq \liminf{{ \frac{1}{\abs{W_n}} }\log \P\left(R^o_{n,\rho_n} \in U_{\bfxi_1^J;m} \right)} \\
    &\overset{(\ref{ineq})}{\geq} - \inf_{\tilde{m} \in U_{\bfxi_1^J;m}} \mc{J}(\tilde{m}) \geq - \mc{J}(m).
\end{align*}
Clearly, if $\mc{J}(m) = +\infty$, then $\liminf{{ \frac{1}{\abs{W_n}} }\log \P\left(R^o_{n,\rho_n} \in U \right)} \geq  - \mc{J}(m)$. Therefore, 
\begin{equation*}
     \liminf{{ \frac{1}{\abs{W_n}} }\log \P\left(R^o_{n,\rho_n} \in U \right)} \geq \sup_{m \in U}(-\mc{J}(m)) = -\inf_{m \in U} \mc{J}(m),
\end{equation*}
which finishes the proof.
\enp

{
Before we move to the proofs of the auxiliary lemmas, let us prove one simple observation, that will allows us to transition our results between the languages of binomial and Poisson point processes.

\bel \label{remarklimitZnZntilde}
It holds that $\lim_{n \rightarrow \infty} \frac{1}{\abs{W_n}} \log \Pi_n(\abs{\om} = n) = 0$.
\enl

\bep
Use Stirling's formula. 
\enp

Particularly, this means that $ \lim{\frac{1}{\abs{W_n}}\log \tilde{Z}_n} = \lim{\frac{1}{\abs{W_n}}\log Z_n}$, if at least one of them exists, and that Corollary \ref{proof:upbound} and the bound (\ref{lb}) can be formulated also for integrals w.r.t.\,the distribution of the binomial point process.}

\subsection{Proofs of auxiliary lemmas}\label{sec:AuxProof3}

It remains to prove the two auxiliary lemmas. In the first proof, we need the following elementary fact. 
\begin{equation}\label{factPois}
    \text{Let } K_n \sim Pois(n)\text{. Then, } \max_{k \in \N_0}\P(K_n = k) =   \P(K_n = n). 
\end{equation}

\bep[Proof of Lemma \ref{lemma:auxlb}]
Take $n$ large enough so that $n\eps > 1$. Consider the representation for Poisson point process $\pi_n$ using $\{X_i\}_i$ iid uniformly distributed random variables in $W_n$ and $K_n \sim Pois(n)$, $K_n$ independent of $\{X_i\}_i$. Denote $\mathbb{X}_i^j = \{X_i,\dots, X_j\}$ and $\mathbb{X}_i^{(j)}$ its periodized version\footnote{For simplicity of notation we omit the dependence to $n$, however,we should keep in mind that for each $n$ we have a different sequence $\{X_i\}_i$. 
} in $W_n$.  We can split
\begin{align*}
&\int \exp(-H_n(\om)) \one\left[n-n\eps < \abs{\om} <n+n\eps \right] \one \left[\onenJ{\bfXi}{n}(\om) \in \oneJ{O}  \right] \d \Pi_n(\om) = J_+ + J_-
\end{align*}
where
\begin{align*}
J_+ &:= \E\exp(-H_n(\mathbb{X}_1^{K_n}))\one\left[K_n \in (n,n+n\eps) \right]\one \left[\onenJ{\bfXi}{n}(\mathbb{X}_1^{K_n}) \in  \oneJ{O}\right], \\
J_- &:= \E\exp(-H_n(\mathbb{X}_1^{K_n}))\one\left[K_n \in \left( n-n\eps,n \right] \right]\one \left[\onenJ{\bfXi}{n}(\mathbb{X}_1^{K_n}) \in  \oneJ{O} \right].
\end{align*}

In the next step, it will be crucial to bound the change caused by adding or deleting $n\eps$ points from the window $W_n$. To do that, we define random variables
\begin{equation*}
    \mc{Z}_{n,+} := \mathbb{X}_1^{n}\left(\bigcup_{j = n+1}^{\lfloor n+n\eps\rfloor}b^{(n)}(X_j, r)\right) \text{ and } \mc{Z}_{n,-} := \mathbb{X}_1^{\lceil n-n\eps \rceil}\left(\bigcup_{j = \lceil n-n\eps \rceil +1}^{n}b^{(n)}(X_j,r)\right).
\end{equation*}
  Then, it holds that
\begin{align}\label{def:Z1Z2}
\begin{split}
    &\mc{Z}_{n,+} \lvert \mathbb{X}_{n+1}^{\lfloor n+n\eps\rfloor} \quad \sim Bi(n,p_n(\mathbb{X}_{n+1}^{\lfloor n+n\eps\rfloor}))  \\
&\mc{Z}_{n,-} \lvert \mathbb{X}_{\lceil n-n\eps \rceil+1}^{n} \sim Bi(\lceil n-n\eps \rceil,p_n(\mathbb{X}_{\lceil n-n\eps \rceil+1}^{n})),
\end{split}
\end{align}
 where the conditional probabilities are given by 
 \begin{align*}
 p_n(\mathbb{X}_{n+1}^{\lfloor n+n\eps\rfloor}) &:= \frac{\la\abs{\bigcup_{j = n+1}^{\lfloor n+n\eps\rfloor}b^{(n)}(X_j, r)}}{n},\\ 
 p_n(\mathbb{X}_{\lceil n-n\eps \rceil+1}^{n}) &:= \frac{\la\abs{\bigcup_{j = \lceil n-n\eps \rceil +1}^{n}b^{(n)}(X_j,r)}}{n}.
 \end{align*}
 Denote $p_\eps:= \eps \la  r^d  v_d$, then clearly $p_n(\mathbb{X}_{n+1}^{\lfloor n+n\eps\rfloor}) \leq 
p_\eps$ and $p_n(\mathbb{X}_{\lceil n-n\eps \rceil+1}^{n}) \leq p_\eps$. 

Define also two enlargements of the sets $O_j$ for $j \in \{1,\dots,J\}$ 
\begin{equation*}
    O_{j,n,+}:= O_j \oplus b(0,\eps c + \frac{2c}{n}\mc{Z}_{n,+})\text{ and }O_{j,n,-}:= O_j \oplus b(0,3\eps c + \frac{2c}{n}\mc{Z}_{n,-}).
\end{equation*} 
Denote $\onenJ{O}{n,+}:= O_{1,n,+}\times \dots \times O_{J,n,+}$ and analogously $\onenJ{O}{n,-}$. Let $h\in \{V,\xi_1,\dots,\xi_J\}$. If $K_n \in (n,n+n\eps)$, then we get that 
\begin{itemize}
    \item[(i)] $\abs{\sum_{i = n+1}^{K_n}h(\mathbb{X}_1^{(K_n)}-X_i)} \leq cn\eps$,
    \item [(ii)] $\abs{\sum_{i = 1}^{n}h(\mathbb{X}_1^{(K_n)}-X_i)-h(\mathbb{X}_1^{(n)}-X_i)} \leq 2c\mathbb{X}_1^{n}\left(\bigcup_{j = n+1}^{K_n}b^{(n)}(X_j, r)\right) \leq 2c\mc{Z}_{n,+}$.
      \end{itemize}      
Particularly $\one \left[\Xi_{j,n}(\mathbb{X}_1^{K_n}) \in O_j\right] 
       \leq \one \left[\Xi_{j,n}(\mathbb{X}_1^{n}) \in O_{j,n,+} \right]$ holds for any $j \in \{1,\dots,J\}$ and therefore also
      \begin{itemize}
    \item[(iii)] $
        \one \left[\onenJ{\bfXi}{n}(\mathbb{X}_1^{K_n}) \in \oneJ{O}\right]  \leq \one \left[\onenJ{\bfXi}{n}(\mathbb{X}_1^{n}) \in \onenJ{O}{n,+} \right].
       $
   \end{itemize}    
Similarly for $K_n \in \left( n-n\eps,n \right]$ we have that 
\begin{itemize}
    \item[(iv)] $\abs{\sum_{i = K_n+1}^{n}h(\mathbb{X}_1^{(K_n)}-X_i)} \leq cn\eps$,
    \item[(v)] \begin{align*}
        \abs{\sum_{i = 1}^{K_n}h(\mathbb{X}_1^{(K_n)}-X_i)-h(\mathbb{X}_1^{(n)}-X_i)}
        &\leq 2c\mathbb{X}_1^{n}\Bigl(\bigcup_{j = \lceil n-n\eps \rceil +1}^{n}b^{(n)}(X_j, r)\Bigr) \\ &\leq 2cn\eps + 2c\mc{Z}_{n,-},
        \end{align*}
    \item[(vi)] $\one \left[\onenJ{\bfXi}{n}(\mathbb{X}_1^{K_n}) \in\oneJ{O} \right] \leq \one \left[\onenJ{\bfXi}{n}(\mathbb{X}_1^{n}) \in \onenJ{O}{n,-} \right].$
\end{itemize}
Furthermore, we get from (\ref{factPois}) that 
\begin{itemize}
    \item[(vii)] $\P(K_n \in (n,n+n\eps) ) \leq n\eps \P(K_n = n)$ and $\P(K_n \in \left(n-n\eps,n\right] ) \leq n\eps \P(K_n = n)$.
\end{itemize}

For our choice of $\delta$ and $\eps$, we get, by conditioning on $\mathbb{X}_{n+1}^{\lfloor n+n\eps\rfloor}$ for $\mc{Z}_{n,+}$ and on $\mathbb{X}_{\lceil n-n\eps \rceil+1}^{n}$ for $\mc{Z}_{n,-}$ and by using the concentration inequality for the tails of binomial distribution (Lemma 1.1 in \cite{PenroseBook}), that
\begin{align} \label{IneqForZ12}
    \begin{split}
         \P \left( \mc{Z}_{n,+} > n\delta \right) &\leq \exp\left(-\frac{n\delta}{2}\log \frac{\delta}{p_\eps}\right), \\
         \P \left( \mc{Z}_{n,-} > \lceil n-n\eps \rceil\delta \right) &\leq \exp\left(-\frac{\lceil n-n\eps \rceil\delta}{2}\log \frac{\delta}{p_\eps}\right) \leq \exp\left(-\frac{n\delta}{4}\log \frac{\delta}{p_\eps}\right).
    \end{split}
\end{align}

We get, using these bounds and the fact that $\mc{Z}_{n,\pm} \leq n\delta$ implies $O_{j,n,\pm} \subset U_j$, $j \in \{1,\dots,J\}$, that
\begin{align*}
    J_+ \overset{\text{(i-iii)}}{\leq} &e^{cn\eps}\E e^{-H_n(\mathbb{X}_1^{n})} e^{2c\mc{Z}_{n,+}} \one\left[K_n \in (n,n+n\eps)\right]  \one \left[\onenJ{\bfXi}{n}(\mathbb{X}_1^{n}) \in \onenJ{O}{n,+} \right] \\
    \overset{\text{(vii)}}{\leq} &e^{cn\eps}n\eps  \P(K_n =n)\E e^{-H_n(\mathbb{X}_1^{n})} e^{2c\mc{Z}_{n,+}}  \one \left[\onenJ{\bfXi}{n}(\mathbb{X}_1^{n}) \in \onenJ{O}{n,+} \right] \\
    \leq\, \, &e^{cn\eps}n\eps \P(K_n =n)  \Bigl[ e^{2cn\delta} \E e^{-H_n(\mathbb{X}_1^{n})}  \one \left[\onenJ{\bfXi}{n}(\mathbb{X}_1^{n}) \in\onenJ{O}{n,+} \right]  \one \left[ \mc{Z}_{n,+} \leq n\delta\right] \\
    & \qquad \qquad \qquad \qquad \qquad \qquad + e^{3cn}\E \one\left[ \mc{Z}_{n,+} > n\delta\right] \Bigr] \\
     \overset{(\ref{IneqForZ12})}{\leq} &e^{3cn\eps}n\eps \left[ e^{2cn\delta} \E e^{-H_n(\mathbb{X}_1^{n})}  \one \left[\onenJ{\bfXi}{n}({\mathbb{X}_1^{n}}) \in \oneJ{U} \right] \one \left[ K_n = n\right] + e^{3cn-\frac{n\delta}{2}\log \frac{\delta}{p_\eps}} \right]
\end{align*}
and 
\begin{align*}
    J_- \overset{\text{(iv-vi)}}{\leq} &e^{3cn\eps}\E e^{-H_n(\mathbb{X}_1^{n})} e^{2c\mc{Z}_{n,-}} \one\left[K_n \in \left( n-n\eps,n \right] \right]\one \left[\onenJ{\bfXi}{n}(\mathbb{X}_1^{n}) \in \onenJ{O}{n,-} \right] \\
    \overset{\text{(vii)}}{\leq} &e^{3cn\eps}n\eps\P(K_n =n)\E e^{-H_n(\mathbb{X}_1^{n})} e^{2c\mc{Z}_{n,-}}  \one \left[\onenJ{\bfXi}{n}(\mathbb{X}_1^{n}) \in\onenJ{O}{n,-}  \right]\\
    \leq\, \, &e^{3cn\eps}n\eps \P(K_n =n) \Bigl[ e^{2cn\delta} \E e^{-H_n(\mathbb{X}_1^{n})}  \one \left[\onenJ{\bfXi}{n}(\mathbb{X}_1^{n}) \in \onenJ{O}{n,-}  \right] \one \left[ \mc{Z}_{n,-} \leq n\delta\right] \\
    & \qquad \qquad \qquad \qquad \qquad \qquad + e^{3cn}\E \one\left[ \mc{Z}_{n,-} > n\delta\right] \Bigr] \\
      \overset{(\ref{IneqForZ12})}{\leq} &e^{3cn\eps}n\eps \left[ e^{2cn\delta} \E e^{-H_n(\mathbb{X}_1^{n})}  \one \left[\onenJ{\bfXi}{n}(\mathbb{X}_1^{n}) \in \oneJ{U} \right] \one \left[ K_n = n\right] + e^{3cn -\frac{n\delta}{4}\log \frac{\delta}{p_\eps}} \right].
\end{align*}
Altogether
\begin{align*}
    \begin{split}
    J_+ + J_- &= \E e^{-H_n(\mathbb{X}_1^{K_n})}\one\left[K_n \in (n\pm n\eps) \right]\one \left[\onenJ{\bfXi}{n}(\mathbb{X}_1^{K_n})) \in \oneJ{O}  \right]   \\
    &\leq 2e^{3cn\eps}n\eps \left[ e^{2cn\delta} \E e^{-H_n(\mathbb{X}_1^{n})}  \one \left[\onenJ{\bfXi}{n}(\mathbb{X}_1^{n}) \in \oneJ{U} \right]  \one \left[ K_n = n\right] + e^{3cn-\frac{n\delta}{4}\log \frac{\delta}{p_\eps}} \right].
    \end{split}
\end{align*}
\enp

\bep[Proof of Lemma \ref{auxlowerbound}]
If $M :=\inf \{I(P) + P^o(V): P \in \mc{P}_s, P^o({\bf1}) = \la, P^o(\oneJ{\bfxi}) \in \la \oneJ{O}\}$ satisfies $M = \infty$, then the claim holds trivially. Therefore, we assume from now on that\footnote{Recall that since $V$ is bounded, we always have $M > - \infty$.}  $M 
< \infty.$

 At first, choose $\delta $ small enough such that $O_j \oplus b(0,5c\delta) \subset U_j$ for each $j \in \{1,\dots,J\}$. For this $\delta$, denote $a(\delta) := \min \Bigl\{\delta, \frac{\delta}{e^2 \la  r^d v_d},\frac{1}{r^d  v_d}\exp{\left(-\frac{12c}{\delta}-16\la c - \frac{4M}{\delta} + \log{\frac{\delta}{\la}}\right)}\Bigr\}$ and choose $\eps \in \left(0,a(\delta)\right)$.

Let $G_n^\eps(\om) = \frac{n}{\la}G^\eps(R^o_{n,\om})$, where $G^\eps(m) := G_1(m) + G^\eps_2(m)$, $G_1(m) = m(V)$ and
\begin{equation*}
G^\eps_2(m):= \infty \cdot \one\left[m({\bf1}) \notin (\la-\la\eps,\la+\la\eps)\right] + \sum_{j=1}^J \infty\cdot \one\left[m(\xi_j)  \in (\la O_j)^c \right].
\end{equation*}
Since $V$ is local and bounded, $G_1$ is continuous and it follows from Lemmas \ref{rem:lsc/uscproperties} and \ref{rem:lsc/uscproperties2} that $G^\eps_2(m)$ is {usc}. Therefore, also $G^\eps$ is {usc} and it follows from Corollary 3.2 from \cite{GeorgiiLDP} that 
\begin{align}\label{epsbound}
\begin{split}
    &\liminf{{ \frac{1}{\abs{W_n}} }\log  \int \exp(-G_n^\eps(\om)) \d \Pi_n(\om)} \geq - \inf \{I(P) + G^\eps(P^o): P \in \mc{P}_s\} \\ &  = - \inf \{I(P) + P^o(V): P \in \mc{P}_s,P^o({\bf 1}) \in (\la\pm \la\eps), P^o(\oneJ{\bfxi})\in \la \oneJ{O} \} \geq - M.
\end{split}
\end{align}
We know from Lemma \ref{lemma:auxlb} that for $n$ large enough so that $n\eps > 1$ we have that
\begin{align}\label{auxbound}
\begin{split}
&\int \exp(-G_n^\eps(\om)) \d \Pi_n(\om) \\&= \E\exp(-H_n(\pi_n))\one\left[\abs{\pi_n} \in (n\pm n\eps) \right]\one \left[\onenJ{\bfXi}{n}(\pi_n) \in \oneJ{O} \right]   \\
    &\leq 2e^{3cn\eps}n\eps \Bigl[ e^{2cn\delta} \E e^{-H_n(\pi_n)}  \one \left[\onenJ{\bfXi}{n}(\pi_n) \in \oneJ{U},\abs{\pi_n} = n\right]   + e^{3cn-\frac{n\delta}{4}\log \frac{\delta}{\eps \la r^d v_d}} \Bigr].
    \end{split}
\end{align}
Taking $\liminf { \frac{1}{\abs{W_n}} } \log$ on both sides of (\ref{auxbound}) and using (\ref{epsbound}) we get that 
\begin{align}
    - M \,\,&\leq  \liminf{{ \frac{1}{\abs{W_n}} }\log  \E\exp(-H_n(\pi_n))\one\left[\abs{\pi_n} \in (n\pm n\eps) \right]\one \left[\onenJ{\bfXi}{n}(\pi_n) \in \oneJ{O}  \right]} \label{ineq1} \\
    \,\,&\leq 3c \la \eps  \nonumber\\
    +\liminf&{{ \frac{1}{\abs{W_n}} }\log \left[ e^{2cn\delta} \E e^{-H_n(\pi_n)}  \one \left[\onenJ{\bfXi}{n}(\pi_n) \in \oneJ{U} \right] \one \left[ \abs{\pi_n} = n\right] + e^{3cn-\frac{n\delta}{4}\log \frac{ \delta}{\eps  \la r^d  v_d}} \right] } \nonumber. 
\end{align}
Since $\eps$ {is} chosen so that $\eps < \frac{1}{r^d v_d} \exp{\left(-\frac{12c}{\delta}-16\la c - \frac{4M}{\delta} + \log{\frac{\delta}{\la}}\right)}$, it follows that
\begin{equation*}
    \lim { \frac{1}{\abs{W_n}} } \log e^{3cn-\frac{n\delta}{4}\log \frac{\delta}{\eps  \la r^d  v_d}}  =3c - \frac{\delta}{4}\log \frac{\delta}{\eps  \la r^d  v_d} < -4c \la \eps - M.
\end{equation*}
This together with (\ref{ineq1}) implies that
\begin{align*}
    &\liminf{{ \frac{1}{\abs{W_n}} }\log \left[ e^{2cn\delta} \E e^{-H_n(\pi_n)}  \one \left[\onenJ{\bfXi}{n}(\pi_n) \in \oneJ{U} \right] \one \left[ \abs{\pi_n} = n\right] + e^{3cn-\frac{n\delta}{4}\log \frac{\delta}{\eps  \la r^d  v_d}} \right] } \\
    & \qquad \qquad = \liminf{{ \frac{1}{\abs{W_n}} }\log  e^{2cn\delta} \E e^{-H_n(\pi_n)}  \one \left[\onenJ{\bfXi}{n}(\pi_n) \in \oneJ{U} \right] \one \left[ \abs{\pi_n} = n\right]},
\end{align*}
and therefore
\begin{align*}
    &- M \leq 3c\la\eps + \liminf{{ \frac{1}{\abs{W_n}} }\log  e^{2cn\delta} \E e^{-H_n(\pi_n)} \one \left[\onenJ{\bfXi}{n}(\pi_n) \in \oneJ{U} \right] \one \left[ \abs{\pi_n} = n\right]  } \\
    &= 3c\la\eps + 2c\delta + \liminf{ { \frac{1}{\abs{W_n}} }\log\int \exp(-H_n(\om)) \one \left[ \abs{\om} = n\right] \one \left[\onenJ{\bfXi}{n}(\om) \in \oneJ{U} \right]\d \Pi_n(\om)}.
\end{align*}
Now, taking $\eps \downarrow 0$ and then $\delta \downarrow 0$ finishes the proof. 
\enp

\section{Proof of Theorem \ref{thm:ldpunbounded}}\label{sec:proofunbounded}

In this section, we prove Theorem \ref{thm:ldpunbounded} as well as Corollaries \ref{cor:LBunboundedscore} and \ref{cor:UnbLBUnbscore}. The proofs rely on Lemmas \ref{lemma:UBunbounded2}, \ref{lemma:LBunboundedV} and  \ref{lemma:LBunboundedscore} presented in this section. The proof of Lemma \ref{lemma:UBunbounded2} is again easier and can be found in Section \ref{sec:proofupunbounded}, while the proofs of the lower bounds from Lemmas \ref{lemma:LBunboundedV} and  \ref{lemma:LBunboundedscore} are substantially more intricate and are left to Section \ref{sec:prooflowunbounded}. 

First, we present Lemmas \ref{lemma:UBunbounded2} and \ref{lemma:LBunboundedV}.

\bel[Unnormalized upper bounds, unbounded case] \label{lemma:UBunbounded2}
Let $J \in \N$. Let $C \subset \mc{M}$ be a closed set.  Assume that the interaction function $V$ is non-negative and {quasi local bounded}. Then, the following upper bounds hold
\begin{align}
\begin{split}\label{UBunbounded}
     \limsup { \frac{1}{\abs{W_n}} }\log\int &\exp(-H_n(\om)) \one \left[ \abs{\om} = n\right] \one \left[R^o_{n,\om} \in C\right]\d \Pi_n(\om)  \\ 
    &\leq - \inf \{I(P) + P^o(V): P \in \mc{P}_s, P^o({\bf1}) = \la, P^o \in C\}, 
\end{split}
\end{align}
\begin{align}\label{upboundZn}
\begin{split}
     \limsup { \frac{1}{\abs{W_n}} }\log\int \exp(-&H_n(\om)) \one \left[ \abs{\om} = n\right] \d \Pi_n(\om)  
\\& \quad \leq - \inf \{I(P) + P^o(V): P \in \mc{P}_s, P^o({\bf1}) = \la\}.
    \end{split}
\end{align}
\enl

\bel[Unnormalized lower bound, unbounded case] \label{lemma:LBunboundedV}
Let $J \in \N$, $\oneJ{U} \subset \R^J$ {be} an open set and let $r,c > 0$. Assume that $V$ is non-negative, $r$-local, increasing and cardinality-bounded by a sequence $\{M_b\}_{b \in \N}$ and assume that the score functions $\oneJ{\bfxi}$ are $r$-local and bounded by the constant $c$. Then, 
\begin{align}
\begin{split} \label{lowbound1}
     \liminf { \frac{1}{\abs{W_n}} }\log&\int \exp(-H_n(\om)) \one \left[\onenJ{\bfXi}{n}(\om) \in \oneJ{U} \right]\d \mathbb{B}_n(\om)  \\ 
    &\geq - \inf \{I(P) + P^o(V): P \in \mc{P}_s, P^o({\bf1}) = \la, P^o(\oneJ{\bfxi}) \in \la \oneJ{U}\}, 
\end{split}
\end{align}
\begin{align}
\begin{split}\label{lowbound2}
     \liminf{ \frac{1}{\abs{W_n}} }\log\int \exp(-H_n(\om))  &\d \mathbb{B}_n(\om) 
    \\ &\geq - \inf \{I(P) + P^o(V): P \in \mc{P}_s, P^o({\bf1}) = \la\}.
\end{split}
\end{align}
\enl

{ Now, we can prove Theorem \ref{thm:ldpunbounded}.}

\bep[Proof of Theorem \ref{thm:ldpunbounded}]
We will again w.l.o.g.\,assume that $V \geq 0$. If $V \geq c'$ for some constant $c' < 0$, then we can work with $V' = V - c'$.

{Using the upper bound (\ref{upboundZn}), the lower bound (\ref{lowbound2})  and Lemma \ref{remarklimitZnZntilde},} we get that the limit $ \lim { \frac{1}{\abs{W_n}} } \log Z_n$ exists and is equal to $- \inf \{I(P) + P^o(V): P \in \mc{P}_s, P^o({\bf1}) = \la\}$. 

{Assume that $\lim \frac{1}{\abs{W_n}} \log Z_n$} is finite. Let $C \subset \mc{M}$ be any closed set. Then, we get, using (\ref{UBunbounded}) from Lemma \ref{lemma:UBunbounded2}, that
\begin{align*}
\begin{split}
    &\limsup{{ \frac{1}{\abs{W_n}} }\log \P(R^o_{n,\rho_n} \in C)} \\
    & = \lim{{ \frac{1}{\abs{W_n}} }\log \frac{1}{Z_n}} + \limsup { \frac{1}{\abs{W_n}} }\log \int \exp(-H_n(\om)) \one \left[ \abs{\om} = n\right] \one \left[R^o_{n,\om} \in C\right]\d \Pi_n(\om) \\ 
    & \leq \inf_{m \in \mc{M}} \tilde{\mc{J}}(m) - \inf \{I(P) + P^o(V): P \in \mc{P}_s, P^o({\bf1}) = \la, P^o \in C\}  = - \inf_{m \in C} \mc{J}(m).
    \end{split}
\end{align*}
We also get, using (\ref{lowbound1}) from Lemma \ref{lemma:LBunboundedV}, that for any $J \in \N$, any $\bfxi_1^J$ local and bounded score functions and {any $\mathbf{U}_1^J$, with $U_i \subset \R$ open,}
\begin{align*}
\begin{split}
    \liminf&{{ \frac{1}{\abs{W_n}} }\log P_n(\bfXi_{1,n}^J \in \mathbf{U}_1^J)} \\ & \geq - \inf \{I(P) + P^o(V): P \in \mc{P}_s, P^o({\bf1}) =\la, P^o(\bfxi_1^J) \in \la \mathbf{U}_1^J\} + \inf_{m \in \mc{M}} \tilde{\mc{J}}(m).
    \end{split}
\end{align*}

Let $U \subset \mc{M}$ be any open set. Then, using the same arguments as in Theorem \ref{thm:ldp}, we can show that 
$
     \liminf{{ \frac{1}{\abs{W_n}} }\log \P\left(R^o_{n,\rho_n} \in U \right)} \geq \sup_{m \in U}(-\mc{J}(m)) = -\inf_{m \in U} \mc{J}(m),
$
and the proof is finished.
\enp

{ To prove Corollaries \ref{cor:LBunboundedscore} and \ref{cor:UnbLBUnbscore}, one needs the following result.}

\bel[Unnormalized lower bounds for unbounded score functions]\label{lemma:LBunboundedscore}
 Let $J \in \N$ and $\mathbf{a} \in \R^J$. Assume that the score functions $\oneJ{\bfxi}$ are non-negative, $r$-local, increasing and cardinality-bounded by a sequence $\{M_b\}_{b \in \N}$ and that the interaction function $V$ fits into one of the two following settings 
\begin{itemize}
    \item[(i)]  $V$ is non-negative, $r$-local and bounded by the constant $c$.
    \item[(ii)]  $V$ is non-negative, $r$-local, increasing and cardinality-bounded by sequence $\{M_b\}_{b \in \N}$.
\end{itemize} 
Then,
\begin{align}
\begin{split} \label{lb1}
     \liminf { \frac{1}{\abs{W_n}} }\log&\int \exp(-H_n(\om)) \one \left[\onenJ{\bfXi}{n}(\om) < \mathbf{a} \right]\d \mathbb{B}_n(\om)  \\ 
    &\geq - \inf \{I(P) + P^o(V): P \in \mc{P}_s, P^o({\bf1}) = \la, P^o(\oneJ{\bfxi}) < \la \mathbf{a}\}. \\
\end{split}
\end{align} 
\enl

\bep[Proof of Corollary \ref{cor:LBunboundedscore}]
If we assume that $V,\xi_1,\dots,\xi_J$ are all non-negative, then the statement follows from (\ref{limitZn}) {in} Lemma \ref{cor:boundedLDP} and (\ref{lb1}) {in} Lemma \ref{lemma:LBunboundedscore} (see Lemma \ref{remarklimitZnZntilde} for the transition from binomial to Poisson process). If $V$ is bounded but possibly negative, then we can work with $V' = V - c$, since both $V$ and $V'$ define the same binomial Gibbs point process. If $\xi_1,\dots,\xi_J$ are bounded from below by $\tilde{c}$, where $\tilde{c} < 0$, then work with $\xi_j' = \xi_j - \tilde{c}$ and $a_j' = a_j - \tilde{c}$.
\enp

\bep[Proof of Corollary \ref{cor:UnbLBUnbscore}]
As in Corollary \ref{cor:LBunboundedscore}, we can assume that $V,\xi_1,\dots,\xi_J$ are all non-negative. The claim then follows from (\ref{lb1}) from Lemma \ref{lemma:LBunboundedscore}.
\enp

\subsection{Proof of Lemma \ref{lemma:UBunbounded2}}
\label{sec:proofupunbounded}
In this section, {we study the unnormalized upper bound in the case with unbounded Hamiltonians} and particularly prove Lemma \ref{lemma:UBunbounded2}. The consequence of not having a bounded interaction function is that we are not able to derive the full upper LDP. We are able to derive an upper bound for the unnormalized version, particularly also for the partition function $Z_n$, however, to derive the full upper LDP, we would need to derive a lower bound for the partition function $Z_n$. For this, as we will see in the next section, our assumptions on $V$ must be stronger.

\bep[Proof of Lemma \ref{lemma:UBunbounded2}]
Clearly (\ref{upboundZn}) is a special case of (\ref{UBunbounded}). To prove the latter, write
\begin{align*}
   \int \exp(-H_n(\om)) \one \left[ \abs{\om} = n\right] \one \left[R^o_{n,\om} \in C\right]\d \Pi_n(\om) &= \int \exp(-G_n(\om))\d \Pi_n(\om) 
\end{align*}
for $G_n(\om) = \frac{n}{\la} G(R^o_{n,\om})$ with $G(m) = G_1(m) + G_2(m)$, where $G_1(m): = m(V)$ and $G_2(m) =  \infty \cdot \one\left[m({\bf1}) \neq\la\right] + \infty \cdot \one\left[m \in C^c \right].$  It holds that $G_2$ is {lsc}. Therefore, to finish the proof, it is enough to show that also $G_1$ is {lsc}. Then, we can use Corollary 3.2 from \cite{GeorgiiLDP} as in the previous cases.

To show that $G_1$ is {lsc}, it is enough to prove that the level sets $L_b = \{m \in \mc{M}: G_1(m) \leq b  \}$ are closed for all $b \in \left[-\infty, \infty \right)$. Since $V \geq 0$, we get that $L_b = \emptyset$ for $b < 0$. Let $b \geq 0$, then we can write 
\begin{align*}
    L_b &= \{m \in \mc{M}: G_1(m) \leq b  \} = \{m \in \mc{M}: m(V) \leq b \} = \bigcap_{r = 1}^\infty \{m \in \mc{M}: m(V_r) \leq b \}
\end{align*}
due to the monotone convergence theorem and our assumptions on $V$. Since $V_r$ are bounded and local, the sets $\{m \in \mc{M}: m(V_r) \leq b \}$ are closed and the proof is finished. 
\enp
We should also add that under the assumptions on $V$ from the previous lemma, we cannot show that $Z_n > 0$. However, the statement holds even if $Z_n = 0$, as then the left-hand sides of (\ref{UBunbounded}) and (\ref{upboundZn}) are $-\infty$.

\section{Lower bound in the unbounded case} \label{sec:prooflowunbounded}
In this section, we derive unnormalized lower bounds in the case with unbounded Hamiltonian {and unbounded score functions}. Particularly, we prove Lemmas \ref{lemma:LBunboundedV} and  \ref{lemma:LBunboundedscore}. Our proofs will follow the proof of Theorem 1.2 from \cite{HJT20LowerLDP} with the strategy of getting rid of the (very unlikely) $b$-dense points (i.e.\,points that have at least $b$ neighbors up to some fixed distance $r$, where $b$ is large). However, since we cannot delete them as in the Poisson case (we have {a} fixed number of points) we must move them to a suitable part of the window.

\subsection{Proof of Lemmas \ref{lemma:LBunboundedV} and  \ref{lemma:LBunboundedscore}}

{Consider the following definition of a $b$-dense point.}

\bede[$b$-dense points]
Let $n \in \N$, $\rho > 0$ 
and $\om \in \Omega_{W_n}$. Then, $x \in \om$ is called a \textbf{$b$-dense point w.r.t radius} $\rho$ in $\om$ , if $\abs{(\om^{(n)}-x)\cap b_\rho} \geq b$. We denote $N^\rho_{n,b}(\om) = \abs{ \{x \in \om: x \text{ is $b$-dense w.r.t radius $\rho$ in } \om\} } $.
\ende

Let us now prove an analogue of Lemma 3.1 from \cite{HJT20LowerLDP} that states that $b$-dense points are sufficiently rare, as $b$ and $n$ grow to infinity. There are two main differences in our and their setting, that need to be considered. Firstly, we will work with binomial, not {the} Poisson point process (which is easily overcome by using the standard relationship between these two processes), and secondly, we have periodic boundary condition. The main idea behind this prove is to bound the number of $b$-dense points by a sum of Poisson random variables and then use (similarly as in the proof of Poisson concentration inequalities) the exponential Markov inequality.

Recall that $\mathbb{B}_n$ denotes the distribution of the binomial point process $B_n$ in $W_n$ with $n$ points.

\bel[Rareness of $b$-dense points] \label{lemma:rarebdensepoints}
Let $\delta,\rho > 0$.  Then, 
\begin{equation*}
    \limsup_{b \rightarrow \infty}  \limsup_{n \rightarrow \infty} { \frac{1}{\abs{W_n}} } \log \mathbb{B}_n(N^\rho_{n,b} > \delta n) = -\infty. 
\end{equation*}
\enl

\bep 
Fix $\delta,\rho > 0$. {Recall that $w_n:=\left(\frac{n}{\la}\right)^{\frac{1}{d}}$ is the} length of the side of the window $W_n$. Fix $L > 3\rho$ and define for $n \in \N$ large enough (such that $w_n  > L$)
\begin{equation*}
a_n = \frac{1}{3} \frac{w_n}{\lfloor \frac{w_n}{L}\rfloor} \qquad \text{ and } \qquad c_n = 3^d \left\lfloor\frac{w_n}{L}\right\rfloor^d.  
\end{equation*}
Then, $a_n > \rho$ and we can\footnote{For the sake of simple notation we exclude the "right boundary" of $W_n$ in this proof, however as it has Lebesgue measure zero, this does not change anything.} split the window {$W_n = \left[-\frac{w_n}{2},\frac{w_n}{2}\right)^d= \bigl[-\frac{1}{2}\left(\frac{n}{\la}\right)^{\frac{1}{d}},\frac{1}{2}\left(\frac{n}{\la}\right)^{\frac{1}{d}}\bigr)^d$} into $c_n$ disjoint cubes with side length $a_n$,
\begin{equation*}
    W_n = \bigcup_{i = 1}^{c_n} Q_{a_n}(z_{i,n}),
\end{equation*}
where $Q_{a_n}(z) = z + Q_{a_n}$, $Q_{a_n} = \left[-\frac{a_n}{2},\frac{a_n}{2}\right)^d$ and $z_{i,n}\in  W_n $ are suitably chosen\footnote{In fact $z_{i,n} \in a_n \Z +e_n$, where $e_n = 0$ if $(c_n)^{1/d}$ is odd and $e_n = \frac{a_n}{2}$ if $(c_n)^{1/d}$ is even.}.

For $\om \in \Omega_{W_n}$ denote as $M_{i,n}(\om) := \om(Q_{a_n}(z_{i,n}))$ and $M'_{i,n}(\om) := \om(Q^{(n)}_{3a_n}(z_{i,n}))$ the number of points in the cubes of side length $a_n$ and $3a_n$, respectively, and centers $z_{i,n}$. Recall that $Q^{(n)}_{r}(z)$ denotes the periodic version of $Q_{r}(z)$ in $W_n$. 
Furthermore, define $M''_{n,b}(\om) := \sum_{i = 1}^{c_n} M_{i,n}(\om) \one \left[ M'_{i,n}(\om) > b\right]$. Then, $N^\rho_{n,b}(\om) \leq M''_{n,b}(\om)$.

Using the relation between $\mathbb{B}_n$ and $\Pi_n$, we can write 
\begin{align}\label{bound0}
\begin{split}
    \log \mathbb{B}_n(N^\rho_{n,b} > \delta n) &= \log \Pi_n(N^\rho_{n,b} > \delta n \, | \, \abs{\om} = n) \\ &\leq \log \Pi_n(N^\rho_{n,b} > \delta n) - \log\Pi_n(\abs{\om} = n).
    \end{split}
\end{align}
Let $t > 0$. By the exponential Markov inequality (recall that $\pi_n$ denotes the Poisson point process with distribution $\Pi_n$) 
\begin{equation} \label{bound1}
     \log \Pi_n(N^\rho_{n,b} > \delta n) \leq \log \Pi_n(M''_{n,b} > \delta n) \leq -\delta t n + \log \E e^{tM''_{n,b}(\pi_n)}.
\end{equation}

\begin{figure}
   \centering
   \includegraphics[width=0.85\linewidth]{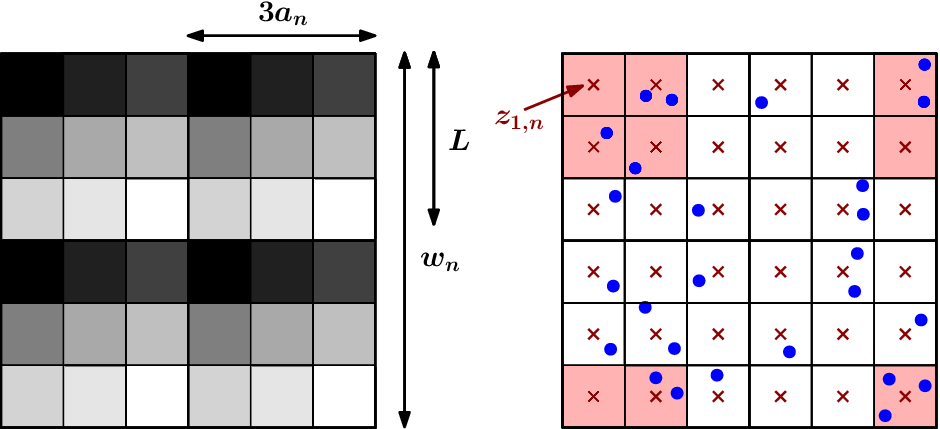}
  \caption{Graphic representation of the partition of the window $W_n$ from the proof of Lemma \ref{lemma:rarebdensepoints} in dimension $d = 2$. Left: The partition corresponding to the choice $\left\lfloor\frac{w_n}{L}\right\rfloor = 2$, i.e.\,$c_n = 3^2 \left\lfloor\frac{w_n}{L}\right\rfloor^2 =36$. Different shades of gray correspond to the partition of the indices $\{1,\dots,36\}$ into subsets $A_1,\dots, A_9$ of size $K_n =\left\lfloor\frac{w_n}{L}\right\rfloor^2= 4$. Right: A point configuration $\om$ (round points) in the partitioned window $W_n$ into cubes $Q_{a_n}(z_{1,n}),\dots, Q_{a_n}(z_{36,n})$. The crosses are the centers of the cubes $z_{i,n}$. Highlighted cubes correspond to the cube $Q_{3a_n}^{(n)}(z_{1,n})$. For this configuration we get that $M_{1,n}(\om) = 0$ and $M'_{1,n}(\om) = 11$.}
   \label{fig:Cubes}\end{figure}

Next, due to our choice of $a_n$, we can\footnote{We can simply take the regular subgrids of $a_n\Z + e_n$.} partition the indices $\{1,\dots, c_n\}$ into subsets $A_1, \dots A_{3^d}$ of the same size $K_n = \left\lfloor\frac{w_n}{L}\right\rfloor^d$ such that $M_{k,n}(\om) \one \left[ M'_{k,n}(\om) > b\right]$ and $M_{l,n}(\om) \one \left[ M'_{l,n}(\om) > b\right]$ are independent whenever $k,l \in A_j$, $k \neq l$ (see Figure~\ref{fig:Cubes}. for the graphical representation). 
Therefore, we can write $M''_{n,b}(\pi_n) = \sum_{j = 1}^{3^d} Z_{j,n}(\pi_n)$, where $Z_{j,n}(\pi_n) = \sum_{i \in A_j} M_{i,n}(\pi_n) \one \left[ M'_{i,n}(\pi_n) > b\right]$. The random variables $Z_{j,n}(\pi_n)$ have the same distribution, particularly they are sums of independent random variables distributed as $M_n \one \left[ M'_n > b\right]$, where $M_n \sim Pois(\la a_n^d)$ and $M'_n \sim Pois(\la (3a_n)^d)$.

Therefore, using Holder's inequality, we can write 
\begin{align}
    \begin{split}\label{bound2}
        \E e^{tM''_{n,b}(\pi_n)} &= \E \prod_{j = 1}^{3^d} e^{tZ_{j,n}(\pi_n)} \leq \prod_{j = 1}^{3^d} (\E e^{3^dtZ_{j,n}(\pi_n)})^{\frac{1}{3^d}}  = \E e^{3^dtZ_{1,n}(\pi_n)} \\
        &= \E \prod_{i \in A_1} e^{3^dtM_{i,n}(\pi_n) \one \left[ M'_{i,n}(\pi_n) > b\right]} = \left( \E e^{3^dtM_n \one \left[ M'_n > b\right]} \right)^{K_n}.
    \end{split}
\end{align}

Denote $M \sim Pois(\la L^d)$ and $M' \sim Pois(\la (3L)^d)$. For $n$ large enough we have that $a_n\leq L$ 
and therefore\footnote{This holds thanks to the fact that sum of two independent Poisson random variables is again Poisson random variable.}
\begin{equation}\label{bound3}
    \left( \E e^{3^dtM_n \one \left[ M'_n > b\right]} \right)^{K_n} \leq \left( \E e^{3^dtM \one \left[ M' > b\right]} \right)^{K_n}.
\end{equation}

Furthermore, since it holds (see Lemma \ref{remarklimitZnZntilde}) that $\lim { \frac{1}{\abs{W_n}} } \log\Pi_n(\abs{\om} = n)=0$, we can altogether bound
\begin{align*}
    \begin{split}
        &\limsup_{n \rightarrow \infty} { \frac{1}{\abs{W_n}} } \log \mathbb{B}_n(N^\rho_{n,b} > \delta n)  \\& \quad \overset{(\ref{bound0})}{\leq}  \limsup_{n \rightarrow \infty} { \frac{1}{\abs{W_n}} } \log \Pi_n(N^\rho_{n,b} > \delta n) - { \frac{1}{\abs{W_n}} }\log\Pi_n(\abs{\om} = n)\\
        & \quad \, = \limsup_{n \rightarrow \infty} { \frac{1}{\abs{W_n}} } \log \Pi_n(N^\rho_{n,b} > \delta n) \overset{(\ref{bound1})}{\leq} -\delta \la t  + \limsup_{n \rightarrow \infty}{ \frac{1}{\abs{W_n}} }\log \E e^{tM''_{n,b}(\pi_n)} \\
       &\,\overset{(\ref{bound2}),(\ref{bound3})}{\leq}  \limsup_{n \rightarrow \infty}{ \frac{1}{\abs{W_n}} }\log \left( \E e^{3^dtM \one \left[ M' > b\right]} \right)^{K_n} -\delta \la t    
        = \frac{\log \E e^{3^dtM \one \left[ M' > b\right]}}{ L^d} -\delta \la t .
    \end{split}
\end{align*}

Therefore, using the dominated convergence theorem, we get that for any $t > 0$
\begin{align*}
    \limsup_{b \rightarrow \infty}  \limsup_{n \rightarrow \infty} { \frac{1}{\abs{W_n}} } \log \mathbb{B}_n(N^\rho_{n,b} > \delta n) \leq -\delta \la t + \limsup_{b \rightarrow \infty} \frac{\log \E e^{3^dtM \one \left[ M' > b\right]}}{L^d}  = -\delta \la t . 
\end{align*}
Taking $t \uparrow \infty$ finishes the proof. 
\enp

Let $r$ denote the interaction range of the interaction function (also all the score functions in this section will later be assumed to be $r$-local). Fix $n_r,K_r\in \N$ such that 
\begin{equation}\label{def:n_r,K_r}
    n_r  > \la (18r)^d \text{ and } K_r > \left(\frac{6r}{\frac{1}{\la^{1/d}}-\frac{12r}{n_r^{1/d}}}\right)^d.
\end{equation}

Define for $n \geq n_r$
\begin{align}\label{def:S_n}
\begin{split}
    S_n(B_n) &:=\{z \in 6r\Z^d \cap W_n: \abs{Q_{6r}(z)\cap B_n} \leq K_r - 1,\, Q_{6r}(z) \subset W_n\} \\
    s_n(B_n) &:= \abs{S_n(B_n)}
    \end{split}
\end{align} the set of centers of cubes with side length $6r$ that lie inside $W_n$ and obtain at most $K_r - 1$ points from $B_n$, and the number of such cubes. As the number of points in these cubes is bounded, we will move the problematic $b$-dense points into some of them.

In the following lemma, we prove that there is a.s.\,at least $n A_r$ such cubes for some small $A_r$ depending on $r$ and $\la$, which ensures that unless we are in a rare event with many $b$-dense points, we will always have enough space to move them.

\bel[Number of sparse cubes] \label{lemma:numberofcubes}
Consider $S_n(B_n)$ and $s_n(B_n)$ defined in (\ref{def:S_n}), $n \geq n_r$. Denote 
\begin{equation}\label{def:A_r}
    A_r = \left( \left(\frac{\frac{1}{\la^{1/d}}-\frac{12r}{n_r^{1/d}}}{6r}\right)^d - \frac{1}{K_r}\right).
\end{equation}
Then, $A_r > 0$ and it holds that 
\begin{equation*}
    \P\left(s_n(B_n) \geq n A_r\right) = 1.
\end{equation*}
\enl

\bep 
Clearly $A_r > 0$ from our choice of $K_r$ and $n_r$. Denote by $$d_n = \abs{\{z \in 6r\Z^d \cap W_n: Q_{6r}(z) \subset W_n\}}$$ the number of cubes with side length $6r$ and centers in $6r\Z^d$, that are subsets of $W_n$. Then, it holds that 
\begin{align*}
    d_n = \left(2 \left \lfloor \frac{\left(\frac{n}{\la}\right)^{1/d} - 6r}{12r} \right \rfloor + 1 \right)^d &\geq \left( \frac{\left(\frac{n}{\la}\right)^{1/d}}{6r} - 2 \right)^d \\
    &\geq \frac{n}{\la}  \left(\frac{1-\frac{12r}{\left(\frac{n}{\la}\right)^{1/d}}}{6r}\right)^d \geq \frac{n}{\la}  \left(\frac{1-\frac{12r}{\left(\frac{n_r}{\la}\right)^{1/d}}}{6r}\right)^d 
\end{align*}
 There can be at most $\left \lfloor\frac{n}{K_r} \right \rfloor$ cubes with $K_r$ or more points. Therefore, 
\begin{equation*}
    s_n(B_n) \geq d_n - \left \lfloor\frac{n}{K_r} \right \rfloor \geq n \left(\frac{\frac{1}{\la^{1/d}}-\frac{12r}{n_r^{1/d}}}{6r}\right)^d -\frac{n}{K_r} = n  A_r,
\end{equation*}
and the proof is finished.
\enp
 
The key ingredient in the proof of the unnormalized lower bounds is the following coupling $(B_n,B_{n,\eps})$ of binomial point processes. Informally, $B_{n,\eps}$ is obtained from $B_n$ by randomly choosing some of its points and moving them independently to some other location.

\bede \label{def:couplingBPP} {
 Let $B_n$ be represented as $\{X_1,\dots, X_n\}$ iid, $X_1 \sim \text{Unif}(W_n)$. Let $\{X''_1,\dots, X''_n\}$ be iid, $X_1''\sim \text{Unif}(W_n)$ and $\{U_1,\dots, U_n\}$ iid, $U_1 \sim \text{Unif}\,(0,1)$, all of these jointly independent. Let $\eps > 0$ and define 
\[
X'_i := 
\begin{cases}
X_i  &\quad U_i \geq \eps, \\
X''_i &\quad U_i < \eps.
\end{cases}
\]
Then, $B_{n,\eps}$ represented by $\{X'_1,\dots, X'_n\}$ is also a binomial point process in $W_n$.} 
\ende

Furthermore, we define, for $b > 0$, the partition of $B_n$ into the set of $b$-dense and non $b$-dense points w.r.t.\,radius $r$, $B_{n,b}^\text{ds} := \{X_i \in B_n: X_i \text{ is $b$-dense w.r.t.\,radius $r$}\}$, and $B_{n,b}^\text{sp} := B_n \setminus B_{n,b}^\text{ds}$ and analogously $B_{n,\eps} := B_{n,\eps,b}^\text{ds}\cup B_{n,\eps,b}^\text{sp}$. For simplicity, we will write $j \in B_{n,b}^\text{ds}$ instead of $j : X_j \in B_{n,b}^\text{ds}$.

The next step is to define an event the has large enough probability (see Lemma \ref{lemma:boundforsetE}) and in which the moved points are those that have too many neighbors and they are moved into suitable subsets of the window (i.e., one of the sparse cubes).

\bede {
Denote $\mc{S}_n = \bigcup_{z \in S_n(B_n)} Q_r(z)$. 
For $n \geq n_r$, $\eps \in (0,1)$ and $b > 0$ define the event
\begin{equation*} 
    E_{n,b,\eps} = \{ B_{n,b}^\text{sp} = B_{n,\eps,b}^\text{sp} \} \cap \{\{X_i'\}_{i \in B_{n,b}^\text{ds}} \subset \mc{S}_n\} \cap \left\{\max_{z \in S_n(B_n)}\abs{Q_r(z) \cap\{X_i'\}_{i \in B_{n,b}^\text{ds}}} \leq 1\right\},
\end{equation*}
where all the $b$-dense points were moved into different $Q_{r}(z)$, $z \in S_n(B_n)$, and all the non $b$-dense points remained the same\footnote{ Here we ignore the event of probability 0 that we "moved" the non $b$-dense points into the same location they were before.}. Notice that if $K_r < b$, then $B_{n,\eps}$ has no $b$-dense points on $E_{n,b,\eps}$.}
\ende

In the next lemma, we derive a lower bound for the probability of $E_{n,b,\eps}$, conditionally on $B_n$. For simplicity, we will denote $N^{\rho}_{n,b}:= N^{\rho}_{n,b}(B_n)$ and 
$s_n:= s_n(B_n)$ in the following proofs.

\bel[Bound for $E_{n,b,\eps}$] \label{lemma:boundforsetE}
Let $n \geq n_r$, $\eps \in (0,1)$ and $b > 0$ . Then, we have that $\P$-a.s.
\begin{equation*}
    \one \left[N^r_{n,b} < s_n \right]\P(E_{n,b,\eps} \, | \, B_n ) \geq  \one \left[N^r_{n,b} < s_n \right]   (1-\eps)^{n} \eps^{N^r_{n,b}} \left(\frac{\la r^d(s_n - N^r_{n,b})}{ n}\right)^{N^r_{n,b}}.
\end{equation*}
\enl

\bep 
Assume that $N^r_{n,b} < s_n$, since the other case is trivial. Notice that $N^r_{n,b} = \abs{B_{n,b}^\text{ds}}$. Denote 
\begin{equation*}
D(\mathbb{X}_n'',B_n):=\bigcap_{i \in B_{n,b}^\text{ds}} \{X_i'' \in \mc{S}_n\}\cap \bigcap_{z \in S_n(B_n)} \{\abs{Q_r(z) \cap \{X_i''\}_{i \in B_{n,b}^\text{ds}}} \leq 1\}.
\end{equation*}  

From our construction of $B_n$ and $B_{n,\eps}$, we can write
\begin{align*}
\begin{split}
    \P(E_{n,b,\eps}\, | \, B_n ) &= \P\Bigl(\bigcap_{j \in B_{n,b}^\text{sp}}\{U_j \geq \eps\}\cap\bigcap_{j  \in B_{n,b}^\text{ds}}\{U_j < \eps\} \cap D(\mathbb{X}_n'',B_n) \, | \, B_n \Bigr) \\
   &= (1-\eps)^{n-N^r_{n,b}}\cdot \eps^{N^r_{n,b}} \cdot  \P\Bigl( D(\mathbb{X}_n'',B_n) \, | \, B_n \Bigr).
   \end{split}
\end{align*}
Furthermore
\begin{align*}
   \P\Bigl( &D(\mathbb{X}_n'',B_n) \, | \, B_n \Bigr) \\
   &= \P\Bigl( \bigcap_{i \in B_{n,b}^\text{ds}} \{X_i'' \in \mc{S}_n\}\cap \bigcap_{z \in S_n(B_n)} \{\abs{Q_r(z) \cap \{X_i''\}_{i \in B_{n,b}^\text{ds}}} \leq 1\} \, | \, B_n \Bigr)
   \\ &=\frac{s_n  r^d }{\frac{n}{\la}} \cdot \frac{(s_n-1)  r^d}{\frac{n}{\la}} \cdot \ldots \cdot \frac{(s_n-N^r_{n,b}+1)  r^d}{\frac{n}{\la}} 
   \geq \left(\frac{\la r^d(s_n - N^r_{n,b})}{n}\right)^{N^r_{n,b}},
\end{align*}
which finishes the proof.
\enp

Recall that $r$ denotes the interaction range of $V$, which is (among others) cardinality bounded by the sequence $\{M_b\}_{b \in \N}$. Recall also the definitions of $n_r,K_r$ and $A_r$ from (\ref{def:n_r,K_r}) and (\ref{def:A_r}). Before we move to the proof of the lower unnormalized bounds, we derive two auxiliary bounds that will allow us to bound (on the level of trajectories) the difference between models w.r.t binomial processes $B_n$ and $B_{n,\eps}$. The proof of this lemma is postponed to Section \ref{sec:ProofsAux}.

\bel[] \label{lemma:LBauxiliarybounds}
Let $c > 0$ and let $\eps \in (0,1)$. Let $b > K_r$ such that $M_b \geq c$. Let $h : \Omega \rightarrow \R$ be a measurable function and denote $h^M:=\min\{h,M\}$. Then, 
\begin{itemize}
    \item[(i)]if $h$ is $r$-local and bounded by the constant $c$, then  $\P$-a.s.
\begin{equation*}
    \abs{\sum_{i = 1}^n h(B_n^{(n)}-X_i) - h(B_{n,\eps}^{(n)}-X_i^{'})} \cdot \one_{E_{n,b,\eps}} \leq 2c(K_r + 1) N_{n,b}^{2r}\one_{E_{n,b,\eps}},   
\end{equation*}
    \item[(ii)] if $h$ is non-negative, $r$-local, increasing and cardinality-bounded by the sequence $\{M_b\}_{b \in \N}$, then  $\P$-a.s.
\begin{equation*}
     \one_{E_{n,b,\eps}}\sum_{i = 1}^n h^{M_b}(B_n^{(n)}-X_i) - h^{M_b}(B_{n,\eps}^{(n)}-X_i^{'})  \geq -M_{K_r} K_r N_{n,b}^{r}\one_{E_{n,b,\eps}}.  
\end{equation*}
\end{itemize}
\enl

{
For the rest of this section, we will use the following notation. For any $M > 0$  denote $V^M := \min\{V,M\}$ and $\xi^M_j:= \min\{\xi_j,M\}$ for any interaction function $V$ and score function $\xi_j$. Furthermore, denote 
\begin{align*}
    H^{M}_n(\om) &= \sum_{x \in \om} V^M(\om^{(n)}-x),\qquad
    \Xi^M_{j,n} (\om)= \frac{1}{n}\sum_{x \in \om} \xi_j^M(\om^{(n)}-x).
\end{align*}}

Now, we move to the proof of the unnormalized lower bounds. This proof has again two parts. Firstly, in Lemma \ref{lemma:Auxlbunbounded}, we approximate the case with unbounded Hamiltonian by the case with bounded Hamiltonian $H_n^{M_b}$ for $b$ large enough, using the coupling from Definition \ref{def:couplingBPP}, Lemma \ref{lemma:LBauxiliarybounds} and restricting ourselves to the set $E_{n,b,\eps}$.

\bel[]\label{lemma:Auxlbunbounded}
Let $J \in \N$, $\oneJ{O},\oneJ{U} \subset \R^J$ {be} open sets such that $O_j \oplus b(0,\Delta) \subset U_j$ for all $j$ and some $\Delta > 0$. Let $r,c > 0$, $\eps \in (0,1)$ and $\delta \in (0,\frac{\Delta}{2c(K_r+1)})$. Let $b > K_r$ and $n \geq n_r$. 
Assume that $V$ is non-negative, $r$-local, increasing and cardinality-bounded by the sequence $\{M_b\}_{b \in \N}$ and assume that the score functions $\oneJ{\bfxi}$ are $r$-local and bounded by the constant $c$. Then, 
\begin{align}\label{auxlowb}
\begin{split}
     &\E e^{-H_n(B_n)} \one \left[\onenJ{\bfXi}{n}(B_n) \in \oneJ{U}\right]  \\ 
        & \geq \E e^{-H^{M_b}_n(B_n)} e^{- M_{K_r}K_r\delta n } \one \left[\onenJ{\bfXi}{n}(B_n) \in\oneJ{O}\right] \one\left[N^{2r}_{n,b} \leq \delta n \right] \P(E_{n,b,\eps} \, | \, B_n ).
        \end{split}
\end{align}
\enl

The proof of this lemma is postponed to {Section \ref{sec:ProofsAux}. We are now ready to prove Lemma \ref{lemma:LBunboundedV} by taking the lower bound (\ref{auxlowb}), using the fact that we have restricted ourselves to a large enough event and the fact that we already have the LDP for the bounded Hamiltonian $H_n^{M_b}$.}

\bep[Proof of Lemma \ref{lemma:LBunboundedV}]
It is enough to prove (\ref{lowbound1}), since  (\ref{lowbound2}) follows from (\ref{lowbound1}). For any $M > 0$ we have that $V^M$ is a bounded and local function and we can\footnote{Recall Lemma \ref{remarklimitZnZntilde}. } therefore apply the results from Section \ref{sec:proof} to obtain that 
\begin{align}\label{lowb1}
\begin{split}
    \liminf { \frac{1}{\abs{W_n}} }\log&\int e^{-H^M_n(\om)} \one \left[\onenJ{\bfXi}{n}(\om) \in \oneJ{O}\right]\d \mathbb{B}_n(\om) \\
    &\geq - \inf \{I(P) + P^o(V^M): P \in \mc{P}_s, P^o({\bf1}) = \la, P^o(\oneJ{\bfxi}) \in \la \oneJ{O} \}  \\
    &\geq - \inf \{I(P) + P^o(V): P \in \mc{P}_s, P^o({\bf1}) = \la, P^o(\oneJ{\bfxi}) \in \la \oneJ{O}\} .
\end{split}
\end{align}

Similarly as in the bounded case, we first prove that for any $O_j, U_j$ open subsets of $\R$, $j \in \{1,\dots,J\}$, satisfying that there exist a $\Delta > 0$ such that $O_j \oplus b(0,\Delta) \subset U_j$, $j \in \{1,\dots,J\}$, we have that 
\begin{align}
\begin{split} \label{lowb2}
     \liminf { \frac{1}{\abs{W_n}} }\log\int &\exp(-H_n(\om)) \one \left[\onenJ{\bfXi}{n}(\om) \in \oneJ{U}\right]\d \mathbb{B}_n(\om)  \\ 
    &\geq - \inf \{I(P) + P^o(V): P \in \mc{P}_s, P^o({\bf1}) = \la, P^o(\oneJ{\bfxi}) \in \la \oneJ{O}\}. \\
\end{split}
\end{align}
The proof of (\ref{lowbound1}) then follows similarly as in Lemma \ref{cor:boundedLDP}.

Let $O_j \subset U_j$, $j \in \{1,\dots,J\}$ be open sets such that there exist a $\Delta > 0$ satisfying $O_j \oplus b(0,\Delta) \subset U_j$, $j \in \{1,\dots,J\}$.

If $L := \inf \{I(P) + P^o(V): P \in \mc{P}_s, P^o({\bf1}) = \la, P^o(\oneJ{\bfxi}) \in \la \oneJ{O}\} = \infty,$ then we are finished. Therefore, assume that $L < \infty$. Let $\eps \in (0,1)$ and $\delta \in (0,\min\{\frac{\Delta}{2c(K_r+1)},A_r\})$. For this $\eps$ and $\delta$ we can find, thanks to Lemma \ref{lemma:rarebdensepoints}, $b > K_r$ such that 
\begin{align} 
\begin{split}\label{choooseb}
    \la (-M_{K_r}K_r\delta + \log(1-\eps) + \delta \log \eps + \delta &\log(r^d\la(A_r-\delta))) -L  \\& >  \limsup_{n \rightarrow \infty} { \frac{1}{\abs{W_n}} } \log \mathbb{P}_n(N^{2r}_{n,b} > \delta n) .
\end{split}
\end{align}

Denote $f_{r,\delta} := \la r^d\left(A_r - \delta \right)$. Then, $0 < \la r^d\left(A_r - \delta \right)< 1$ and thanks to Lemma \ref{lemma:numberofcubes} we get that for our choice of $\delta$ and $n$, it holds that $\delta n < s_n$. Therefore, it follows from Lemmas \ref{lemma:boundforsetE} and \ref{lemma:numberofcubes} that 
\begin{align*}
\begin{split}
    \one\left[N^r_{n,b} \leq \delta n \right]  \P(E_{n,b,\eps} \, | \, B_n ) &\geq \one\left[N^r_{n,b} \leq \delta n \right]  (1-\eps)^{n} \eps^{N^r_{n,b}} \left(\frac{\la  r^d(s_n - N^r_{n,b})}{n}\right)^{N^r_{n,b}} \\
    &\geq \one\left[N^r_{n,b} \leq \delta n \right]  (1-\eps)^{n}  \eps^{n\delta}\left(\frac{\la r^d(s_n - n\delta)}{n}\right)^{N^r_{n,b}} \\ 
    &\geq \one\left[N^r_{n,b} \leq \delta n \right]  (1-\eps)^{n}  \eps^{n\delta}  f_{r,\delta}^{N^r_{n,b}} \\
    &\geq \one\left[N^r_{n,b} \leq \delta n \right]  (1-\eps)^{n}  \eps^{n\delta}  f_{r,\delta}^{n \delta}.
\end{split}
\end{align*}

 Since $N^r_{n,b} \leq N^{2r}_{n,b}$, particularly $\one\left[N^r_{n,b} \leq \delta n \right] \geq \one\left[N^{2r}_{n,b} \leq \delta n \right]$, it follows that
\begin{equation}\label{lowb3}
    \one\left[N^{2r}_{n,b} \leq \delta n \right]  \P(E_{n,b,\eps} \, | \, B_n ) \geq \one\left[N^{2r}_{n,b} \leq \delta n \right]  (1-\eps)^{n}  \eps^{n\delta}  f_{r,\delta}^{n \delta}.
\end{equation}
 Putting together (\ref{auxlowb}) and (\ref{lowb3}) with $V \geq 0$ and $f_{r,\delta} < 1 $, we get that
\begin{align}
        &\E \exp(-H_n(B_n)) \one \left[\onenJ{\bfXi}{n}(B_n) \in \oneJ{U}\right] \label{lowb4} \\
        &\geq \E e^{-H^{M_b}_n(B_n)} e^{-M_{K_r}K_r\delta n } \one \left[\onenJ{\bfXi}{n}(B_n) \in \oneJ{O}\right] \one\left[N^{2r}_{n,b} \leq \delta n \right]  (1-\eps)^{n} \eps^{n\delta}  f_{r,\delta}^{n \delta} \nonumber\\
        &\geq \left(\E e^{-H^{M_b}_n(B_n)} e^{- M_{K_r}K_r\delta n } \one \left[\onenJ{\bfXi}{n}(B_n) \in \oneJ{O} \right] (1-\eps)^{n}  \eps^{n\delta}  f_{r,\delta}^{n \delta} \right)- \P(N^{2r}_{n,b} > \delta n ). \nonumber
\end{align}

Now, we take $\liminf { \frac{1}{\abs{W_n}} } \log$ on both sides of (\ref{lowb4}), and put it together with (\ref{lowb1}). Thanks to our choice of $b$ (see (\ref{choooseb})) we get that
\begin{align*}
           &\liminf { \frac{1}{\abs{W_n}} } \log \E \exp(-H_n(B_n)) \one \left[\onenJ{\bfXi}{n}(B_n) \in \oneJ{U}\right] \\
           &\overset{(\ref{lowb4})}{\geq} \liminf { \frac{1}{\abs{W_n}} } \log \Bigl[ \E e^{-H^{M_b}_n(B_n)} e^{-  M_{K_r}K_r\delta n } \one \left[\onenJ{\bfXi}{n}(B_n) \in \oneJ{O}\right] (1-\eps)^{n}  \eps^{n\delta}  f_{r,\delta}^{n \delta} \\ & \hspace{10cm} - \P(N^{2r}_{n,b} > \delta n ) \Bigr] \\
           &\overset{(\ref{choooseb})}{=} \liminf { \frac{1}{\abs{W_n}} } \log \left[ \E e^{-H^{M_b}_n(B_n)} e^{-  M_{K_r}K_r\delta n } \one \left[\onenJ{\bfXi}{n}(B_n) \in \oneJ{O}\right] (1-\eps)^{n}  \eps^{n\delta}  f_{r,\delta}^{n \delta} \right]\\
           &\,\, = \la (- M_{K_r}K_r\delta + \log(1-\eps) + \delta \log\eps +  \delta \log f_{r,\delta}) \\ & \hspace{5cm} + \liminf { \frac{1}{\abs{W_n}} } \log \left[ \E e^{-H^{M_b}_n(B_n)} \one \left[\onenJ{\bfXi}{n}(B_n) \in \oneJ{O}\right] \right]\\
           &\overset{(\ref{lowb1})}{\geq} \la(- M_{K_r}K_r\delta + \log(1-\eps) + \delta \log(\eps \la r^d\left(A_r - \delta \right)))  - L.
\end{align*}
Taking first $\delta \downarrow 0$ and then $\eps \downarrow 0$ finishes the proof of (\ref{lowb2}). 
\enp

To derive the full LDP for the sequence of empirical measures driven by the processes $\rho_n$, we needed to assume that the score functions are bounded and local. However, using similar arguments as in the previous lemmas, we can {prove Lemma \ref{lemma:LBunboundedscore}, i.e., prove} unnormalized lower bounds for score functions that are non-negative, increasing, $r$-local and cardinality bounded. This allows us to derive additional large-deviation lower bounds for $\P_n(\onenJ{\bfXi}{n}< \mathbf{a})$, in models with both bounded and unbounded Hamiltonians. {First, we present an analogous result to Lemma \ref{lemma:Auxlbunbounded}.}

\bel[] \label{lemma:Auxlbunbounded2}
Let $J \in \N$ and  $0<\mathbf{a}_1<\mathbf{a}_2$, where $\mathbf{a}_i:= (a_{1,i},\dots,a_{J,i}) \in \R^J$. Let $n \geq n_r$. Let  $\eps \in (0,1)$, $\delta \in (0,\min_{j \in \{1,\dots,J\}} \frac{a_{j,2}-a_{j,1}}{M_{K_r}\cdot K_r})$ and  $b > K_r$.  

Assume that the score functions $\oneJ{\bfxi}$ are non-negative, $r$-local, increasing and cardinality-bounded by a sequence $\{M_b\}_{b \in \N}$ and the interaction function $V$ fits into one of the two following settings
\begin{itemize}
    \item[(i)]  $V$ is non-negative, $r$-local and bounded by a constant $c$.
    \item[(ii)]  $V$ is non-negative, $r$-local, increasing and cardinality-bounded by $\{M_b\}_{b \in \N}$.
\end{itemize}
Assume that $b$ is chosen large enough so that $M_b \geq c$. Then, 
\begin{align}\label{lbound2}
\begin{split}
     \E & e^{-H_n(B_n)} \one \left[\onenJ{\bfXi}{n}(B_n) < \mathbf{a}_2 \right]   \\ &\geq \E e^{-H^{M_b}_n(B_n)} e^{- \tilde{c}_r\delta n } \one \left[\bfXi^{J,M_b}_{1,n}(B_n) < \mathbf{a}_1\right] \one\left[N^{2r}_{n,b} \leq \delta n \right] \P(E_{n,b,\eps} \, | \, B_n ),
        \end{split}
\end{align}
where $\tilde{c}_r = \max\{2c(K_r+1),M_{K_r} K_r\}$ and $\bfXi^{J,M_b}_{1,n}:= (\Xi_{1,n}^{M_b},\dots,\Xi_{J,n}^{M_b})^T$.
\enl

 The proof Lemma \ref{lemma:Auxlbunbounded2} is again postponed to Section \ref{sec:ProofsAux}.

 {Now, we can prove Lemma \ref{lemma:LBunboundedscore}, i.e., prove the }full unnormalized lower bound for the case with unbounded score functions. As the proof is very similar to the proof of Lemma \ref{lemma:LBunboundedV}, we present only a sketch of it.

\bep[Proof of Lemma \ref{lemma:LBunboundedscore}]
The theorem holds trivially for $\mathbf{a} \leq 0$. 
Recall that for any $M > 0$ we have that $V^M$ and $\xi_j^M$ are bounded and local and we can\footnote{Recall Lemma \ref{remarklimitZnZntilde}. } therefore apply the results from Section \ref{sec:proof} to obtain that 
\begin{align*}
\begin{split}
    \liminf { \frac{1}{\abs{W_n}} }\log&\int e^{-H^M_n(\om)} \one \left[\bfXi^{J,M}_{1,n}(\om) < \mathbf{a}\right]\d \mathbb{B}_n(\om) \\
    &\geq - \inf \{I(P) + P^o(V^M): P \in \mc{P}_s, P^o({\bf1}) = \la , P^o(\bfxi^{J,M}_1) < \la \mathbf{a}\} \\
    &\geq - \inf \{I(P) + P^o(V^M): P \in \mc{P}_s, P^o({\bf1}) = \la, P^o(\oneJ{\bfxi}) < \la \mathbf{a} \}  \\
    &\geq - \inf \{I(P) + P^o(V): P \in \mc{P}_s, P^o({\bf1}) = \la, P^o(\oneJ{\bfxi}) < \la \mathbf{a}\} .
\end{split}
\end{align*}
Similarly as in the bounded case, we first prove that for any $0<\mathbf{a}_1<\mathbf{a}_2$ we have that 
\begin{align}
\begin{split} \label{lb3}
     \liminf { \frac{1}{\abs{W_n}} }\log\int &\exp(-H_n(\om)) \one \left[\onenJ{\bfXi}{n}(\om) < \mathbf{a}_2\right]\d \mathbb{B}_n(\om)  \\ 
    &\geq - \inf \{I(P) + P^o(V): P \in \mc{P}_s, P^o({\bf1}) = \la, P^o(\oneJ{\bfxi}) < \la \mathbf{a}_1\}. \\
\end{split}
\end{align}
The proof of (\ref{lb3}) follows from Lemmas \ref{lemma:numberofcubes} and \ref{lemma:boundforsetE} and (\ref{lbound2}) in the same way that (\ref{lowb2}) followed from Lemmas \ref{lemma:numberofcubes} and \ref{lemma:boundforsetE} and (\ref{auxlowb}) in the proof of Lemma \ref{lemma:LBunboundedV}. The proof of (\ref{lb1}) then follows from (\ref{lb3}) similarly as in Lemma \ref{cor:boundedLDP}. 
\enp

\subsection{Proofs of auxiliary lemmas}\label{sec:ProofsAux}

\bep[Proof of Lemma \ref{lemma:LBauxiliarybounds}] Assume that we are on $E_{n,b,\eps}$. First consider (i). Then, we can (see the justification below) bound
\begin{align*}
& \abs{\sum_{i = 1}^n h(B_n^{(n)}-X_i) - h(B_{n,\eps}^{(n)}-X_i^{'})}  \\
         &\leq \sum_{i \in B_{n,b}^\text{ds}} \abs{h(B_n^{(n)}-X_i) - h(B_{n,\eps}^{(n)}-X''_i)} + \sum_{i \in B_{n,b}^\text{sp}} \abs{h(B_n^{(n)}-X_i) - h(B_{n,\eps}^{(n)}-X_i)} \\
        &\leq 2cN^r_{n,b}  + \sum_{i \in B_{n,b}^\text{sp}: B_{n,\eps}^{(n)}\cap b(X_i,r) \not\subset B_{n}^{(n)}\cap b(X_i,r)} 2c +  \sum_{i \in B_{n,b}^\text{sp}: B_{n,\eps}^{(n)}\cap b(X_i,r) \subsetneqq B_{n}^{(n)}\cap b(X_i,r)}2c \\
        &\leq 2cN^{r}_{n,b} + 2cN^{r}_{n,b} + 2c(K_r-1)N^{2r}_{n,b} \leq 2c(K_r+1)N^{2r}_{n,b}.
    \end{align*}
In the third inequality we use the fact that if the new point $X''_j \in Q_r(z_j)$ interacts with some $X_i$, then (from the construction of $S_n$) $B(X_i,r) \subset Q_{6r}(z_j)$ and therefore also $X_i$ interacts with at most $K_r$ other points in $B_{n,\eps}$. Furthermore, each moved point can influence at most $K_r-1$ points. We also use the fact that if $X_i \in B_n$ such that it is a neighbor of a $b$-dense point w.r.t. radius $r$ in $B_n$ (and is therefore affected by the change between $B_n$ and $B_{n,\eps}$), then it is a $b$-dense point w.r.t. radius $2r$ in $B_n$. In other words, $\abs{\{X_i \in B_{n,b}^\text{sp}: B_{n,\eps}^{(n)}\cap b(X_i,r) \subsetneqq B_{n}^{(n)}\cap b(X_i,r) \}}\leq N^{2r}_{n,b}.$ Finally, we use the fact that $N^{r}_{n,b} \leq N^{2r}_{n,b}$.

Now, consider (ii). Then, we have the following bounds, to be justified below,
\begin{align*}
&\sum_{i = 1}^n h^{M_b}(B_n^{(n)}-X_i) - h^{M_b}(B_{n,\eps}^{(n)}-X_i^{'}) \\
         &= \sum_{i \in B_{n,b}^\text{ds}} h^{M_b}(B_n^{(n)}-X_i) - h^{M_b}(B_{n,\eps}^{(n)}-X''_i) \\
         & \hspace{2.5cm} + \sum_{i \in B_{n,b}^\text{sp}} h^{M_b}(B_n^{(n)}-X_i) - h^{M_b}(B_{n,\eps}^{(n)}-X_i) \\
        &\geq N^r_{n,b} \cdot (0-M_{K_r}) + 0 \\
        & \hspace{2.5cm} + \sum_{i \in B_{n,b}^\text{sp}: B_{n,\eps}^{(n)}\cap b(X_i,r) \not\subset B_{n}^{(n)}\cap b(X_i,r)} h^{M_b}(B_n^{(n)}-X_i) - h^{M_b}(B_{n,\eps}^{(n)}-X_i) \\
        &\geq -M_{K_r} N^r_{n,b} - (K_r - 1)M_{K_r}N^r_{n,b} = - M_{K_r} K_r N^r_{n,b}.
    \end{align*}
    In the first inequality, we use the fact that $h \geq 0$ and increasing and that all the moved (previously $b$-dense) points now interact with at most $K_r$ points. In the second inequality we again use the fact that if the new point $X''_j \in Q_r(z_j)$ interacts with some $X_i$, then (from the construction of $S_n$) $B(X_i,r) \subset Q_{6r}(z_j)$ and therefore also $X_i$ interacts with at most $K_r$ other points in $B_{n,\eps}$. Furthermore, each moved point can influence at most $K_r-1$ points.  
\enp

\bep[Proof of Lemma \ref{lemma:Auxlbunbounded}]
Consider the following observations
\begin{itemize}
    \item[(a)] $B_{n,\eps}$ has no $b$-dense points on $E_{n,b,\eps}$
    \item[(b)] on $E_{n,b,\eps}$ we can bound, using (ii) from Lemma \ref{lemma:LBauxiliarybounds} and the fact that $N^r_{n,b} \leq N^{2r}_{n,b}$ 
    \begin{align*}
        H_n^{M_b}(B_n) & - H_n^{M_b}(B_{n,\eps}) \geq - M_{K_r} K_r N^r_{n,b} \geq - M_{K_r} K_r N^{2r}_{n,b},
    \end{align*}
    \item[(c)] using (i) from Lemma \ref{lemma:LBauxiliarybounds} we get that for any $j\in \{1,\dots,J\}$ on $E_{n,b,\eps}$  
    \begin{equation*}
    \abs{\Xi_{j,n}(B_n) - \Xi_{j,n}(B_{n,\eps})} \leq \frac{2c(K_r + 1) N_{n,b}^{2r}}{n}
    \end{equation*} 
    and therefore we have, thanks to our choice of $\delta$, that
    \begin{equation*}
        \one\left[N^{2r}_{n,b} \leq \delta n \right] \one \left[ \onenJ{\bfXi}{n}(B_{n,\eps}) \in  \oneJ{U} \right] \geq \one\left[N^{2r}_{n,b} \leq \delta n \right] \one \left[\onenJ{\bfXi}{n}(B_n) \in  \oneJ{O}\right].
    \end{equation*} 
\end{itemize}
Therefore, we can bound
\begin{align*}
    \begin{split} 
        \E &e^{-H_n(B_n)} \one \left[\onenJ{\bfXi}{n}(B_n) \in \oneJ{U}\right]= \E e^{-H_n(B_{n,\eps})} \one \left[\onenJ{\bfXi}{n}(B_{n,\eps}) \in \oneJ{U}\right] \\
        &\geq \E e^{-H_n(B_{n,\eps})} \one \left[\onenJ{\bfXi}{n}(B_{n,\eps}) \in \oneJ{U}\right] \one_{E_{n,b,\eps}} \\
        &\overset{(a)}{=} \E e^{-H^{M_b}_n(B_{n,\eps})} \one \left[\onenJ{\bfXi}{n}(B_{n,\eps}) \in \oneJ{U}\right] \one_{E_{n,b,\eps}} 
        \\
        &\overset{(b)}{\geq} \E e^{-H^{M_b}_n(B_n)} e^{- M_{K_r} K_r N^{2r}_{n,b}} \one \left[\onenJ{\bfXi}{n}(B_{n,\eps}) \in \oneJ{U}\right] \one_{E_{n,b,\eps}} \\
        &\geq \E e^{-H^{M_b}_n(B_n)} e^{- M_{K_r} K_r N^{2r}_{n,b}} \one \left[\onenJ{\bfXi}{n}(B_{n,\eps}) \in \oneJ{U}\right] \one\left[N^{2r}_{n,b} \leq \delta n \right] \one_{E_{n,b,\eps}}  \\
        &\overset{(c)}{\geq} \E e^{-H^{M_b}_n(B_n)} e^{- M_{K_r} K_r\delta n } \one \left[\onenJ{\bfXi}{n}(B_n) \in  \oneJ{O} \right] \one\left[N^{2r}_{n,b} \leq \delta n \right] \one_{E_{n,b,\eps}}  \\
        & = \E e^{-H^{M_b}_n(B_n)} e^{- M_{K_r} K_r\delta n }\one \left[\onenJ{\bfXi}{n}(B_n) \in  \oneJ{O} \right] \one\left[N^{2r}_{n,b} \leq \delta n \right] \P(E_{n,b,\eps} \, | \, B_n ),
    \end{split}
\end{align*}
and the proof is finished.
\enp

\bep[Proof of Lemma \ref{lemma:Auxlbunbounded2}]
First, consider the setting (ii), i.e. $V\geq 0$, increasing, $r$-local and cardinality-bounded. Then, 
\begin{itemize}
    \item[(a)] $B_{n,\eps}$ has no $b$-dense points on $E_{n,b,\eps}${,}
    \item[(b)] on $E_{n,b,\eps}$ we can bound, using (ii) from Lemma \ref{lemma:LBauxiliarybounds}
    $$H_n^{M_b}(B_n)  - H_n^{M_b}(B_{n,\eps}) \geq - M_{K_r} K_r N^r_{n,b},$$ 
    \item[(c)] analogously $\bfXi^{J,M_b}_{1,n}(B_n) - \bfXi^{J,M_b}_{1,n}(B_{n,\eps}) \geq - \frac{1}{n}M_{K_r} K_r N^r_{n,b}$ on $E_{n,b,\eps}$  and therefore 
    $$
       \one_{E_{n,b,\eps}} \one \left[ \bfXi^{J,M_b}_{1,n}(B_{n,\eps}) < \mathbf{a}_2\right] \geq \one_{E_{n,b,\eps}} \one \left[\bfXi^{J,M_b}_{1,n}(B_n) < \mathbf{a}_2 - \frac{M_{K_r} K_r N^r_{n,b}}{n}\,  \right]{,}
    $$
    \item[(d)] thanks to our choice of $\delta$, it holds that $a_{j,2} - M_{K_r} K_r \delta \geq a_{j,1}$ for any $j\in \{1,\dots,J\}$,
    \item[(e)] clearly $N^r_{n,b} \leq N^{2r}_{n,b}$, particularly $\one\left[N^r_{n,b} \leq \delta n \right] \geq \one\left[N^{2r}_{n,b} \leq \delta n \right]$. 
\end{itemize}

Therefore, we can bound
\begin{align*}
    \begin{split} 
        &\E e^{-H_n(B_n)} \one \left[\onenJ{\bfXi}{n}(B_n) < \mathbf{a}_2\right]  = \E e^{-H_n(B_{n,\eps})} \one \left[\onenJ{\bfXi}{n}(B_{n,\eps}) < \mathbf{a}_2\right] \\
        &\geq \E e^{-H_n(B_{n,\eps})} \one \left[\onenJ{\bfXi}{n}(B_{n,\eps}) < \mathbf{a}_2\right] \one_{E_{n,b,\eps}} \\
        &\overset{(a)}{=} \E e^{-H^{M_b}_n(B_{n,\eps})} \one \left[\bfXi^{J,M_b}_{1,n}(B_{n,\eps}) < \mathbf{a}_2\right] \one_{E_{n,b,\eps}} \\
        &\overset{(b),(c)}{\geq}  \E e^{-H^{M_b}_n(B_n)} e^{- M_{K_r} K_r N^r_{n,b}} \one \left[\bfXi^{J,M_b}_{1,n}(B_n) < \mathbf{a}_2 - \frac{M_{K_r} K_r N^r_{n,b}}{n}\right] \one_{E_{n,b,\eps}} \\
        &\geq \E e^{-H^{M_b}_n(B_n)} e^{- M_{K_r} K_r N^r_{n,b}} \one \left[\bfXi^{J,M_b}_{1,n}(B_n) < \mathbf{a}_2 - \frac{M_{K_r} K_r N^r_{n,b}}{n}\right] \one\left[N^r_{n,b} \leq \delta n \right] \one_{E_{n,b,\eps}}  \\
        &\geq \E e^{-H^{M_b}_n(B_n)} e^{- M_{K_r} K_r\delta n } \one \left[\bfXi^{J,M_b}_{1,n}(B_n) < \mathbf{a}_2 - M_{K_r} K_r \delta \right] \one\left[N^r_{n,b} \leq \delta n \right] \one_{E_{n,b,\eps}}  \\
        & = \E e^{-H^{M_b}_n(B_n)- M_{K_r} K_r\delta n } \one \left[\bfXi^{J,M_b}_{1,n}(B_n) < \mathbf{a}_2 - M_{K_r} K_r \delta \right] \one\left[N^r_{n,b} \leq \delta n \right] \P(E_{n,b,\eps} \, | \, B_n )\\
        & \overset{(d)}{\geq} \E e^{-H^{M_b}_n(B_n)} e^{- M_{K_r} K_r\delta n } \one \left[\bfXi^{J,M_b}_{1,n}(B_n) < \mathbf{a}_1\right] \one\left[N^r_{n,b} \leq \delta n \right] \P(E_{n,b,\eps} \, | \, B_n )\\
        & \overset{(e)}{\geq} \E e^{-H^{M_b}_n(B_n)} e^{- \tilde{c}_r\delta n } \one \left[\bfXi^{J,M_b}_{1,n}(B_n) < \mathbf{a}_1\right] \one\left[N^{2r}_{n,b} \leq \delta n \right] \P(E_{n,b,\eps} \, | \, B_n ).
    \end{split}
\end{align*}

Now, assume that the setting (i) holds. Then, (a),(c),(d) and (e) still hold and furthermore we get that
\begin{itemize}
    \item[(f)] using (i) from Lemma \ref{lemma:LBauxiliarybounds} we can bound
    \begin{align*}
    \one_{E_{n,b,\eps}}(H_n(B_n) & - H_n(B_{n,\eps}))  \geq -2c(K_r+1)N^{2r}_{n,b}\one_{E_{n,b,\eps}} \geq - \tilde{c}_r N^{2r}_{n,b}\one_{E_{n,b,\eps}}. 
\end{align*}
\end{itemize}

Notice that in this setting $H_n^{M_b} = H_n$ from the choice of $b$. Therefore, it holds that
\begin{align*}
    \begin{split} 
        \E & e^{-H_n(B_n)} \one \left[\onenJ{\bfXi}{n}(B_n) < \mathbf{a}_2\right] 
        \overset{(a)}{\geq} \E e^{-H_n(B_{n,\eps})} \one \left[\bfXi^{J,M_b}_{1,n}(B_{n,\eps}) < \mathbf{a}_2\right] \one_{E_{n,b,\eps}} \\
        &\overset{(c),(f)}{\geq} \E e^{-H_n(B_n)} e^{- \tilde{c}_r N^{2r}_{n,b}} \one \left[\bfXi^{J,M_b}_{1,n}(B_n) < \mathbf{a}_2 - \frac{M_{K_r} K_r N^r_{n,b}}{n}\right] \one_{E_{n,b,\eps}} \\
        &\geq \E e^{-H_n(B_n)} e^{- \tilde{c}_r N^{2r}_{n,b}} \one \left[\bfXi^{J,M_b}_{1,n}(B_n) < \mathbf{a}_2 - \frac{M_{K_r} K_r N^r_{n,b}}{n}\right] \one\left[N^{2r}_{n,b} \leq \delta n \right] \one_{E_{n,b,\eps}}  \\
        &\overset{(e)}{\geq} \E e^{-H_n(B_n)} e^{-\tilde{c}_r\delta n } \one \left[\bfXi^{J,M_b}_{1,n}(B_n) < \mathbf{a}_2 - M_{K_r} K_r \delta \right] \one\left[N^{2r}_{n,b} \leq \delta n \right] \one_{E_{n,b,\eps}}  \\
        & = \E e^{-H_n(B_n)} e^{- \tilde{c}_r\delta n } \one \left[\bfXi^{J,M_b}_{1,n}(B_n) < \mathbf{a}_2 - M_{K_r} K_r \delta \right] \one\left[N^{2r}_{n,b} \leq \delta n \right] \P(E_{n,b,\eps} \, | \, B_n )\\
        & \overset{(d)}{\geq} \E e^{-H_n(B_n)} e^{-\tilde{c}_r\delta n } \one \left[\bfXi^{J,M_b}_{1,n}(B_n) < \mathbf{a}_1\right] \one\left[N^{2r}_{n,b} \leq \delta n \right] \P(E_{n,b,\eps} \, | \, B_n )\\
        & = \E e^{-H_n^{M_b}(B_n)} e^{- \tilde{c}_r\delta n } \one \left[\bfXi^{J,M_b}_{1,n}(B_n) < \mathbf{a}_1\right] \one\left[N^{2r}_{n,b} \leq \delta n \right] \P(E_{n,b,\eps} \, | \, B_n ).
    \end{split}
\end{align*}
\enp

\section{Proof of Theorem \ref{cor:HCLDP}} \label{sec:HC}
We start by proving the unnormalized upper bound for the hard-core process. This is again just a straightforward application of the theory from \cite{GeorgiiLDP}.

\bel[Unnormalized upper bound for hard-core interaction]\label{lemma:HCAuxUB}
Let $C \subset \mc{M}$ be closed. Then, the following inequality holds
\begin{align*}
\begin{split}
      \limsup { \frac{1}{\abs{W_n}} }\log \int &\exp(-H_n(\om)) \one \left[ \abs{\om} = n\right] \one \left[R^o_{n,\om} \in C\right]\d \Pi_n(\om) \\
    & \leq - \inf \{I(P) + P^o(V): P \in \mc{P}_s, P^o({\bf1}) = \la, P^o \in C\}.
    \end{split}
\end{align*}
\enl 

\bep
As in Lemma \ref{lemma:aux1}, we can write
\begin{align*}
      \limsup { \frac{1}{\abs{W_n}} }\log \int  \exp(-H_n(\om)) \one \left[ \abs{\om} = n\right] &\one \left[R^o_{n,\om} \in C\right]\d \Pi_n(\om) 
     \\ &\leq  - \inf \{I(P) + G_{lsc}(P^o): P \in \mc{P}_s\},
\end{align*}
with $ G(m) = G_1(m) + G_2(m) + G_3(m):=m(V) + \infty \cdot \one\left[m({\bf1}) \neq \la \right] + \infty \cdot \one\left[m \in C^c \right].$ Now, it is enough to prove that $G_1$ is a {lsc} function. 

Let $b \in \R$. Then, for $ b < 0$ we have that $L_b := \{m \in \mc{M}: G_1(m) \leq b  \} = \emptyset$ is a closed set. Let $b \geq 0$, then
$$
L_b = \{m \in \mc{M}: G_1(m) \leq b  \} = \{m \in \mc{M}: m(V) = 0 \} =  \{m \in \mc{M}: m(V^1) = 0\},
$$ 
which is again closed as $V^1$ is a local bounded function. Altogether, all level sets $L_b$ are closed and therefore $G_1$ is a {lsc} function and the proof is finished, since both $G_2$ and $G_3$ are also {lsc}. 
\enp

{ Now, we state the unnormalized lower bound. As before, the proof of the lower bound will be substantially more demanding and is postponed to Section \ref{sec:HCLBound}.}

\bel[Unnormalized lower bound] \label{lemma:HCUnnormLB}
Let $J \in \N$ and $ c,r > 0$. Let $\oneJ{U} \subset \R^J$ be an open set. Assume that the score functions $\oneJ{\bfxi}$ are bounded by the constant $c$ and $r$-local. Then, it holds that 
\begin{align}
\begin{split}\label{HCUnnormLB}
    &\liminf{ { \frac{1}{\abs{W_n}} }\log\int \exp(-H_n(\om)) \one \left[ \abs{\om} = n\right] \one \left[\onenJ{\bfXi}{n}(\om)  \in \oneJ{U}\right]\d \Pi_n(\om)} \\
    & \qquad  \qquad  \geq - \inf \{I(P) + P^o(V): P \in \mc{P}_s, P^o({\bf1}) = \la, P^o(\oneJ{\bfxi}) \in \la \oneJ{U}\}.
    \end{split}
\end{align}
\enl

\bep[Proof of Theorem \ref{cor:HCLDP}]
The claim follows easily from Lemmas \ref{lemma:HCAuxUB} and \ref{lemma:HCUnnormLB}. 
\enp

\subsection{Proof of Lemma \ref{lemma:HCUnnormLB}} \label{sec:HCLBound}
The proof of the lower bound follows a similar pattern as the proof in Section \ref{sec:prooflowunbounded}, particularly we define a suitable coupling of point processes (see Definition \ref{def:Coupling}) and an event on which we get rid of the problematic points.

The difference is that now, a problematic point is any point that has one or more $R$-neighbors. Therefore, we approximate $V$ by a local bounded interaction function $V^s$, $s \geq 0$, defined as
\begin{equation}\label{def:Vs}
    V^s(\om) =  s\cdot \one \left[\abs{\om \cap b(0,R)} \geq 2 \right].
\end{equation}

The corresponding Hamiltonian is denoted as $ H^s_n(\om) = \sum_{x \in \om\cap W_n} V^s(\om^{(n)}-x)$. 
We can prove that under the Hamiltonian $H_n^s$, the problematic points are sufficiently rare as $s \rightarrow \infty$.

\bel[] \label{lemma:chooseS}
Let $\delta > 0$. Then,  
\begin{equation*}
    \limsup_{s \rightarrow \infty} \limsup_{n \rightarrow \infty} { \frac{1}{\abs{W_n}} } \log \E \left[ \exp{(-H^s_n(\pi_n))} \cdot \one\left[N_n(\pi_n) > n\delta\right] \right] = -\infty.
\end{equation*}
\enl

\bep
First,
$\E \left[ e^{-H^s_n(\pi_n)} \cdot \one\left[N_n(\pi_n) > n\delta\right] \right] =  \E \left[ e^{-sN_n(\pi_n) } \cdot \one\left[N_n(\pi_n) > n\delta\right] \right] \leq e^{-n s\delta}.$
Therefore,
\begin{align*}
    \limsup_{s \rightarrow \infty} \limsup_{n \rightarrow \infty} { \frac{1}{\abs{W_n}} } \log \E \left[ \exp{(-H^s_n(\pi_n))} \cdot \one\left[N_n(\pi_n) > n\delta\right] \right] \leq \limsup_{s \rightarrow \infty} (-s \delta \la )  = -\infty.
\end{align*}
\enp

We now prove the unnormalized lower bound. Our strategy is the following. 

\begin{enumerate}
    \item First we prove the unnormalized lower bound for the Poisson case (see Lemma \ref{lemma:HCAuxLBI}), approximating our unbounded interaction $V$ by its bounded counterpart $V^s$ (recall (\ref{def:Vs})) for $s$ large enough and using a concentration inequality (see Lemma \ref{lemma:AuxConcBound}) for the number of $r$-neighbors of points that have at least one $R$-neighbor.
    \item Then, we convert the problem from the case of the binomial point process to the case with Poisson point process, while controlling the price of this change by restricting ourselves to having the number of points being $\eps$-close to $n$ (see Lemma \ref{lemma:HCAuxLBII}). 
\end{enumerate}

The proofs of Lemmas \ref{lemma:AuxConcBound}, \ref{lemma:HCAuxLBI} and \ref{lemma:HCAuxLBII} are postponed to the end of this section. These three lemmas are the key results to prove the unnormalized lower bound in Lemma \ref{lemma:HCUnnormLB}. The full LDP for local bounded functionals of the hard-core process is summarized in Theorem \ref{cor:HCLDP}.

Recall that $R$ denotes the radius of the hard-core interaction. Let us stress that it is crucial (opposite to the case from Section \ref{sec:proof}) to distinguish the interaction radius $R$ of the interaction function $V$ and the interaction radius $r$ of the score functions, as $R$ is tied to the intensity $\la$ by the bound (\ref{ass:HCmodelLambda}), whereas the interaction radius of the score functions does not have to satisfy such requirement.

We will work with the following coupling of Poisson point processes.

\bede \label{def:Coupling}{
Let $ \delta \in (0,1)$ and $n \in \N$. Let $\mathbb{U},\mathbb{X},K_n,\mathbb{Y},M_n^\delta$ be jointly independent such that
\begin{enumerate}
    \item $\mathbb{U} = \{U_1,U_2,\dots\}$ is a sequence of iid Unif\,$([0,1])$ random variables,
    \item $\mathbb{X}, \mathbb{Y}$ are sequences of iid Unif\,$(W_n)$ random vectors,
    \item $K_n \sim Pois(n)$ and $M_n^\delta \sim Pois(\delta n)$.
\end{enumerate}
Then, $\mathbb{X}_1^{K_n}:= \{X_1,\dots, X_{K_n}\}$ is a Poisson point process in $W_n$ with intensity $\la$ and $\mathbb{Y}_1^{M_n^\delta}$ is a Poisson point process in $W_n$ with intensity $\delta \la$. Denote by $\text{thin}_\delta(\mathbb{X}_1^{K_n})$ the process given by thinning $\mathbb{X}_1^{K_n}$ based on the rule
\begin{equation*}
    X_i \text{ is deleted } \iff U_i < \delta.
\end{equation*}
Then, the process $\mathbb{Z}_1^{L_n}
:= \text{thin}_\delta(\mathbb{X}) \cup \mathbb{Y}_1^{M_n^\delta} $ has the same distribution as $\mathbb{X}_1^{K_n}$.}
\ende

{
Let $\delta \in (0,1)$, $n \in \N$. Let $\varphi = \{(x_1,u_1),\dots,(x_K,u_K)\} \subset W_n \times [0,1]$ be any finite subset of size $K$, $K \in \N$, and denote $\varphi' := \{x_1,\dots,x_K\}$. Define a partition of $\varphi = T_\delta(\varphi) \cup T^c_\delta(\varphi)$, where
\begin{align*}
T_\delta(\varphi) &:= \{(x_i,u_i) \in \varphi: \varphi'(b^{(n)}(x_i,R)) \leq 1\}, \\
T^c_\delta(\varphi) &:= \{(x_i,u_i) \in \varphi: \varphi'(b^{(n)}(x_i,R)) \geq 2\}.
\end{align*}
Then, define
\begin{equation}\label{def:SetEndelta}
   E_{n,\delta} := \{\varphi \subset W_n \times [0,1]: T_\delta(\varphi) \subset W_n \times [\delta,1] \text{ and } T^c_\delta(\varphi) \subset W_n \times [0,\delta ) \}
\end{equation}}
which corresponds to the event of deleting all points having at least one $R$-neighbor.

We start with the concentration inequality for the number of neighbors. To simplify the notation, let us denote for any $r > 0$ and $\eps \in (0,1)$
\begin{equation*}
    a_{r,\eps}:= \frac{3}{2} r^d v_d \la \eps \left(1+\frac{\eps}{2}\right).
\end{equation*}

\bel[Concentration inequality for the number of neighbors]\label{lemma:AuxConcBound}
Let  $\gamma,r > 0$  and  $n \in \N$. Let $\eps > 0$ be small enough so that $\gamma > e^2  a_{r,\eps}$ and let $0 < \delta < \frac{\eps}{2}$. Denote
\begin{equation*}
    \one_{\eps,\delta}(\mathbb{X}_1^{K_n},K_n,\mathbb{U}_1^{K_n}):=\one\left[N_n(\mathbb{X}_1^{K_n}) \leq n\delta\right]
      \one\left[K_n \in (n\pm n\frac{\eps}{2})\right]
      \one\left[(\mathbb{X}_1^{K_n},\mathbb{U}_1^{K_n}) \in E_{n,\delta}\right].
\end{equation*}Then, for $n$ large enough it holds that 
\begin{align*}
\begin{split}
    \E \Biggl[
    \one_{\eps,\delta}(\mathbb{X}_1^{K_n},K_n,\mathbb{U}_1^{K_n}) &\cdot \one\Bigl[\mathbb{X}_{\{i \in \{1,\dots,K_n\}: U_i \geq \delta \}} \Bigl(\bigcup_{i \in \{1,\dots,K_n\}: U_i < \delta } b^{(n)}(X_i,r)\Bigr) > n\gamma \Bigr] \Biggr] \\
    &\leq \exp{\left(-\frac{n\gamma}{2}\cdot \log\left(\frac{\gamma}{a_{r,\eps}}\right)\right)}.
    \end{split}
\end{align*}
\enl

Then, we state the auxiliary unnormalized lower bound for the hard-core Hamiltonian $H_n$ in the Poisson case. Let $\oneJ{O} \subset \R^J$ be any open set, we denote 
\begin{equation*}
\oneJ{O}\oplus b(0,\Delta):= (O_1\oplus b(0,\Delta),\dots, O_J\oplus b(0,\Delta))^T.
\end{equation*}

\bel[] \label{lemma:HCAuxLBI}
Let $J \in \N$ and $\Delta,c,r > 0$.  Let $\oneJ{O} \subset \R^J$ be an open set. Assume that the score functions $\oneJ{\bfxi}$ are bounded by the constant $c$ and $r$-local. 

Denote $L := \inf \{I(P) + P^o(V): P \in \mc{P}_s, P^o({\bf1}) = \la, P^o(\oneJ{\bfxi}) \in \la \oneJ{O}\}$ and assume that $L < \infty$. {Choose $\eps_0 = \eps_0(\Delta, r, c) > 0$} small enough so that for all $\eps \in (0,\eps_0)$, it holds that
\begin{itemize}
    \item[(i)] $\eps <  \min\{ \frac{1}{2}, \frac{\Delta}{6c}  \}$ such that $ e^2 a_{r,\eps} < \frac{\Delta}{6c} $
    \item[(ii)] $-\la \frac{\Delta}{12c} \log \left(\frac{\Delta}{6ca_{r,\eps}} \right) < -L + b,$ where $b = -\frac{\la}{4} + \frac{\la}{4}\cdot \log\frac{1}{4} + \frac{3\la}{2}\log(\frac{3}{4})$.
\end{itemize}
Then, for all $\eps \in (0,\eps_0) $, it holds that
\begin{align*}
\begin{split} 
    &\liminf{ { \frac{1}{\abs{W_n}} }\log\E e^{-H_n(\mathbb{X}_1^{K_n})} \one \left[ K_n \in (n\pm n\eps)\right] \one \left[\onenJ{\bfXi}{n}(\mathbb{X}_1^{K_n})  \in \oneJ{O}\oplus b(0,\Delta) \right]} \geq -L.
    \end{split}
\end{align*}
\enl

At last, we present the key inequalities to get from the binomial case to the Poisson case. Recall the equality (\ref{factPois}) and the definition of the random variables $\mc{Z}_{n,+}, \mc{Z}_{n,-}$ (note that they depend on $\eps$) from (\ref{def:Z1Z2}).

\bel[] \label{lemma:HCAuxLBII}
Let $J \in \N$ and $\Delta, c,r > 0$. Let $0< \eps < 1$. Let $\oneJ{O},\oneJ{U} \subset \R^J$ be open sets such that $ \oneJ{O} \oplus b(0,2\Delta) \subset \oneJ{U}$. Assume that the score functions $\oneJ{\bfxi}$ are bounded by the constant $c$ and $r$-local. Denote $$I:= \E e^{-H_n(\mathbb{X}_1^{n})}\one\left[K_n = n \right]\one \left[\onenJ{\bfXi}{n}(\mathbb{X}_1^{n}) \in \oneJ{U}\right].$$ Then, for $n$ large enough,
\begin{align}
    I  &\geq \frac{1}{n\eps} \E e^{-H_n(\mathbb{X}_1^{K_n})}\one\left[K_n \in \left[n,n+n\eps\right) \right]\one \left[\onenJ{\bfXi}{n}(\mathbb{X}_1^{K_n}) \in \oneJ{O}\oplus b(0,\Delta)\right] \nonumber \\
    & \qquad \qquad \qquad \qquad -  \frac{1}{n\eps}\P( \mc{Z}_{n,+} > n \frac{\Delta}{6c} ),  \label{IneqI1} \\
  I &\geq \frac{\left( 1 - \la v_d R^d\right)^{n\eps } }{n\eps}  \E \Bigl[ e^{-H_n(\mathbb{X}_1^{K_n})}\one\left[K_n \in \left( n-n\eps,n\right)\right]\one \left[\onenJ{\bfXi}{n}(\mathbb{X}_1^{K_n}) \in \oneJ{O}\oplus b(0,\Delta)\right]  \Bigr] \nonumber \\
    & \qquad \qquad \qquad \qquad -  \frac{1}{n\eps}\P( \mc{Z}_{n,-} > n \frac{\Delta}{6c} ). \label{IneqI2}
\end{align}
\enl

Now, we are ready to prove the unnormalized lower bound. 

\bep[Proof of Lemma \ref{lemma:HCUnnormLB}] First we prove that for any $\Delta > 0$ and any $\oneJ{O} \subset \R^J$ open set such that $ \oneJ{O} \oplus b(0,2\Delta) \subset \oneJ{U}$ the following weaker lower bound holds
\begin{align}
\begin{split}\label{weakerbound}
    &\liminf{ { \frac{1}{\abs{W_n}} }\log\int \exp(-H_n(\om)) \one \left[ \abs{\om} = n\right] \one \left[\onenJ{\bfXi}{n}(\om)  \in \oneJ{U}\right]\d \Pi_n(\om)} \\
    & \qquad  \qquad  \geq - \inf \{I(P) + P^o(V): P \in \mc{P}_s, P^o({\bf1}) = \la, P^o(\oneJ{\bfxi}) \in \la \oneJ{O}\}.
    \end{split}
\end{align}
The proof of (\ref{HCUnnormLB}) then follows from (\ref{weakerbound}) in the same way as (\ref{lb}) followed from Lemma \ref{auxlowerbound} in the proof of Lemma \ref{cor:boundedLDP}.

Let $\Delta > 0$ and  $\oneJ{O} \subset \R^J$ be an open set such that $\oneJ{O} \oplus b(0,2\Delta) \subset \oneJ{U}$. If $L = \infty$, where we denote $L := \inf \{I(P) + P^o(V): P \in \mc{P}_s, P^o({\bf1}) = \la, P^o(\oneJ{\bfxi}) \in \la \oneJ{O}\}$, then the bound holds trivially. Assume therefore that $L < \infty$. Choose $\eps_1 $ small enough so that $\eps_1 < \eps_0$  from Lemma \ref{lemma:HCAuxLBI} (for $c,r,\Delta$ from our setting) and so that for all $\eps \in (0,\eps_1)$ we have that  $-\la \frac{\Delta}{24c} \log \left(\frac{\Delta}{6 c p_\eps} \right) < -L + \la \log \left(1-\la v_d R^d\right)$, where $p_\eps:= \eps \la  r^d  v_d$.

Take any $\eps \in (0,\eps_1)$ and denote for simplicity 
\begin{align*}
    &I_-:= \E \Bigl[ e^{-H_n(\mathbb{X}_1^{K_n})}\one\left[K_n \in \left( n-n\eps,n\right)\right]\one \left[\onenJ{\bfXi}{n}(\mathbb{X}_1^{K_n}) \in \oneJ{O}\oplus b(0,\Delta)\right]  \Bigr],\\
    &I_+:= \E \Bigl[ e^{-H_n(\mathbb{X}_1^{K_n})}\one\left[K_n \in \left[n,n+n\eps\right) \right]\one \left[\onenJ{\bfXi}{n}(\mathbb{X}_1^{K_n}) \in \oneJ{O}\oplus b(0,\Delta)\right]\Bigr].
\end{align*}

Using {the} bounds (\ref{IneqI1}) and (\ref{IneqI2}) from Lemma \ref{lemma:HCAuxLBII} and the concentration inequalities (\ref{IneqForZ12}), we can bound with {$\eta := 1 - \lambda v_dR^d$},
\begin{align*}
     \begin{split}
         &\liminf { \frac{1}{\abs{W_n}} } \log \Bigl( \E e^{-H_n(\mathbb{X}_1^{n})}\one\left[K_n = n \right]\one \left[\onenJ{\bfXi}{n}(\mathbb{X}_1^{n}) \in \oneJ{U}  \right] \Bigr) \\
         &\geq \liminf { \frac{1}{\abs{W_n}} } \log \Biggl( \frac{\eta^{n\eps } }{2n\eps}  I_- -  \frac{1}{2n\eps}\P( \mc{Z}_{n,-} > n \frac{\Delta}{6c} ) + \frac{1}{2n\eps} I_+ -  \frac{1}{2n\eps}\P( \mc{Z}_{n,+} > n \frac{\Delta}{6c} ) \Biggr)\\
         &\geq 0+ \liminf { \frac{1}{\abs{W_n}} } \log \Bigl( \eta^{n\eps } (I_- + I_+ ) -  \P( \mc{Z}_{n,-} > n \frac{\Delta}{6c} ) -  \P( \mc{Z}_{n,+} > n \frac{\Delta}{6c} ) \Bigr)\\
         &\overset{(\ref{IneqForZ12})}{\geq}
         \liminf { \frac{1}{\abs{W_n}} } \log \Bigl( \eta^{n\eps } (I_- + I_+ ) - 2 \exp\Bigl(-\frac{n\Delta}{24c}\log \frac{\Delta}{6c p_\eps}\Bigr)  \Bigr) \\
         &\geq \eps \la \log \eta +
         \liminf { \frac{1}{\abs{W_n}} } \log \Biggl(  (I_- + I_+ ) -  2\frac{\exp\left(-\frac{n\Delta}{24c}\log \frac{\Delta}{6c p_\eps}\right)}{\eta^{n\eps }} \Biggr).
     \end{split}
\end{align*}
Recall the choice of $\eps_1$ from the beginning of this proof. We can easily see that 
\begin{align*}
    \begin{split}
        \limsup { \frac{1}{\abs{W_n}} } \log \Biggl( 2\frac{\exp\left(-\frac{n\Delta}{24c}\log \frac{\Delta}{6c p_\eps}\right)}{\eta^{n\eps }} \Biggr)&= -\la \frac{\Delta}{24c}\log \frac{\Delta}{6c p_\eps} - \eps \la \log \left(\eta\right) \\ &\leq - \la \frac{\Delta}{24c}\log \frac{\Delta}{6c p_\eps} -  \la \log \left( \eta\right) < -L.
    \end{split}
\end{align*}
On the other hand, Lemma \ref{lemma:HCAuxLBI} gives us that
\begin{align*}
    \begin{split}
        \liminf & { \frac{1}{\abs{W_n}} } \log  (I_- + I_+ )\\ &=\liminf{ { \frac{1}{\abs{W_n}} }\log\E e^{-H_n(\mathbb{X}_1^{K_n})} \one \left[ K_n \in (n\pm n\eps)\right] \one \left[\onenJ{\bfXi}{n}(\mathbb{X}_1^{K_n})  \in \oneJ{O}\oplus b(0,\Delta)  \right]},
    \end{split}
\end{align*}
{and the right-hand side is bounded below by $-L$}.
Therefore, 
\begin{align*}
    \liminf { \frac{1}{\abs{W_n}} } \log \Biggl(  (I_- + I_+ ) & -  2\frac{\exp\left(-\frac{n\Delta}{24c}\log \frac{\Delta}{6c p_\eps}\right)}{\left( 1 - \la v_d R^d\right)^{n\eps }} \Biggr) \\
    &= \liminf { \frac{1}{\abs{W_n}} } \log  (I_- + I_+ )  \geq - L.
\end{align*}
Altogether for any $\eps \in (0,\eps_1)$ 
\begin{align*}
         \liminf { \frac{1}{\abs{W_n}} } \log  \Bigl( \E e^{-H_n(\mathbb{X}_1^{n})}  \one\left[K_n = n , \onenJ{\bfXi}{n}(\mathbb{X}_1^{n}) \in \oneJ{U}  \right] \Bigr)
         &\geq \eps \la \log \left( 1 - \la v_d R^d\right) -L.
\end{align*}
Taking $\eps \downarrow 0$ finishes the proof. 
\enp

\subsection{Proofs of auxiliary lemmas}

\bep[Proof of Lemma \ref{lemma:AuxConcBound}]We will prove the result by conditioning on the sequence $\mathbb{U}$ and using the joint independence from the construction of our coupling.

Fix a sequence $\{u_i\}_{i \in \N} \in \left[0,1\right]^{\N}$ and denote 
$$
n({\bf u_1^{K_n}}) := \abs{\{i \in {1,\dots,K_n}: u_i \geq \delta\}}.
$$
We will w.l.o.g.\,assume that $\{i \in {1,\dots,K_n}: u_i \geq \delta\} = \{1,\dots,n({\bf u_1^{K_n}})\}$. Under the event corresponding to the indicator $\one_{\eps,\delta}(\mathbb{X}_1^{K_n},K_n,{\bf u_1^{K_n}})$  and thanks to the choice of $\delta$, it holds that 
$
    n + n\frac{\eps}{2} \geq  n({\bf u_1^{K_n}}) \geq K_n - n\delta \geq n - n\frac{\eps}{2} - n\delta \geq n - n\eps.
$

Denote $\one_{(a,b)}(n({\bf u_1^{K_n}})) := \one\left[a \leq n({\bf u_1^{K_n}}) \leq b \right]$. We can bound 
\begin{align*}
    \begin{split}
          \E \Biggl[
    \one_{\eps,\delta}&(\mathbb{X}_1^{K_n},K_n,{\bf u_1^{K_n}}) \cdot \one\Bigl[\mathbb{X}_{\{i \in \{1,\dots,K_n\}: U_i \geq \delta \}}\Bigl(\bigcup_{i \in \{1,\dots,K_n\}: u_i < \delta } b^{(n)}(X_i,r)\Bigr) > n\gamma \Bigr] \Biggr]  \\
    &\leq  \E  \Biggl[
    \one_{\eps,\delta}(\mathbb{X}_1^{K_n},K_n,{\bf u_1^{K_n}}) \cdot \one\Bigl[\mathbb{X}_1^{n({\bf u_1^{K_n}})} \Bigl(\bigcup_{n({\bf u_1^{K_n}}) + 1}^{n+n\frac{\eps}{2}} b^{(n)}(X_i,r)\Bigr) > n\gamma \Bigr] \Biggr] \\
    &=  \E  \Biggl[
    \one_{(n - n\eps,n + n\frac{\eps}{2})}(n({\bf u_1^{K_n}})) \cdot \one_{\eps,\delta}(\mathbb{X}_1^{K_n},K_n,{\bf u_1^{K_n}}) \\ & \hspace{5cm} \cdot \one\Bigl[\mathbb{X}_1^{n({\bf u_1^{K_n}})} \Bigl(\bigcup_{n({\bf u_1^{K_n}}) + 1}^{n+n\frac{\eps}{2}} b^{(n)}(X_i,r)\Bigr) > n\gamma \Bigr] \Biggr]\\
    &\leq  \E  \Biggl[
    \one_{(n - n\eps,n + n\frac{\eps}{2})}(n({\bf u_1^{K_n}})) \cdot \one\Bigl[\mathbb{X}_1^{n({\bf u_1^{K_n}})} \Bigl(\bigcup_{n({\bf u_1^{K_n}}) + 1}^{n+n\frac{\eps}{2}} b^{(n)}(X_i,r)\Bigr) > n\gamma \Bigr] \Biggr].
    \end{split}
\end{align*}

Altogether we get, using the concentration inequality for binomial distribution (Lemma 1.1 in \cite{PenroseBook}), that
\begin{align*}
    \begin{split}
        &\E  \Biggl[
    \one_{(n - n\eps,n + n\frac{\eps}{2})}(n({\bf u_1^{K_n}})) \cdot \one\Bigl[\mathbb{X}_1^{n({\bf u_1^{K_n}})} \Bigl(\bigcup_{n({\bf u_1^{K_n}}) + 1}^{n+n\frac{\eps}{2}} b^{(n)}(X_i,r)\Bigr) > n\gamma \Bigr] \Biggr] \\
    &=   \E\Biggl[ \E  \Biggl[
    \one_{(n - n\eps,n + n\frac{\eps}{2})}(n({\bf u_1^{K_n}})) \\
    & \hspace{3cm} \cdot \one\Bigl[\mathbb{X}_1^{n({\bf u_1^{K_n}})} \Bigl(\bigcup_{n({\bf u_1^{K_n}}) + 1}^{n+n\frac{\eps}{2}} b^{(n)}(X_i,r)\Bigr) > n\gamma \Bigr]  \lvert \mathbb{X}_{n({\bf u_1^{K_n}}) + 1}^{n+n\frac{\eps}{2}} ,K_n  \Biggr]\Biggr]\\
     &=   \E\Biggl[ \one_{(n - n\eps,n + n\frac{\eps}{2})}(n({\bf u_1^{K_n}})) \\
     & \hspace{3cm} \cdot \E  \Biggl[
     \one\Bigl[\mathbb{X}_1^{n({\bf u_1^{K_n}})} \Bigl(\bigcup_{n({\bf u_1^{K_n}}) + 1}^{n+n\frac{\eps}{2}} b^{(n)}(X_i,r)\Bigr) > n\gamma \Bigr]  \lvert \mathbb{X}_{n({\bf u_1^{K_n}}) + 1}^{n+n\frac{\eps}{2}} ,K_n  \Biggr]\Biggr]\\
     &\leq  \E\Biggl[ \one_{(n - n\eps,n + n\frac{\eps}{2})}(n({\bf u_1^{K_n}})) \cdot \exp{\Bigl(-\frac{n \gamma}{2}\cdot \log \Bigl( \frac{n \gamma}{n({\bf u_1^{K_n}}) \cdot p\Bigl(\mathbb{X}_{n({\bf u_1^{K_n}}) + 1}^{n+n\frac{\eps}{2}}\Bigr)}\Bigr) \Bigr)}\Biggr],
    \end{split}
\end{align*}

where we have used the fact that 
\begin{equation*}
    \mathbb{X}_1^{n({\bf u_1^{K_n}})} \Bigl(\bigcup_{n({\bf u_1^{K_n}}) + 1}^{n+n\frac{\eps}{2}} b^{(n)}(X_i,r)\Bigr) \,\lvert \,\, \mathbb{X}_{n({\bf u_1^{K_n}}) + 1}^{n+n\frac{\eps}{2}} ,K_n \, \sim \, Bi\Bigl(n({\bf u_1^{K_n}}), p\Bigl(\mathbb{X}_{n({\bf u_1^{K_n}}) + 1}^{n+n\frac{\eps}{2}} \Bigr)\Bigr),
\end{equation*}
where $p\Bigl(\mathbb{X}_{n({\bf u_1^{K_n}}) + 1}^{n+n\frac{\eps}{2}} \Bigr):= \frac{\abs{\bigcup_{n({\bf u_1^{K_n}}) + 1}^{n+n\frac{\eps}{2}} b^{(n)}(X_i,r)}}{\abs{W_n}}$.

Since 
$n + n\frac{\eps}{2} \geq n({\bf u_1^{K_n}})  \geq n - n\eps \implies n({\bf u_1^{K_n}}) \cdot p\Bigl(\mathbb{X}_{n({\bf u_1^{K_n}}) + 1}^{n+n\frac{\eps}{2}} \Bigr) \leq n a_{r,\eps},
$, the proof is finished.
\enp

Now, we prove Lemma \ref{lemma:HCAuxLBI}, where we need the following inequality. It shows that the event $\{(\mathbb{X}_1^{K_n},\mathbb{U}_1^{K_n}) \in E_{n,\delta}\}$ where all the points with at least one $R$-neighbor are deleted (see (\ref{def:SetEndelta})), has a sufficiently large probability. 
Let $0 < \eps,\delta < 1 $, then a.s.
\begin{align}
        &\one\left[N_n(\mathbb{X}_1^{K_n}) \leq n\delta\right]
     \one\left[K_n \in (n\pm n\frac{\eps}{2})\right] \cdot
     \E \Bigl[ \one\left[(\mathbb{X}_1^{K_n},\mathbb{U}_1^{K_n}) \in E_{n,\delta}\right] \lvert (\mathbb{X}_1^{K_n},\mathbb{Y}_1^{M_n^\delta}) \Bigr] \nonumber \\
     &=\one\left[N_n(\mathbb{X}_1^{K_n}) \leq n\delta\right]
     \one\left[K_n \in (n\pm n\frac{\eps}{2})\right] \cdot
     \E \Bigl[ \one\left[(\mathbb{X}_1^{K_n},\mathbb{U}_1^{K_n}) \in E_{n,\delta}\right] \lvert \mathbb{X}_1^{K_n} \Bigr] \nonumber \\
     & = \one\left[N_n(\mathbb{X}_1^{K_n}) \leq n\delta\right]
     \one\left[K_n \in (n\pm n\frac{\eps}{2})\right] \cdot\delta^{N_n(\mathbb{X}_1^{K_n})} \cdot (1-\delta)^{K_n - N_n(\mathbb{X}_1^{K_n})} \nonumber \\
     &\geq \one\left[N_n(\mathbb{X}_1^{K_n}) \leq n\delta\right]
     \one\left[K_n \in (n\pm n\frac{\eps}{2})\right] \cdot\delta^{n\delta} \cdot (1-\delta)^{n+n\eps}. \label{IneqForEndelta}
\end{align}

{
We also need the following auxiliary bound.

\bel
Assume that we are in the situation of Lemma \ref{lemma:HCAuxLBI}. Then, it holds that
\begin{align}
    \begin{split}\label{vz2}
         &\liminf \frac{1}{\abs{W_n}} \log \E \Biggl[ e^{-H^s_n(\mathbb{X}_1^{K_n})}  \one \left[\onenJ{\bfXi}{n}(\mathbb{X}_1^{K_n})  \in \oneJ{O} \right] \\
         & \hspace{4cm} \cdot \one\left[N_n(\mathbb{X}_1^{K_n}) \leq n\delta\right]
     \one\left[K_n \in (n\pm n\frac{\eps}{2})\right]
    \Biggr] {\geq} -L. 
    \end{split}
\end{align}
\enl

\bep
First, note that thanks to Corollary 3.2 from \cite{GeorgiiLDP} we have that for any $s > 0$ 
\begin{align} \label{vz1}
    \begin{split}
        \liminf&{ \frac{1}{\abs{W_n}}\log\E e^{-H^s_n(\mathbb{X}_1^{K_n})} \one \left[ K_n \in (n\pm n\frac{\eps}{2})\right] \one \left[\onenJ{\bfXi}{n}(\mathbb{X}_1^{K_n})  \in \oneJ{O} \right]} \\
        & \geq - \inf \{I(P) + P^o(V^s): P \in \mc{P}_s, P^o({\bf1}) \in (\la \pm \la \frac{\eps}{2}), P^o(\oneJ{\bfxi}) \in \la\oneJ{O}\} \\
        & \geq - \inf \{I(P) + P^o(V^s): P \in \mc{P}_s, P^o({\bf1}) = \la , P^o(\oneJ{\bfxi}) \in \la \oneJ{O}\} \\
        & \geq - \inf \{I(P) + P^o(V): P \in \mc{P}_s, P^o({\bf1}) = \la , P^o(\oneJ{\bfxi}) \in \la \oneJ{O}\} 
        = -L.
    \end{split}
\end{align}

For any $\delta > 0$ we can, thanks to (\ref{vz1}) and Lemma \ref{lemma:chooseS}, choose $s_\delta$ such that for any $s \geq s_\delta$, it holds that 
\begin{align*}
    \begin{split}
         &\liminf \frac{1}{\abs{W_n}} \log \E \Biggl[ e^{-H^s_n(\mathbb{X}_1^{K_n})}  \one \left[\onenJ{\bfXi}{n}(\mathbb{X}_1^{K_n})  \in \oneJ{O}, N_n(\mathbb{X}_1^{K_n}) \leq n\delta, K_n \in (n\pm n\frac{\eps}{2})\right]
    \Biggr] \\
       &= \liminf \frac{1}{\abs{W_n}}\log \Biggl[ \E e^{-H^s_n(\mathbb{X}_1^{K_n})} \one \left[ K_n \in (n\pm n\frac{\eps}{2})\right] \one \left[\onenJ{\bfXi}{n}(\mathbb{X}_1^{K_n})  \in \oneJ{O} \right]  \\
       &\hspace{0.5cm} - \E e^{-H^s_n(\mathbb{X}_1^{K_n})} \one \left[ K_n \in (n\pm n\frac{\eps}{2})\right] \one \left[\onenJ{\bfXi}{n}(\mathbb{X}_1^{K_n})  \in \oneJ{O} \right] \one\left[N_n(\mathbb{X}_1^{K_n}) > n\delta\right] \Biggr] \\
       & =  \liminf \frac{1}{\abs{W_n}}\log  \E e^{-H^s_n(\mathbb{X}_1^{K_n})} \one \left[ K_n \in (n\pm n\frac{\eps}{2})\right] \one \left[\onenJ{\bfXi}{n}(\mathbb{X}_1^{K_n})  \in \oneJ{O} \right] \overset{(\ref{vz1})}{\geq} -L. 
    \end{split}
\end{align*}
\enp
}

{ Now, we can prove Lemma \ref{lemma:HCAuxLBI}}

\bep[Proof of Lemma \ref{lemma:HCAuxLBI}]
Fix $\eps \in (0,\eps_0)$ and recall the notation $\one_{\eps,\delta}(\mathbb{X}_1^{K_n},K_n,\mathbb{U}_1^{K_n})$ from Lemma \ref{lemma:AuxConcBound} and the coupling from Definition \ref{def:Coupling}. To simplify the notation, denote
 \begin{equation*}
     Z_{r,\delta}(\mathbb{X},\mathbb{U)}:= \mathbb{X}_{\{i \in \{1,\dots,K_n\}: U_i \geq \delta \}} \left(\bigcup_{i \in \{1,\dots,K_n\}: U_i < \delta } b^{(n)}(X_i,r)\right).
 \end{equation*}
Let $\delta < \frac{\eps}{2}$, $s \geq s^\delta$ and $\gamma = \frac{\Delta}{6c}$.

Using the bound from Lemma \ref{lemma:AuxConcBound} and our choice of $\eps$ and $\gamma$, we get that 
\begin{align*}
       \limsup & { \frac{1}{\abs{W_n}} } \log \E \Bigl[
    \one_{\eps,\delta}(\mathbb{X}_1^{K_n},K_n,\mathbb{U}_1^{K_n}) \one\left[Z_{r,\delta}(\mathbb{X},\mathbb{U)} > n\gamma \right] \Bigr] \\
    & \qquad \qquad  \leq \limsup { \frac{1}{\abs{W_n}} } \cdot \left(-\frac{n\gamma}{2}\cdot \log\left(\frac{\gamma}{a_{r,\eps}}\right)\right) = -\la \frac{\Delta}{12c} \cdot \log \left(\frac{\Delta}{6ca_{r,\eps}} \right)  < -L + b .
\end{align*}

Then, we can bound
\begin{align*}
        \E  & e^{-H_n(\mathbb{X}_1^{K_n})} \one \left[ K_n \in (n\pm n\eps)\right] \one \left[\onenJ{\bfXi}{n}(\mathbb{X}_1^{K_n})  \in \oneJ{O}\oplus b(0,\Delta) \right]\\
        = &\E  e^{-H_n(\mathbb{Z}_1^{L_n})} \one \left[ L_n \in (n\pm n\eps)\right] \one \left[\onenJ{\bfXi}{n}(\mathbb{Z}_1^{L_n})  \in \oneJ{O}\oplus b(0,\Delta) \right]\\
        \geq & \E \Bigl[ e^{-H_n(\mathbb{Z}_1^{L_n})} \one \left[ L_n \in (n\pm n\eps), \onenJ{\bfXi}{n}(\mathbb{Z}_1^{L_n})  \in \oneJ{O}\oplus b(0,\Delta),   \abs{\mathbb{Y}_1^{M_n^\delta}}=0\right] \one_{\eps,\delta}(\mathbb{X}_1^{K_n},K_n,\mathbb{U}_1^{K_n})  \Bigr]\\
        \overset{(a)}{\geq} & \E \Bigl[ e^{-H^s_n(\mathbb{X}_1^{K_n})} \one \left[ L_n \in (n\pm n\eps)\right] \one \left[\onenJ{\bfXi}{n}(\mathbb{Z}_1^{L_n})  \in \oneJ{O}\oplus b(0,\Delta) \right]  \one_{\eps,\delta}(\mathbb{X}_1^{K_n},K_n,\mathbb{U}_1^{K_n}) \one \left[\abs{\mathbb{Y}_1^{M_n^\delta}}=0\right]\Bigr] \\
                \overset{(b)}{\geq} &   \E \Bigl[ e^{-H^s_n(\mathbb{X}_1^{K_n})}  \one \left[\onenJ{\bfXi}{n}(\mathbb{Z}_1^{L_n})  \in \oneJ{O}\oplus b(0,\Delta) \right] \one_{\eps,\delta}(\mathbb{X}_1^{K_n},K_n,\mathbb{U}_1^{K_n})  \one \left[\abs{\mathbb{Y}_1^{M_n^\delta}}=0\right] \one\left[Z_{r,\delta}(\mathbb{X},\mathbb{U)} \leq n\gamma \right] \Bigr] 
        \end{align*}
        Therefore,
\begin{align*}                \E  e^{-H_n(\mathbb{X}_1^{K_n})} & \one \left[ K_n \in (n\pm n\eps)\right] \one \left[\onenJ{\bfXi}{n}(\mathbb{X}_1^{K_n})  \in \oneJ{O}\oplus b(0,\Delta) \right]\\
         \overset{(c)}{\geq} &   \E \Bigl[ e^{-H^s_n(\mathbb{X}_1^{K_n})}  \one \left[\onenJ{\bfXi}{n}(\mathbb{X}_1^{K_n})  \in \oneJ{O} \right] \one_{\eps,\delta}(\mathbb{X}_1^{K_n},K_n,\mathbb{U}_1^{K_n}) \\
        & \hspace{4.5cm} \cdot \one \left[\abs{\mathbb{Y}_1^{M_n^\delta}}=0\right]   \one\left[Z_{r,\delta}(\mathbb{X},\mathbb{U)} \leq n\gamma \right] \Bigr]\\
         \geq &   \E \Bigl[ e^{-H^s_n(\mathbb{X}_1^{K_n})}  \one \left[\onenJ{\bfXi}{n}(\mathbb{X}_1^{K_n})  \in \oneJ{O} \right] \one_{\eps,\delta}(\mathbb{X}_1^{K_n},K_n,\mathbb{U}_1^{K_n}) \one \left[\abs{\mathbb{Y}_1^{M_n^\delta}}=0\right] \Bigr]  \\
        & \hspace{4cm} -   \E \Bigl[   \one_{\eps,\delta}(\mathbb{X}_1^{K_n},K_n,\mathbb{U}_1^{K_n})  \cdot  \one\left[Z_{r,\delta}(\mathbb{X},\mathbb{U)}>n\gamma \right] \Bigr], \\
\end{align*}
where we have used that under the corresponding indicators, it holds that 
\begin{itemize}
    \item[(a)] $H_n(\mathbb{Z}_1^{L_n}) = 0 \leq H^s_n(\mathbb{X}_1^{K_n})$,
    \item[(b)] 
    \begin{align*}
    \one \left[ L_n \in (n\pm n\eps)\right] \one_{\eps,\delta}(\mathbb{X}_1^{K_n},K_n,\mathbb{U}_1^{K_n}) &\one \left[\abs{\mathbb{Y}_1^{M_n^\delta}}=0\right] \\ &= \one_{\eps,\delta}(\mathbb{X}_1^{K_n},K_n,\mathbb{U}_1^{K_n}) \one \left[\abs{\mathbb{Y}_1^{M_n^\delta}}=0\right],
    \end{align*}
    \item[(c)] $\abs{\Xi_{j,n}(\mathbb{Z}_1^{L_n}) - \Xi_{j,n}(\mathbb{X}_1^{K_n})} \leq \frac{2c}{n}Z_{r,\delta}(\mathbb{X},\mathbb{U)} + c\delta < \Delta$.
\end{itemize}
Furthermore, we can also bound 
\begin{align*}
    \begin{split}
        &\E \Bigl[ e^{-H^s_n(\mathbb{X}_1^{K_n})}  \one \left[\onenJ{\bfXi}{n}(\mathbb{X}_1^{K_n})  \in \oneJ{O} \right] \one_{\eps,\delta}(\mathbb{X}_1^{K_n},K_n,\mathbb{U}_1^{K_n}) \one \left[\abs{\mathbb{Y}_1^{M_n^\delta}}=0\right] \Bigr] \\ 
        & = \E \Bigl[ e^{-H^s_n(\mathbb{X}_1^{K_n})}  \one \left[\onenJ{\bfXi}{n}(\mathbb{X}_1^{K_n})  \in \oneJ{O} \right]\one \left[\abs{\mathbb{Y}_1^{M_n^\delta}}=0\right] 
        \\  &\hspace{3.5cm} \cdot \E \Bigl[  \one_{\eps,\delta}(\mathbb{X}_1^{K_n},K_n,\mathbb{U}_1^{K_n})  \lvert (\mathbb{X}_1^{K_n},\mathbb{Y}_1^{M_n^\delta}) \Bigr] \Bigr]\\
         & = \E \Bigl[ e^{-H^s_n(\mathbb{X}_1^{K_n})}  \one \left[\onenJ{\bfXi}{n}(\mathbb{X}_1^{K_n})  \in \oneJ{O} \right]\one \left[\abs{\mathbb{Y}_1^{M_n^\delta}}=0\right]  \one\left[N_n(\mathbb{X}_1^{K_n}) \leq n\delta\right]
     \\
        &  \hspace{3.5cm} \cdot \one\left[K_n \in (n\pm n\frac{\eps}{2})\right] 
     \E \Bigl[ \one\left[(\mathbb{X}_1^{K_n},\mathbb{U}_1^{K_n}) \in E_{n,\delta}\right] \lvert (\mathbb{X}_1^{K_n},\mathbb{Y}_1^{M_n^\delta}) \Bigr] \Bigr] \\ 
     &\overset{(\ref{IneqForEndelta})}{\geq} \E \Bigl[ e^{-H^s_n(\mathbb{X}_1^{K_n})}  \one \left[\onenJ{\bfXi}{n}(\mathbb{X}_1^{K_n})  \in \oneJ{O} \right]\one \left[\abs{\mathbb{Y}_1^{M_n^\delta}}=0\right] \\
        & \hspace{4cm} \cdot \one\left[N_n(\mathbb{X}_1^{K_n}) \leq n\delta\right]
     \one\left[K_n \in (n\pm n\frac{\eps}{2})\right]
     \delta^{n\delta}  (1-\delta)^{n+n\eps} \Bigr] \\
      &= \delta^{n\delta} \cdot (1-\delta)^{n+n\eps} e^{-n\delta}\E \Bigl[ e^{-H^s_n(\mathbb{X}_1^{K_n})}  \one \left[\onenJ{\bfXi}{n}(\mathbb{X}_1^{K_n})  \in \oneJ{O} \right] \\ & \hspace{4.5cm} \cdot \one\left[N_n(\mathbb{X}_1^{K_n}) \leq n\delta\right]
     \one\left[K_n \in (n\pm n\frac{\eps}{2})\right]
     \Bigr],
    \end{split} 
\end{align*}
which (together with $\eps < \frac{1}{2}$ and $\delta < \frac{\eps}{2} < \frac{1}{4}$) leads to the fact that
\begin{align*}
    \begin{split}
         \liminf &{ \frac{1}{\abs{W_n}} } \log \E \Bigl[ e^{-H^s_n(\mathbb{X}_1^{K_n})}  \one \left[\onenJ{\bfXi}{n}(\mathbb{X}_1^{K_n})  \in \oneJ{O} \right] \one_{\eps,\delta}(\mathbb{X}_1^{K_n},K_n,\mathbb{U}_1^{K_n}) \one \left[\abs{\mathbb{Y}_1^{M_n^\delta}}=0\right] \Bigr] \\
         &\geq \la(\delta \log \delta +  (1+\eps) \log (1-\delta)- \delta)\\ 
         & \hspace{2cm} + \liminf { \frac{1}{\abs{W_n}} } \log \E \Bigl[ e^{-H^s_n(\mathbb{X}_1^{K_n})}  \one \left[\onenJ{\bfXi}{n}(\mathbb{X}_1^{K_n})  \in \oneJ{O} \right] \\ & \hspace{5cm} \cdot \one\left[N_n(\mathbb{X}_1^{K_n}) \leq n\delta\right]
     \one\left[K_n \in (n\pm n\frac{\eps}{2})\right]
    \Bigr] 
    \\ &\overset{(\ref{vz2})}{\geq} - L + \la(\delta \log \delta +  (1+\eps) \log (1-\delta)- \delta) \geq -L + b.
    \end{split}
\end{align*}

Altogether we have proved that 
\begin{align*}
    \begin{split}
        &\liminf { \frac{1}{\abs{W_n}} } \log \E e^{-H_n(\mathbb{X}_1^{K_n})} \one \left[ K_n \in (n\pm n\eps)\right] \one \left[\onenJ{\bfXi}{n}(\mathbb{X}_1^{K_n})  \in \oneJ{O}\oplus b(0,\Delta) \right] \\
        &\geq \liminf { \frac{1}{\abs{W_n}} } \log \Biggl[ \E  e^{-H^s_n(\mathbb{X}_1^{K_n})}  \one \left[\onenJ{\bfXi}{n}(\mathbb{X}_1^{K_n})  \in \oneJ{O} \right] \one_{\eps,\delta}(\mathbb{X}_1^{K_n},K_n,\mathbb{U}_1^{K_n}) \one \left[\abs{\mathbb{Y}_1^{M_n^\delta}}=0\right]   \\
       & \hspace{3cm}   -  \E \Bigl[  \one_{\eps,\delta}(\mathbb{X}_1^{K_n},K_n,\mathbb{U}_1^{K_n})  \cdot  \one\left[Z_{r,\delta}(\mathbb{X},\mathbb{U)} > n\gamma \right]  \Bigr] \Biggr] \\
       &=  \liminf { \frac{1}{\abs{W_n}} } \log \E \Bigl[ e^{-H^s_n(\mathbb{X}_1^{K_n})}  \one \left[\onenJ{\bfXi}{n}(\mathbb{X}_1^{K_n})  \in \oneJ{O} \right] \one_{\eps,\delta}(\mathbb{X}_1^{K_n},K_n,\mathbb{U}_1^{K_n}) \one \left[\abs{\mathbb{Y}_1^{M_n^\delta}}=0\right] \Bigr]  \\
       &\geq - L + \la(\delta \log \delta +  (1+\eps) \log (1-\delta)- \delta).
    \end{split}
\end{align*}
For our fixed $\eps$ this inequality holds for any $\delta < \frac{\eps}{2}$.Taking $\delta \downarrow 0$ finishes the proof. 
\enp

{ Finally, we prove Lemma \ref{lemma:HCAuxLBII}.}

\bep[Proof of Lemma \ref{lemma:HCAuxLBII}]
We start with the proof of (\ref{IneqI1}).

Let $n$ be large enough, then 
\begin{itemize}
    \item[(a)] if $K_n \geq n$, then $H_n(\mathbb{X}_1^{n}) \leq H_n(\mathbb{X}_1^{K_n})$,
    \item[(b)] if $\one\left[K_n \in \left[n,n+n\eps\right) \right]\one \left[ \mc{Z}_{n,+} \leq n \frac{\Delta}{6c} \right] = 1$ then $$\abs{\Xi_{j,n}(\mathbb{X}_1^{n})-\Xi_{j,n}(\mathbb{X}_1^{K_n}))} \leq c \eps + \frac{2c}{n} \mc{Z}_{n,+} < \Delta.$$
\end{itemize}

Therefore, the following bounds hold
\begin{align*}
    &\E e^{-H_n(\mathbb{X}_1^{n})}\one\left[K_n = n \right]\one \left[\onenJ{\bfXi}{n}(\mathbb{X}_1^{n}) \in \oneJ{U}\right] \\ &\overset{(\ref{factPois})}{\geq} \frac{1}{n\eps} \E e^{-H_n(\mathbb{X}_1^{n})}\one\left[K_n \in \left[n,n+n\eps\right) \right]\one \left[\onenJ{\bfXi}{n}(\mathbb{X}_1^{n}) \in \oneJ{U} \right] \\
    &\overset{(a)}{\geq} \frac{1}{n\eps} \E e^{-H_n(\mathbb{X}_1^{K_n})}\one\left[K_n \in \left[n,n+n\eps\right) \right]\one \left[\onenJ{\bfXi}{n}(\mathbb{X}_1^{n}) \in \oneJ{U} \right] \one \left[ \mc{Z}_{n,+} \leq n \frac{\Delta}{6c} \right]\\
    &\overset{(b)}{\geq}  \frac{1}{n\eps} \E e^{-H_n(\mathbb{X}_1^{K_n})}\one\left[K_n \in \left[n,n+n\eps\right) \right]\one \left[\onenJ{\bfXi}{n}(\mathbb{X}_1^{K_n}) \in \oneJ{O}\oplus b(0,\Delta) \right] \one \left[ \mc{Z}_{n,+} \leq n \frac{\Delta}{6c} \right]\\
    &\geq \frac{1}{n\eps} \E e^{-H_n(\mathbb{X}_1^{K_n})}\one\left[K_n \in \left[n,n+n\eps\right) \right]\one \left[\onenJ{\bfXi}{n}(\mathbb{X}_1^{K_n}) \in \oneJ{O}\oplus b(0,\Delta)\right] \\ & \hspace{9cm} -  \frac{1}{n\eps}\P( \mc{Z}_{n,+} > n \frac{\Delta}{6c} ).
\end{align*}

To prove (\ref{IneqI2}) is a bit trickier, because (a) no longer holds. Using the fact that 
\begin{itemize}
    \item[(c)]if  $\one\left[K_n \in \left(n-n\eps,n\right) \right]\one \left[ \mc{Z}_{n,-} \leq n \frac{\Delta}{6c} \right] = 1$ then 
$$\abs{\Xi_{j,n}(\mathbb{X}_1^{n})-\Xi_{j,n}(\mathbb{X}_1^{K_n}))} \leq 3c \eps + \frac{2c \mc{Z}_{n,-}}{n} < \Delta,
$$
\end{itemize}
 we can bound 
\begin{align*}
    \begin{split}
         &\E e^{-H_n(\mathbb{X}_1^{n})}\one\left[K_n = n \right]\one \left[\onenJ{\bfXi}{n}(\mathbb{X}_1^{n}) \in \oneJ{U} \right] \\&\overset{(\ref{factPois})}{\geq} \frac{1}{n\eps} \E e^{-H_n(\mathbb{X}_1^{n})}\one\left[K_n \in \left(n-n\eps,n\right) \right]\one \left[\onenJ{\bfXi}{n}(\mathbb{X}_1^{n}) \in \oneJ{U} \right] \\
    &\overset{(c)}{\geq} \frac{1}{n\eps} \E e^{-H_n(\mathbb{X}_1^{n})}\one\left[K_n \in \left(n-n\eps,n\right) \right]\one \left[\onenJ{\bfXi}{n}(\mathbb{X}_1^{K_n}) \in \oneJ{O}\oplus b(0,\Delta) \right] \\ & \hspace{9cm} -  \frac{1}{n\eps}\P( \mc{Z}_{n,-} > n \frac{\Delta}{6c} ).
    \end{split}
\end{align*}
 Now, we find a lower bound for the first term. Denote for $J \in \N$ and ${\bf y}_1^J \in (\R^d)^J$
\begin{equation*}
    E_n({\bf y}_1^J):= \{\om \in \Omega_{W_n}: \forall x \in \om \quad \abs{\om^{(n)} \cap b(x,R)} = 1 \text{ and } x \notin \bigcup_{j = 1}^J b^{(n)}(y_j,R)\}.
\end{equation*}
Recall that we assume that $\la R^d v_d < 1$. Let $k \in \{\lceil n-n\eps \rceil, \dots, n-1\}$ and denote $A_k(\mathbb{X}_1^{k}):=\bigcup_{j = 1}^k b^{(n)}(X_j,R)$ then a.s.
\begin{align*}
    \begin{split}
        &\E \Bigl[ \one \left[\mathbb{X}_{k+1}^n \in E_n(\mathbb{X}_1^{k})\right] \lvert \mathbb{X}_1^{k} \Bigr] \\
        &= \frac{1}{\abs{W_n}^{n-k}} \int_{W_n \setminus A_k(\mathbb{X}_1^{k})} \dots \int_{W_n\setminus A_k(\mathbb{X}_1^{k})} \one \left[ \abs{y_i - y_j} > R,\, \forall i\neq j\right] \d y_1  \dots \d y_{n-k}\\
        &\geq \dots \geq \frac{\left(\abs{W_n} - (n-1)v_d R^d \right)^{n-k}}{\abs{W_n}^{n-k}} = \left( 1 - \la \frac{n-1}{n} v_d R^d\right)^{n-k} \geq \left( 1 - \la v_d R^d\right)^{n-k}.
    \end{split}
\end{align*}
Using this, we can bound
\begin{align*}
        &\E e^{-H_n(\mathbb{X}_1^{n})}\one\left[K_n \in \left(n-n\eps,n\right) \right]\one \left[\onenJ{\bfXi}{n}(\mathbb{X}_1^{K_n}) \in \oneJ{O}\oplus b(0,\Delta)\right] \\
        & = \E \Bigl[e^{-H_n(\mathbb{X}_1^{K_n})}\one\left[K_n \in \left(n-n\eps,n\right) \right]\one \left[\onenJ{\bfXi}{n}(\mathbb{X}_1^{K_n}) \in \oneJ{O}\oplus b(0,\Delta)\right] \\& \hspace{8cm} \cdot \one \left[\mathbb{X}_{K_n+1}^n \in E_n(\mathbb{X}_1^{K_n})\right] \Bigr] \\
        & = \sum_{k = \lceil n-n\eps \rceil}^{n-1} \E \Bigl[ e^{-H_n(\mathbb{X}_1^{k})}\one\left[K_n = k\right]\one \left[\onenJ{\bfXi}{n}(\mathbb{X}_1^{k}) \in \oneJ{O}\oplus b(0,\Delta)\right] \\ & \hspace{8cm} \cdot \E \Bigl[ \one \left[\mathbb{X}_{k+1}^n \in E_n(\mathbb{X}_1^{k})\right] \lvert \mathbb{X}_1^{k}\Bigr] \Bigr] \\
        &\geq  \sum_{k = \lceil n-n\eps \rceil}^{n-1} \E \Bigl[ e^{-H_n(\mathbb{X}_1^{k})}\one\left[K_n = k\right]\one \left[\onenJ{\bfXi}{n}(\mathbb{X}_1^{k}) \in \oneJ{O}\oplus b(0,\Delta)\right] \left( 1 - \la v_d R^d\right)^{n-k} \Bigr] \\
        &\geq  \sum_{k = \lceil n-n\eps \rceil}^{n-1} \E \Bigl[ e^{-H_n(\mathbb{X}_1^{k})}\one\left[K_n = k\right]\one \left[\onenJ{\bfXi}{n}(\mathbb{X}_1^{k}) \in \oneJ{O}\oplus b(0,\Delta)\right] \left( 1 - \la v_d R^d\right)^{n\eps } \Bigr] \\
        &= \left( 1 - \la v_d R^d\right)^{n\eps }  \E \Bigl[ e^{-H_n(\mathbb{X}_1^{K_n})}\one\left[K_n \in \left( n-n\eps,n\right)\right]\one \left[\onenJ{\bfXi}{n}(\mathbb{X}_1^{K_n}) \in \oneJ{O}\oplus b(0,\Delta)\right]  \Bigr]. 
\end{align*}

Altogether, we get that (\ref{IneqI2}) holds. \enp

\section{Proof of Theorem \ref{lemma:difbc}} \label{sec:DifHam}

First, in Lemma \ref{lemma:hatLDP}, we prove that under some general conditions on the two Hamiltonians $H_n$ and $\hat{H}_n$, also the new Gibbs measure $\hat{P}_n$ satisfies the {unnormalized large-deviation bounds} with the same rate function as $P_n$. {This lemma is then used to prove Theorem \ref{lemma:difbc}.}

\bel[]\label{lemma:hatLDP}
Let $J \in \N$. Let $\oneJ{U} \subset \R^J$ be an open set and $\oneJ{C} \subset \R^J$ be a closed set. Assume that the Hamiltonian $H_n$ is given\footnote{Recall Definition \ref{def:BGPP}.} by an $r$-local and bounded (by a constant $c$) interaction function $V$. Assume that the score functions $\oneJ{\bfxi}$ are such that the following two inequalities hold
\begin{align}\label{as:lb}
\begin{split}
        \liminf { \frac{1}{\abs{W_n}} }&\log\int \exp(-H_n(\om)) \one \left[\onenJ{\bfXi}{n}(\om) \in \oneJ{U} \right]\d \mathbb{B}_n(\om)  \\ 
    &\geq - \inf \{I(P) + P^o(V): P \in \mc{P}_s, P^o({\bf1}) = \la, P^o(\oneJ{\bfxi}) \in \la \oneJ{U}\}, \\
\end{split}
\end{align}
\begin{align}\label{as:ub}
\begin{split}
        \limsup { \frac{1}{\abs{W_n}} }&\log\int \exp(-H_n(\om))  \one \left[\onenJ{\bfXi}{n}(\om)  \in \oneJ{C}\right]\d \mathbb{B}_n(\om)  \\ 
    &\leq - \inf \{I(P) + P^o(V): P \in \mc{P}_s, P^o({\bf1}) = \la, P^o(\oneJ{\bfxi})  \in \la \oneJ{C}\}.
    \end{split}
\end{align}
 Assume that for any $\eps \in (0,1)$ there exists $\{E_{n,\eps}\}_{n \in \N} \subset \mc{F}$ such that 
\begin{enumerate}
    \item[(i)] 
   $\sup_{\om \in E_{n,\eps}} \abs{H_n(\om)-\hat{H}_n(\om)} \leq c n \eps $ holds for any  $ n \in \N$,
    \item[(ii)]  $\limsup_{n \rightarrow \infty}{{ \frac{1}{\abs{W_n}} } \log \P (B_n \in E^c_{n,\eps})} = - \infty$,
    \item[(iii)] $\hat{H}_n(\om) \geq -c\abs{\om \cap W_n} $ holds for any $ \om \in E^c_{n,\eps}$ and $ n \in \N$.
    \end{enumerate} 
{ Then, 
\begin{align}\label{auxlb}
\begin{split}
    \liminf \frac{1}{\abs{W_n}}&\log \int \exp(-\hat{H}_n(\om)) \one \left[\onenJ{\bfXi}{n}(\om) \in \oneJ{U} \right]\d \mathbb{B}_n(\om)  \\ 
    &\geq - \inf \{I(P) + P^o(V): P \in \mc{P}_s, P^o({\bf1}) = \la, P^o(\oneJ{\bfxi}) \in \la \oneJ{U}\},
    \end{split}
\end{align}
\begin{align}\label{auxub}
\begin{split}
        \limsup \frac{1}{\abs{W_n}}&\log\int \exp(-\hat{H}_n(\om))  \one \left[\onenJ{\bfXi}{n}(\om)  \in \oneJ{C}\right]\d \mathbb{B}_n(\om)  \\ 
    &\leq - \inf \{I(P) + P^o(V): P \in \mc{P}_s, P^o({\bf1}) = \la, P^o(\oneJ{\bfxi})  \in \la \oneJ{C}\}.
    \end{split}
\end{align}}
\enl

\bep
We start with (\ref{auxlb}). Denote
$$
A(\oneJ{U}):= \inf \{I(P) + P^o(V): P \in \mc{P}_s, P^o({\bf1}) = \la, P^o(\oneJ{\bfxi}) \in \la \oneJ{U}\}.
$$
If $A(\oneJ{U}) = \infty$, then (\ref{auxlb}) holds trivially. Therefore, assume\footnote{Recall that $A(\oneJ{U})$ is bounded from below, since $I$ is non-negative and $V$ is bounded. } that $A(\oneJ{U}) < \infty$.

Let $\eps \in (0,1)$ and $n \in \N$, then 
\begin{align*}
        \E e^{-\hat{H}_n(B_n)} & \one \left[\onenJ{\bfXi}{n}(B_n) \in \oneJ{U} \right] \\
        &\geq \E e^{-H_n(B_n)} e^{-\abs{H_n(B_n)-\hat{H}_n(B_n)}} \one \left[\onenJ{\bfXi}{n}(B_n) \in \oneJ{U} \right] \one \left[B_n \in E_{n,\eps}\right] \\
        &\overset{(i)}{\geq} e^{-n \eps c} \E e^{-H_n(B_n)} \one \left[\onenJ{\bfXi}{n}(B_n) \in \oneJ{U} \right] \one \left[B_n \in E_{n,\eps}\right] \\
        &=  e^{-n \eps c} \E e^{-H_n(B_n)} \one \left[\onenJ{\bfXi}{n}(B_n) \in \oneJ{U} \right] \\ 
        & \qquad \qquad \qquad \qquad - e^{-n \eps c} \E e^{-H_n(B_n)} \one \left[\onenJ{\bfXi}{n}(B_n) \in \oneJ{U} \right] \one \left[B_n \in E^c_{n,\eps}\right]\\
        &\geq  e^{-n \eps c} \E e^{-H_n(B_n)} \one \left[\onenJ{\bfXi}{n}(B_n) \in \oneJ{U} \right]  - e^{n c} \P(B_n \in E^c_{n,\eps}).
\end{align*}
Thanks to (ii) we have that 
\begin{equation} \label{limsup}
    \limsup { \frac{1}{\abs{W_n}} }\log \left( e^{n c} \P(B_n \in E^c_{n,\eps})\right) = c \la +  \limsup { \frac{1}{\abs{W_n}} }\log \P(B_n \in E^c_{n,\eps}) = - \infty.
\end{equation}
On the other hand 
\begin{align}
\begin{split}\label{liminf}
    \liminf { \frac{1}{\abs{W_n}} }\log & \left(e^{-n \eps c} \E e^{-H_n(B_n)} \one \left[\onenJ{\bfXi}{n}(B_n) \in \oneJ{U} \right]\right) \\
    &= - c \eps \la +  \liminf { \frac{1}{\abs{W_n}} }\log \left(\E e^{-H_n(B_n)} \one \left[\onenJ{\bfXi}{n}(B_n) \in \oneJ{U} \right]\right) \\
    &\overset{(\ref{as:lb})}{\geq} -c \eps \la - A(\oneJ{U}) > - \infty. 
\end{split}
\end{align}
Therefore, 
\begin{align*}
    \begin{split}
         \liminf { \frac{1}{\abs{W_n}} }&\log\int \exp(-\hat{H}_n(\om)) \one \left[\onenJ{\bfXi}{n}(\om) \in \oneJ{U} \right]\d \mathbb{B}_n(\om) \\
         \geq \quad  &\liminf { \frac{1}{\abs{W_n}} }\log \left(e^{-n \eps c} \E e^{-H_n(B_n)} \one \left[\onenJ{\bfXi}{n}(B_n) \in \oneJ{U} \right]  - e^{n c} \P(B_n \in E^c_{n,\eps})\right)\\
         \overset{(\ref{limsup}),(\ref{liminf})}{=}  &\liminf { \frac{1}{\abs{W_n}} }\log \left(e^{-n \eps c} \E e^{-H_n(B_n)} \one \left[\onenJ{\bfXi}{n}(B_n) \in \oneJ{U} \right] \right) \\
          \overset{(\ref{as:lb})}{\geq} \quad &-c \eps \la - 
 \inf \{I(P) + P^o(V): P \in \mc{P}_s, P^o({\bf1}) = \la, P^o(\oneJ{\bfxi}) \in \la \oneJ{U}\}.
    \end{split}
\end{align*}
Taking $\eps \rightarrow 0$ finishes the proof of (\ref{auxlb}).

Now, we turn to the upper bound (\ref{auxub}). Take again $\eps \in (0,1)$ and $n \in \N$ and denote for simplicity $\one_1^J(B_n):= \one \left[\onenJ{\bfXi}{n}(B_n) \in \oneJ{C} \right]$, then 
\begin{align}
    \begin{split}\label{auxupbound}
        \E & e^{-\hat{H}_n(B_n)} \one_1^J(B_n)\\
        &=  \E e^{-\hat{H}_n(B_n)} \one_1^J(B_n)\one \left[B_n \in E_{n,\eps}\right] + \E e^{-\hat{H}_n(B_n)} \one_1^J(B_n)\one \left[B_n \in E^c_{n,\eps}\right] \\ 
        &\overset{(iii)}{\leq}\E e^{-H_n(B_n)} e^{\abs{H_n(B_n)-\hat{H}_n(B_n)}} \one_1^J(B_n) \one \left[B_n \in E_{n,\eps}\right]+ e^{nc} \P(B_n \in E^c_{n,\eps}) \\
        &\overset{(i)}{\leq} e^{n \eps c} \E e^{-H_n(B_n)} \one_1^J(B_n) \one \left[B_n \in E_{n,\eps}\right] + e^{nc} \P(B_n \in E^c_{n,\eps}) \\
        &\leq  e^{n \eps c} \E e^{-H_n(B_n)} \one_1^J(B_n) + e^{nc} \P(B_n \in E^c_{n,\eps}) .
    \end{split}
\end{align}
Clearly $\limsup { \frac{1}{\abs{W_n}} } \log e^{n \eps c} \E e^{-H_n(B_n)} \one_1^J(B_n) 
\overset{(\ref{limsup})}{\geq}\limsup { \frac{1}{\abs{W_n}} } \log e^{nc} \P(B_n \in E^c_{n,\eps}) $ and therefore 
\begin{align*}
    \limsup { \frac{1}{\abs{W_n}} } \log \Bigl(e^{n \eps c} \E e^{-H_n(B_n)} &\one_1^J(B_n) + e^{nc} \P(B_n \in E^c_{n,\eps}) \Bigr) \\ &= \limsup { \frac{1}{\abs{W_n}} } \log \left(e^{n \eps c} \E e^{-H_n(B_n)} \one_1^J(B_n)\right).
\end{align*}
Using this and (\ref{auxupbound}), we get that
\begin{align*}
\begin{split}
        \limsup { \frac{1}{\abs{W_n}} }\log\int &\exp(-\hat{H}_n(\om))  \one \left[\onenJ{\bfXi}{n}(\om)  \in \oneJ{C}\right]\d \mathbb{B}_n(\om)  \\
        &\leq c \eps \la + \limsup { \frac{1}{\abs{W_n}} } \log  \E e^{-H_n(B_n)} \one \left[\onenJ{\bfXi}{n}(B_n) \in \oneJ{C} \right]\\
    &\overset{(\ref{as:ub})}{\leq}  c \eps \la - \inf \{I(P) + P^o(V): P \in \mc{P}_s, P^o({\bf1}) = \la, P^o(\oneJ{\bfxi})  \in \la \oneJ{C}\}.
    \end{split}
\end{align*}
Taking $\eps \rightarrow 0$ finishes the proof of (\ref{auxub}).
\enp

\bere Before we proceed, we have two remarks regarding the assumptions of the previous lemma.
 \begin{enumerate}
     \item The constant in the second assumption for $\hat{H}_n$ could depend on $\eps$, however in such a way so that $c_\eps \eps \rightarrow 0$ as $\eps \rightarrow 0$.
     \item For the proof of the lower bound (\ref{auxlb}) we could change assumption (ii) so that we only require that the $\limsup$ is small enough. Also notice that we do not use assumption (iii). On the other hand for the upper bound (\ref{auxub}) we need (ii) as it is.
\end{enumerate}
\enre

Now, we will use Lemma \ref{lemma:hatLDP} to {prove Theorem \ref{lemma:difbc}. In the proof, we need the following notation. Let $r > 0$ and denote} the $r$-boundary of the window $W_n$ as 
$$
\partial_rW_n:= W_n \setminus \left[-\frac{1}{2}\left(\frac{n}{\la}\right)^{\frac{1}{d}}+r,\frac{1}{2}\left(\frac{n}{\la}\right)^{\frac{1}{d}}-r\right]^d.
$$

\bep[Proof of Theorem \ref{lemma:difbc}] 
{Since $V$ is bounded, Theorem \ref{thm:ldp} and Corollaries \ref{lemma:UBunbounded1} and \ref{cor:LBunboundedscore} tell us that the unnormalized bounds (\ref{as:lb}) and (\ref{as:ub}) hold for all three cases (a-c).} It is therefore enough to prove that $\hat{H}^1$ and $\hat{H}^2$ satisfy the assumptions (i)-(iii) from Lemma \ref{lemma:hatLDP} with 
\begin{equation*}
 E_{n,\eps}:= \{\om \in \Omega_{W_n}: \abs{\om \cap \partial_r W_n} \leq n \eps\},\, n \in \N,\, \eps \in (0,1).   
\end{equation*}
The rest of the proof then follows from (\ref{auxlb}) and (\ref{auxub}). 

As $V$ is bounded, clearly assumption (iii) holds for $\hat{H}_n^1$ with $c \geq c_1$ and for $\hat{H}_n^2$ with $c \geq 2c_1$. Regarding the assumption (i), from the assumptions on $V$ we get that
\begin{itemize}
    \item[(a)] $x \in W_n \setminus \partial_r W_n \implies V(\om^{(n)}-x) = V((\om\cap W_n)\cup(\gamma \cap W_n^c) - x)$,
    \item[(b)]  $x \in W_n \setminus \partial_r W_n \implies V((\gamma \cap W_n^c) - x)=0$.
\end{itemize}
Therefore, for $\om \in E_{n,\eps}$
\begin{align*}
    \abs{H_n(\om)-\hat{H}^1_n(\om)} &= \abs{\sum_{x \in \om\cap W_n} V(\om^{(n)} - x)- V(((\om\cap W_n)\cup(\gamma \cap W_n^c)) - x)} \\
    &\leq 2c \abs{\om \cap \partial_r W_n} \leq 2c n\eps
\end{align*}
and also
\begin{align*}
     \abs{H_n(\om)&-\hat{H}^2_n(\om)} \\
     &= \abs{\sum_{x \in \om\cap W_n} V(\om^{(n)} - x)- V(((\om\cap W_n)\cup(\gamma \cap W_n^c)) - x) {-} V((\gamma \cap W_n^c) - x)}\\
     &\leq 3c \abs{\om \cap \partial_r W_n} \leq 3c n\eps.
\end{align*}
To verify (ii), realize that $\abs{B_n \cap \partial_r W_n} \sim Bi(n,p_n)$ with parameter 
$$
p_n =\frac{\abs{\partial_r W_n}}{\abs{W_n}}= 1-\Bigl(1-2r\Bigl(\frac{\la}{n}\Bigr)^{\frac{1}{d}}\Bigr)^d.
$$
Using the concentration inequality for the binomial distribution (Lemma 1.1 in \cite{PenroseBook}) we get that for $n$ large enough (so that $\eps \geq e^2 p_n$)
\begin{equation*}
    \P(B_n \in E^c_{n,\eps}) = \P(\abs{B_n \cap \partial_r W_n} > n \eps ) \leq \exp\left(-\frac{n\eps}{2} \log \left(\frac{\eps}{p_n}\right)\right).
\end{equation*}
Therefore,
\begin{equation*}
    \limsup { \frac{1}{\abs{W_n}} }\log  \P(B_n \in E^c_{n,\eps}) \leq 
    - \frac{\eps \la \log \eps}{2} + \frac{\eps\la}{2} \lim\log\Bigl(1-\Bigl(1-2r\Bigl(\frac{\la}{n}\Bigr)^{\frac{1}{d}}\Bigr)^d\Bigr) = -\infty
\end{equation*}
and the proof is finished. 
\enp
\section*{Acknowledgments}
CH was supported by a research grant (VIL69126) from VILLUM FONDEN.
MP was supported by the Charles University Grant Agency (project no.\,70524), a visit grant by the Danish Data Science Academy (DDSA-V-2023-1701) and by the Mobility Fund of Charles University (project FM/c/2023-2-028).
\bibliographystyle{abbrv}
\bibliography{lit}

\end{document}